\documentclass{article}
\usepackage{graphicx} 
\usepackage{hyperref}   
\hypersetup{
	colorlinks=true,
	linkcolor=blue,
	citecolor=blue,
	urlcolor=blue,
}
\usepackage[square,numbers]{natbib}
\usepackage{xcolor}
\usepackage{amsthm}
\usepackage{subcaption}
\usepackage{dirtytalk}
\usepackage{algorithm}
\usepackage{algpseudocode}
\usepackage{nicematrix}
\usepackage{graphicx}
\usepackage{amsfonts}
\usepackage{amsmath,mathtools,amssymb,amsfonts}

\newcommand{\norm}[1]{\left\lVert#1\right\rVert}
\newcommand{\xddots}{%
  \raise 4pt \hbox {.}
  \mkern 6mu
  \raise 1pt \hbox {.}
  \mkern 6mu
  \raise -2pt \hbox {.}
}
\DeclareMathOperator*{\SubjectTo}{Subject\phantom{a}to:}
\DeclareMathOperator*{\Minimize}{Minimize:}
\DeclareMathOperator*{\Maximize}{Maximize:}
\DeclareMathOperator*{\argmin}{arg\,min}

\newtheorem{theorem}{\bf{Theorem}}

\newtheorem{assumption}{Assumption} 
\newtheorem{proposition}{\bf{Proposition}}
\newtheorem{remark}{\bf{Remark}}
\newtheorem{definition}{Definition}

\title{Degradation-based Energy Management for Microgrids in the Presence of Energy Storage Elements}
\author{Satish Vedula}

\date{}

\begin{document}

\maketitle
\begin{centering}
Department of Electrical and Computer Engineering, the Center for Advanced Power Systems, Florida State University
E-mail: \{svedula\}@fsu.edu
\end{centering}

\section{Abstract}
Integration of Inverter-based Resources (IBRs) such as solar-powered plants which lack the intrinsic characteristics such as the inertial response of the traditional synchronous-generator (SG) based sources presents a new challenge in the form of analyzing the grid stability under their presence. With the increasing penetration of the IBRs into the traditional grid, the non-renewable means of power generation which are dependent on coal, and natural gas is reducing drastically during the presence of power generated through renewable sources. For example, solar power is available for approximately from 9 AM-5 PM. However, the result of the rise in power consumption after 6 PM and the reverting back to the non-renewable source of power generation during that period puts immense stress on the grid, testing the ramp limitations of the SGs. Failure to meet the required power demand due to SG ramp limitations leads to failure of the power grid and other catastrophes. Numerous mitigation techniques exist in order to address the ramping issues with adding the energy storage elements (ESE) to the grid being one. ESEs have higher ramping capabilities compared to the traditional SGs. Also, the ESEs can store the energy and supply it to the grid when required making them extremely responsive to high ramp situations. However, the rate of degradation of the ESEs is faster than the SGs. Failure to adaptively manage the degradation of the ESEs might lead to the risk of sudden failure and eventually trigger a grid collapse. This raises an important issue of addressing the degradation of the ESEs while meeting the required power demand objectives and constraints. This work proposes a battery degradation-aware model predictive energy management strategy and it is tested via a numerical simulation on multiple physical systems such as Shipboard Power Systems (SPS) and Hybrid Electric vehicles (HEV). Moreover, the risk arising due to the fault in the IBR is also studied by means of a numerical simulation. Overall, the goal of this study is to make the existing power grid more robust, resilient, and risk-free from component degradation and eventual failures.  

\section{Introduction and Background}

\subsection{Motivation}
\textit{Feedback Control} is the technique that drives the mechanical, electrical, or biological systems towards equilibrium and maintains the system at equilibrium \cite{FrankLewis}. Origin of control dates back to ancient Egypt. Egyptians first designed a water clock around 1500 B.C Greeks improvised it around 300 B.C, and they called it float regulator. The main goal of the float regulator designed by the Greeks was to maintain a constant water level in the tank. Present-day ball in flush toilets works on a similar principle. Ancient Greeks also deployed this in various applications such as automatic wine dispensers.  Its principle is the modern-day equivalent of open-loop control. This concept of an open loop-controlled water clock soon moved across the world. The Chinese and the Korean water clocks were designed using similar principles. During $1200 A.D$, Arab engineers further improvised the existing float regulator and introduced an on/off control.

A major breakthrough in the field of controls was during the Industrial Revolution. Advanced mills, boilers, and steam engines were developed during this period. These advancements required new control methods, thus, temperature control, pressure control, etc were introduced. James Watt invented the steam engine in 1769 A.D. One prominent feedback control introduced by James Watt in 1788 A.D called fly-ball governor was used in steam engine speed regulation. The difference between the actual speed and the desired speed (error signal) dictated the amount of steam output to the engine \cite{CKang}. British astronomer G.B. Airy in $1840 A.D$ developed a closed loop feedback for a telescope that automatically adjusted the telescope concerning earth rotation. This gave birth to mathematical control theory. He was the first to discover that improper design of closed-loop systems will lead to oscillations and instability.

In 1868 A.D. James Maxwell published a paper titled On Governors \cite{JCMaxwell}. In this paper, Maxwell remarked that every physical system can be expressed as a differential equation. He analyzed fly-ball governors' stability by expressing it in differential equation form. Maxwell successfully presented the stability conditions for systems up to third order. Maxwell linearized the differential equations and formulated the \textit{characteristic equation} to determine the roots of the system. The stability condition proposed was the possible parts of all the roots shall be negative \cite{CKang}. He was unable to present the stability conditions for the fifth and higher order system completely and thus left it for future research enthusiasts. In 1877 A.D. English mathematician Edward J Routh and In 1895 German mathematician Adolf Hurwitz independently formulated a numerical method to determine the stable roots of the characteristic equation. This method became widely popular and is used even today in stability analysis under the name of \textit{Routh Hurwitz Criteria}. Even though many researchers during that time found Maxwell's paper extremely hard to interpret, it was Norbert Wiener who highlighted the significance of Maxwell's paper in 1948 A.D. In 1892 A.D. a Russian mathematician named A.M. Lyapunov studied the stability of nonlinear differential equations using energy functions. The work of Lyapunov was introduced to the Western world during the 1960s by a Hungarian-American mathematician named R.E. Kalman. Lyapunov's work is pivotal in the development of \textit{Modern Control Theory}.

While the control development was taking place during the late $19^{th}$ century, during the same time frame researchers like Thomas Edison and Nikola Tesla were making developmental strides in the development of the Alternating Current (AC) and Direct Current (DC). Edison and Tesla had differing views on electricity distribution, with Edison advocating for DC and Tesla for AC. Edison famously promoted the use of DC through public demonstrations. The rivalry between Edison's General Electric Company and George Westinghouse's company, which backed Tesla's AC system, intensified the competition between the two inventors. Tesla's work on AC power laid the groundwork for the modern electrical grid, which is based predominantly on AC distribution. The development of the first modern power grid is often credited to systems that emerged in the late 19$^{th}$ and early 20$^{th}$ centuries. Here are a few notable examples:
\begin{itemize}
    \item \textbf{Pearl Street Station (1882)}: This was the first centralized power plant in the United States, located in lower Manhattan, New York City. Thomas Edison's Edison Illuminating Company built it and began operations on September 4, 1882. The Pearl Street Station provided direct current (DC) electricity to customers within a one-mile radius using underground conductors.
    \item \textbf{London Electric Supply Corporation (1882)}: In the UK, the London Electric Supply Corporation (LESC) established one of the earliest central power stations in the world. It was located in Deptford, London, and began supplying electricity to street lamps and a few businesses in 1882. The system used both AC and DC, with AC for long-distance transmission and DC for local distribution.
\end{itemize}

\begin{figure}[h!]
      \centering
      \includegraphics[width=0.95\textwidth]{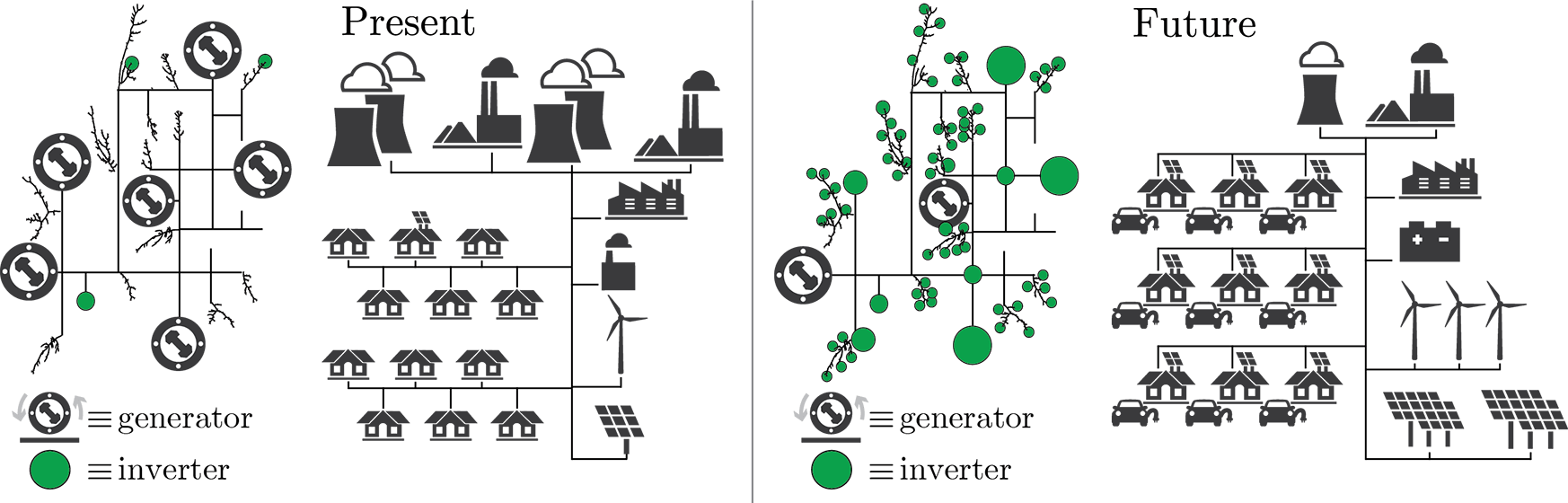}
	 \caption{Present and future trends in the IBR penetration in the microgrids \cite{9729134}.}
     \label{IBR_Present_Future} 
\end{figure}

The late 20$^{th}$ century saw the rapid growth of the solid-state electronic industry which led foundation for the modern power inverters. Inverters were primarily used in specialized applications such as uninterruptible power supplies (UPS), telecommunications, and industrial motor drives. The 1980s saw the commercialization of solar photovoltaic (PV) technology and the growth of wind power as a viable renewable energy sources. Inverters became essential components in solar PV and wind turbine systems, converting DC power generated by solar panels and wind turbines into AC power suitable for grid integration. Early inverters were relatively simple and lacked many of the advanced features found in modern inverters. The 2010s saw explosive growth in solar PV and wind power installations worldwide, fueled by declining costs, supportive policies, and growing environmental awareness. Inverter-based resources became mainstream, with utility-scale solar and wind projects incorporating advanced inverters capable of grid support functions such as voltage regulation and reactive power control. Distributed energy resources (DERs), including rooftop solar systems and community solar installations, proliferated, driving demand for residential and commercial-scale inverters. Figure \ref{IBR_Present_Future} shows the penetration of the IBRs into the modern grid.

\begin{figure}[h!]
      \centering
      \includegraphics[width=0.85\textwidth]{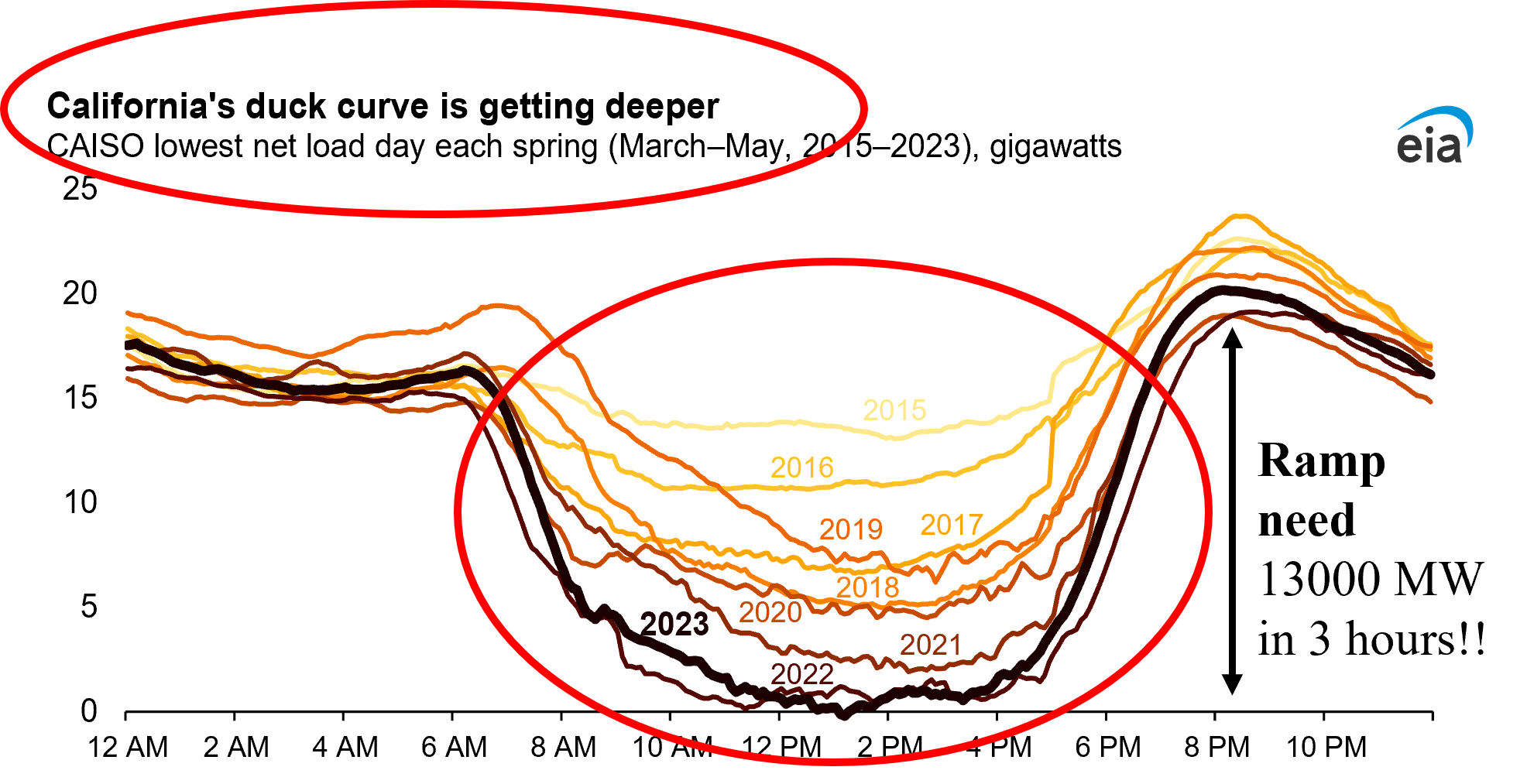}
	 \caption{California's power generation through non-renewable means from years 2015-2022 \cite{duck_curve}.}
     \label{duck_curve} 
\end{figure}

The ``California Duck Curve" shown in Figure \ref{duck_curve} refers to a graphical representation of the electricity demand pattern in California throughout the day. The curve gets its name because the shape resembles a duck, with a steep drop in demand during midday followed by a sharp increase in the evening hours. This pattern emerged as a result of the increasing adoption of solar power in California. The steep drop during midday is due to the abundance of solar power generated from rooftop solar panels and large-scale solar farms during sunny hours. As the sun rises, solar generation ramps up, leading to a decrease in demand for electricity from other sources. However, as the sun sets and solar generation declines, there is a rapid increase in demand as people return home from work, and switch on lights, appliances, and air conditioning, hence the ``belly" of the duck curve. The California Independent System Operator (CAISO) manages the grid in California and closely monitors this phenomenon. While the duck curve poses challenges for grid operators in managing the balance between supply and demand, it also presents opportunities for integrating energy storage, demand response programs, and other flexible resources to help smooth out the curve and ensure grid stability. Additionally, efforts are underway to optimize the use of renewable energy, manage energy demand more efficiently, and invest in grid infrastructure to accommodate the changing dynamics of energy generation and consumption.

Flattening the California Duck Curve requires a multi-faceted approach involving various strategies and technologies. Here are some solutions that can help flatten the curve:
\begin{itemize}
    \item \textbf{Energy Storage}: Deploying energy storage systems, such as batteries, pumped hydro storage, or thermal storage, allows excess energy generated during periods of high solar output to be stored and used during times of high demand, thereby smoothing out fluctuations in supply and demand.
    \item \textbf{Demand Response Programs}: Implementing demand response programs incentivizes consumers to adjust their electricity usage in response to grid conditions. This can involve shifting electricity consumption to times of low demand through time-of-use pricing, smart thermostats, or automated control systems.
    \item \textbf{Distributed Energy Resources (DERs)}: Encouraging the adoption of distributed energy resources, such as rooftop solar panels, behind-the-meter batteries, and electric vehicle (EV) chargers, can decentralize energy generation and consumption, reducing strain on the grid during peak periods.
    \end{itemize}

    By implementing these solutions in a coordinated manner, California and other regions facing similar challenges can work towards flattening the Duck Curve shown in Figure \ref{duck_curve}, ensuring grid stability, and maximizing the benefits of renewable energy integration.

    Adding energy storage to the grid offers several advantages, including:
\begin{itemize}
    \item \textbf{Grid Stability and Reliability}: Energy storage systems provide grid operators with the flexibility to balance supply and demand in real-time, helping to stabilize the grid and ensure reliable electricity delivery. They can store excess energy during periods of low demand and release it during peak demand, reducing strain on the grid and minimizing the risk of blackouts or grid instability.

    \item \textbf{Integration of Renewable Energy}: Energy storage helps to address the intermittency and variability of renewable energy sources, such as solar and wind power. By storing surplus renewable energy when it's abundant and releasing it when needed, energy storage enables greater integration of renewables into the grid, reducing reliance on fossil fuels and lowering greenhouse gas emissions.

    \item \textbf{Peak Shaving and Load Management}: Energy storage systems can be used to "shave" peaks in electricity demand by discharging stored energy during times of high demand, thereby reducing the need for expensive peaking power plants and alleviating strain on the grid infrastructure. Additionally, they enable more efficient management of energy loads, optimizing electricity usage and reducing overall energy costs.

    \item \textbf{Grid Ancillary Services}: Energy storage can provide valuable ancillary services to the grid, such as frequency regulation, voltage support, and reactive power control. These services help to maintain grid stability, improve power quality, and enhance the overall reliability of the electricity system.

    \item \textbf{Backup Power and Resilience}: Energy storage systems can serve as backup power sources during grid outages or emergencies, providing critical electricity supply to homes, businesses, and critical infrastructure. This enhances grid resilience and helps to minimize the impact of disruptions on communities and economies.

    \item \textbf{Electric Vehicle Integration}: Energy storage can facilitate the integration of electric vehicles (EVs) into the grid by providing charging infrastructure and managing EV charging patterns. By storing renewable energy for EV charging during off-peak hours and supporting bi-directional charging capabilities, energy storage maximizes the environmental and economic benefits of electric transportation.
\end{itemize}

\begin{figure}
    \centering
    \includegraphics[width=\linewidth]{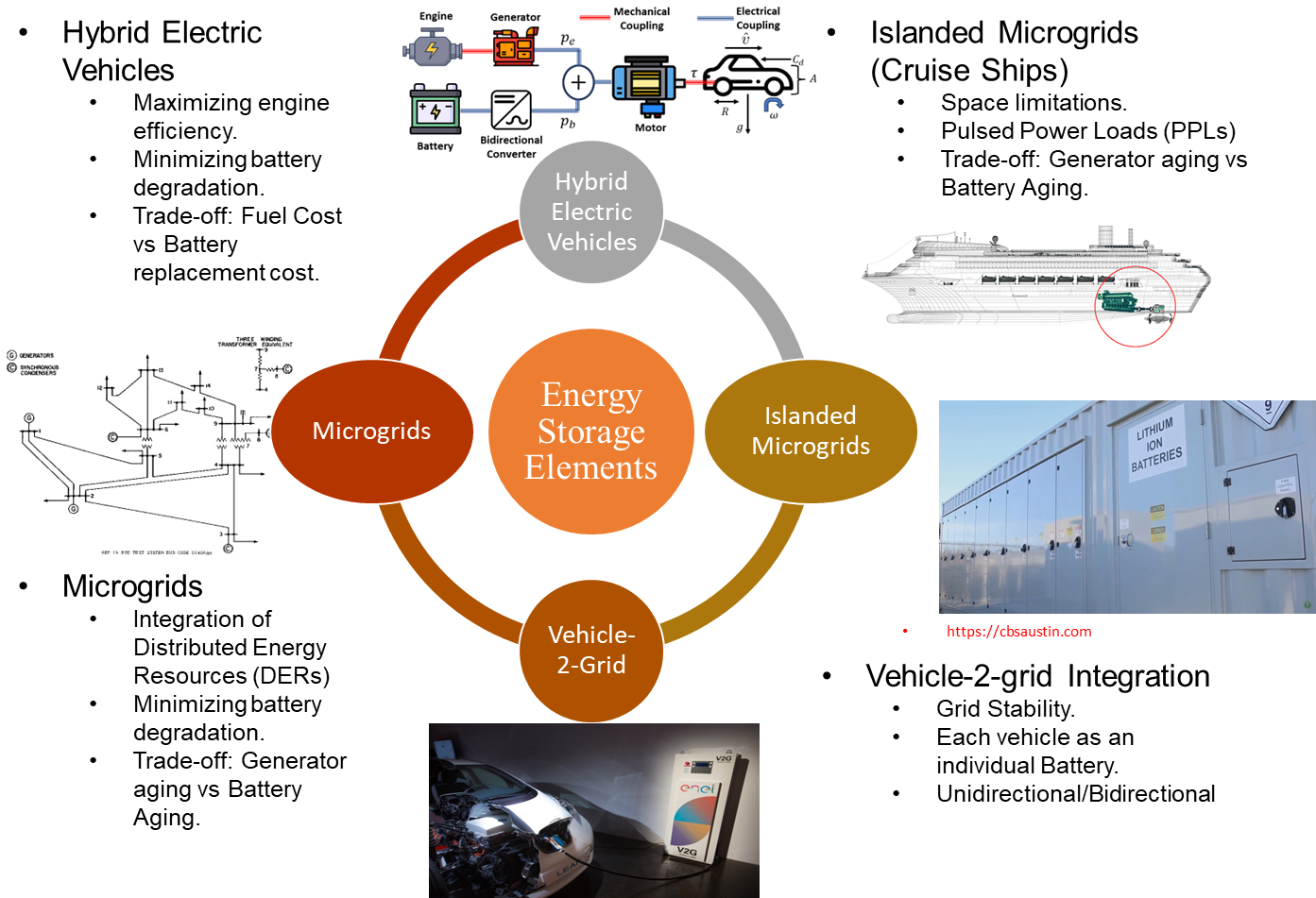}
    \caption{Energy storage elements presence in today's grid.}
    \label{Battery_Overview}
\end{figure}

Overall, adding energy storage to the grid offers a wide range of benefits, including enhanced grid flexibility, renewable energy integration, improved reliability, and cost savings, making it a valuable tool for transitioning towards a more sustainable and resilient energy future. Overall, it can be seen from Figure \ref{Battery_Overview} that the presence of the energy storage elements in the grid raises an important concern in terms of managing their degradation since generators and batteries have different lifespans. Generators, especially those powered by fossil fuels such as diesel or gasoline, can have a long lifespan if properly maintained. Industrial-grade generators can last anywhere from 10,000 to 30,000 hours of operation or more. Whereas, the lifespan of a battery depends on its chemistry, usage patterns, depth of discharge (DoD), temperature, and maintenance and is usually around 5-8 years.

\subsection{Background}
\subsubsection{What is Energy Dispatch Process?}
\begin{figure}[ht]
	\centering
	\includegraphics[width=0.85\textwidth]{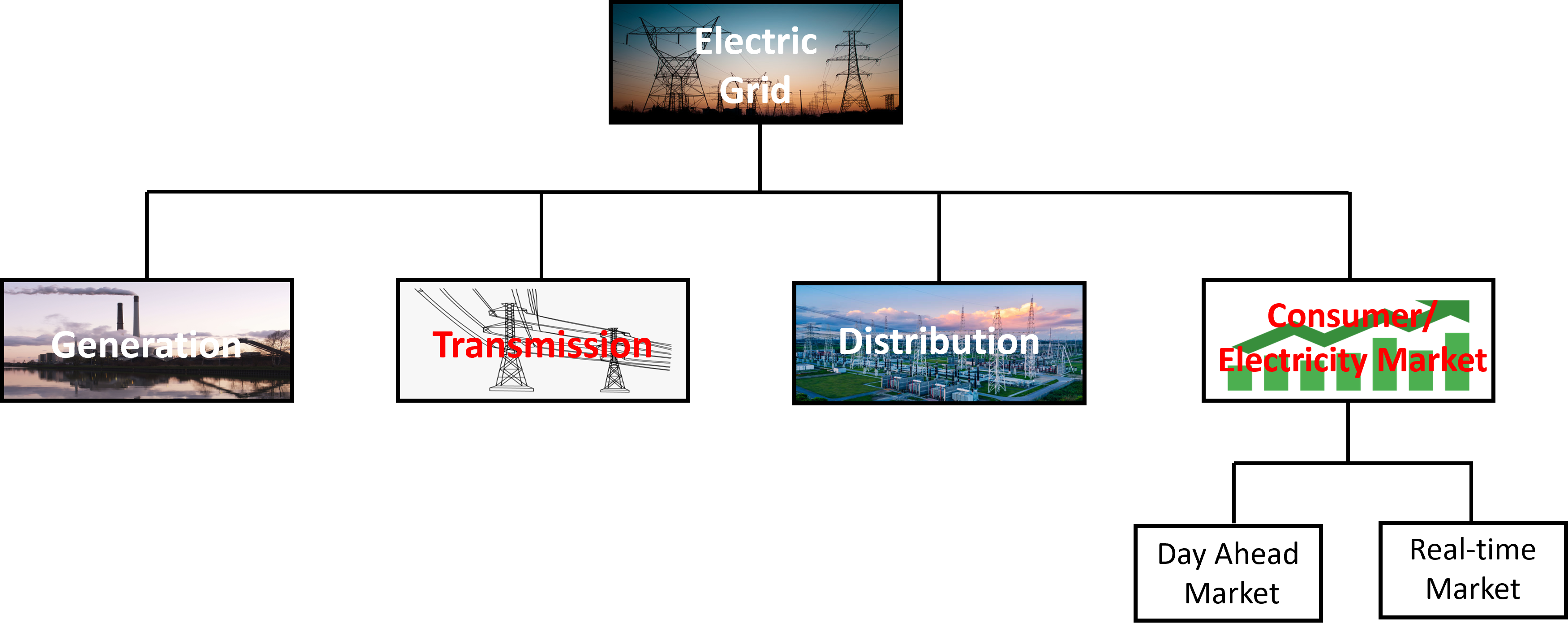}
	\caption{Electric grid subsystems}
 \label{Electric grid}
\end{figure}

\textit{Modern electric grid} is a complex system that is mainly divided into four subsystems namely: \textit{generation}, \textit{transmission}, \textit{distribution} and \textit{consumer market} as shown in Figure \ref{Electric grid}. In recent years, power generation has been intensively shifting towards green and renewable energy sources. Previously during many instances, much of the generated energy was left unused. However, introducing energy storage elements such as batteries, super-capacitors, etc leads to the storage of unused power and use when necessary. Thus creating a more economical grid operation. 

\subsubsection{Power Flow Problem}
The power-flow model must oblige to power generated equals power demanded criteria. Thus, the power-flow model can be represented as:    
\begin{align} \label{Power Flow Model}
f(p_g,p_L)=0,
\end{align}
 where $p_g$ and $p_L$ represent the generator sets (gen-sets) and the load powers. In (\ref{Power Flow Model}), generators are considered to inject power in a unidirectional form. Another representation which is usually helpful as an optimization constraint is:
\begin{equation}
    \sum_{i=1}^{n_g}p_{g_i} = p_L,
\end{equation}
where $n_g$ is the number of generators. 
\subsubsection{Optimal Energy Dispatch}
\begin{figure}[h!]
	\centering
	\includegraphics[width=0.8\textwidth]{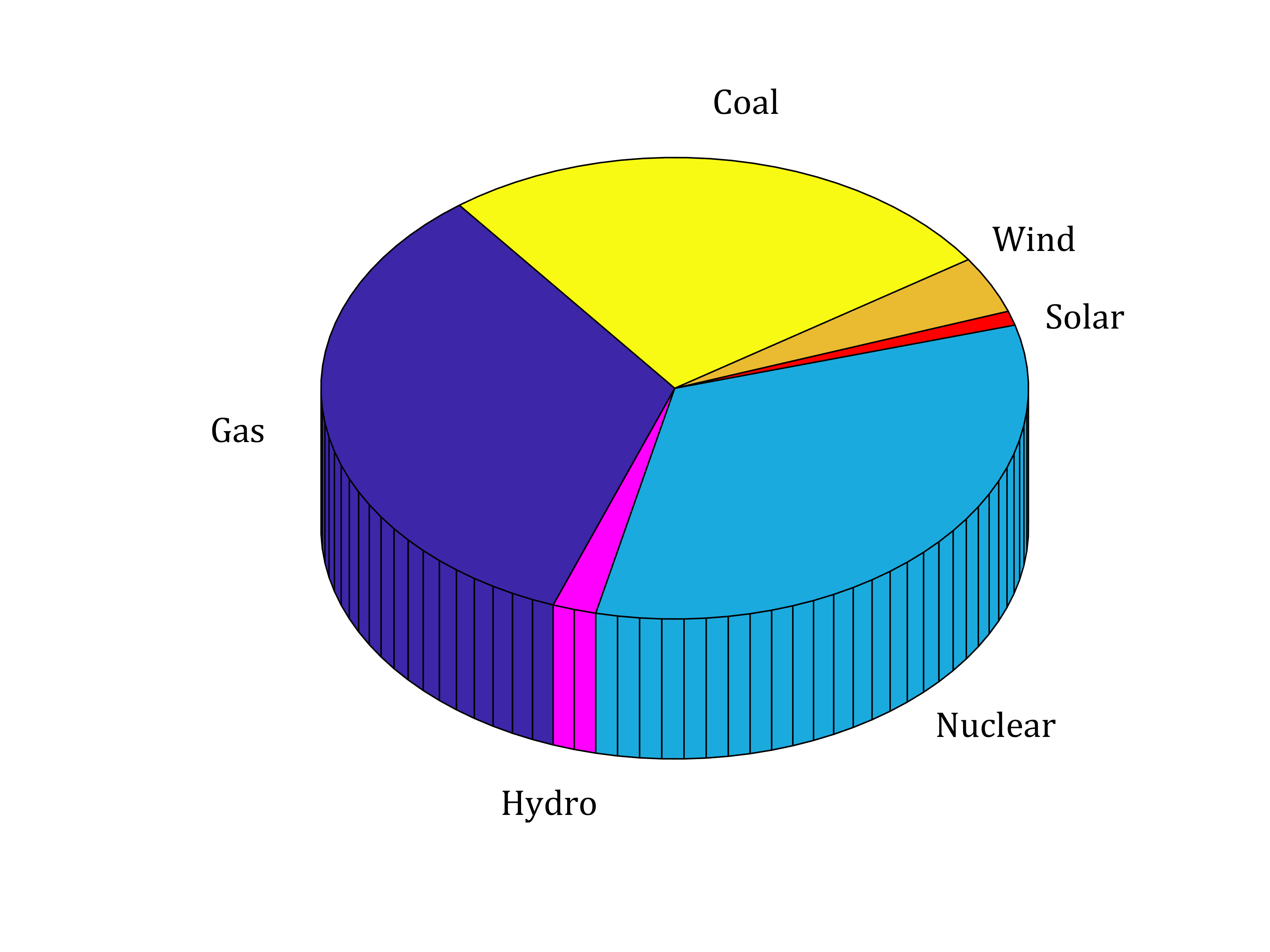}
	\caption{Pie chart showing American power generating sources and their contributions}
    \label{Energy distribution}
\end{figure}
In a power grid, the power generation units have varying costs of generation based on their location from loads. In the United States up to 75$\%$ of the power generation is done through Gas, coal, and Nuclear powers shown in Figure \ref{Energy distribution}. Usually, the capacity for power generation is more than the power demand for normal operation. The main motive is to minimize the fuel cost of operation of each and individual generating unit and also make sure that the power flow problem constraint is satisfied. This is called a \textit{optimal power flow} problem. Much of the discussion in this section is based on \cite{principles}. The power flow problem here is posed as an optimization problem that minimizes the specific objective in accordance with imposed constraints, which could be generator power limitations ramp rate limitations etc. Our main focus here is to study about \textit{optimal economic energy dispatch}. The \textit{cost function} for generator operation is usually modeled as a quadratic cost. One of the main reasons is that quadratic functions are convex in nature. Thus the solution is feasible and convergence is guaranteed. The following is the representation of the cost function:
\begin{equation}
    C_i(p_i) = a_ip_i^2+b_ip_i+c_i.
\end{equation}
The \textit{incremental fuel cost} represents the cost of production of power at the next increment. It is represented as follows:
\begin{equation}
    \frac{dC_i}{dp_i} = 2a_ip_i+b_i.
\end{equation}
If we neglect the losses and only consider the generator power limitation. The optimization problem that minimizes the cost is given as follows:
\begin{equation}\label{Economic Dispatch}
\begin{aligned}
\Minimize  \quad & \sum_{i}^{n_g}  C_i(p_{g_i})\\
\SubjectTo \quad & \sum_{i=1}^{n_g}p_{g_i} = p_L,\\
& \underline{p}_{g_i} \preceq p_{g_i} \preceq \overline{p}_{g_i},
\end{aligned}
\end{equation}
where $\underline{p}_{g_i} \in \mathbb{N}$ and $\overline{p}_{g_i} \in \mathbb{N}$ represent lower and upper limitations on generators. $n_g$ represents the number of generators. Since generators are assumed to operate unidirectionally, the lower limit is usually $\underline{p}_{g_i} \geq 0$. In some cases, individual generators might reach their operational capacity, which means they supply power equivalent to $\overline{p}_{g_i}$. In such cases that specific generator is set to supply maximum power and the process is repeated with the remaining generators.


\subsection{Assessment of Energy Dispatch Process}
In the field of power systems, \textit{Operating Reserve} is the generational capacity available to meet the power demand in case of disruption of supply at any given time interval. Usually, power systems are designed to meet minimum supplier load requirements and fraction of peak load. Operating reserves consist of two methods namely: \textit{spinning reserve} and \textit{non-spinning reserve} \cite{operatingreserves}. Spinning reverse is meeting the power deficit left by the failed generator by increasing the power output from remaining generators in the grid. In some cases, due to generator power limitations, the power deficit left because of generator failure cannot be met by increasing the output from remaining generators operating in the same grid. In such situations, the deficit power is supplied from the neighbouring grid. This method is called non-spinning reserve. Operating reserves play crucial role in ensuring day ahead approach withstands the unforeseen uncertainties.

\textit{Thermal Ratings} of the transmission are also be considered as a risk factor in propagation of power in a power system. Thermal limit is determined by the maximum power the transmission material or the conductor can carry without burning the material. The temperature of these transmission line must be kept below a rated maximum temperature for safe operation. By abiding to it ageing of the conductors can be delayed. Usually transmission line thermal ratings can be adjusted by adjusting the power loading through that line. There are two types of thermal line rating methods namely: static line rating, which accommodates for seasonal weather changes and dynamic line rating, which changes with real-time weather conditions. In \cite{DADouglass} Douglass provides insight into dynamic thermal line rating methods with forecasting. This thermal requirement can be maintained in real-time by direct measurement of conductor temperature and sagging from on line data.

\subsection{Assessment of Current Electricity Market}
\textit{Day Ahead Market} and \textit{Real-time Market} fall under \textit{Locational Marginal Pricing} (LMP) methods. The concept of LMP was introduced in 1980 by Caramanis, Bohn, and Schweppe at the MIT energy laboratory \cite{Caramanis}. The pricing is done based on the principle that the price of one unit of power varies at different points in a given network \cite{JBastian}. During this period various international governments were privatizing utility companies. This is also known as \textit{deregulation}. To have a competitive edge over their contemporaries, the utility companies started pricing the electricity based on \textit{spot pricing of electricity}. This introduction of pricing separated transmission and distribution from the electricity pricing market, creating a whole new branch in power distribution based on electricity pricing. LMP is a pricing method used by numerous existing utility companies or \textit{independent system operators} (ISOs). One of the famous ISOs is the New York ISO. These ISOs were created by the Federal Energy Regulatory Commission (FERC) as a set towards privatization. Two of the LMP methods that are hugely popular are day-ahead market pricing and real-time market pricing. The idea behind the day ahead market is to predict the load requirements for future time and adjust the electricity pricing accordingly. \textit{Forecasting} the LMPs on a long or short-term basis is termed as day ahead approach. In \cite{JBastian} authors gave a brief insight into forecasting energy prices in a competitive market. Usually, the forecast is done using statistical methods such as time series analysis or methods based on simulations. With increased computational capabilities much of the simulation methods are based on optimization. The main goal of the optimization problem is to minimize the electricity cost subject to some probabilistic constraints. \textit{Power Exchange} (PX) is a common platform established worldwide by most of the local country governments for managing the electricity power bids for buying and selling \cite{2020_Devnathshah}. In simple words, PX is an internet-based online platform that acts as a unifying place for buyers and sellers.

\textit{Real-time Energy Market} allows electricity suppliers and consumers to sell buy and operate during the day in real-time. The electricity demand bid can be met within one hour of requesting it. One of the main advantages of a real-time market is that operating reserves can be eliminated from the energy dispatch process. But on the other hand, the price per MW of electricity skyrockets during peak operating hours. This is the reason alternatives such as adding ESS to the grid are gaining popularity. By adding ESS to the grid we can operate on a day-ahead basis and the excess energy can be stored or sold off to another power retailer. The prominent advantage of using a real-time market is the risk is mitigated. Operating reserves always pose an imminent threat of mishandling the excess power. In the real-time market, the power demand is met in real-time by spinning reserves, thus widely minimizing the risk of excess power wastage. Cyber attacks on the power grid can be minimized by operating in real-time. This is studied in detail in the next section.

\textit{Unit Commitment} (UC) problem is the basis for the current operation of the day-ahead market. The idea behind it is to decide which generating units connect over the upcoming 24 hours. It involves integer decision unit committed (1) or not committed (0). Most of the present ISOs run security-constrained unit commitment (SCUC) problems. \textit{Security constrained economic dispatch} (SCED) is analogous to SCUC. SCED is a mathematical model that enforces laid constraints while delivering power economically. SCED can be used either in real-time (5 minutes ahead) or day-ahead markets (24 hours ahead). Qin in \cite{2013_Qin} presented a real-time electricity market in which risk is used to model the system's overall security level. The term \textit{Risk Based SCED} (RB-SCED) was coined by Qin. The key difference between the risk associated with the power grid compared to financial institutions lies in the source of uncertainty. According to IEEE standards \textit{risk can be calculated as the product of the probability of a contingency occurrence multiplied by the consequence of that contingency} \cite{2013_Qin}.

\textit{Contingency Analysis} \cite{VJMishra} in power systems is associated with planning and utilizing operating reserves. SCED is one such method as discussed before. Power systems are operated in a fashion that does not overload them in real time or at any given instance. This is referred to as maintaining system security. North American Electric Reliability Corporation (NREC) presides over the reliability issues in the power grid arising due to supply failures. The grid operational standards must meet those set by NREC to ensure stability. Power system stability is defined as \textit{the ability of an electric power system, for a given initial operating condition, to regain a state of operating equilibrium after being subjected to a physical disturbance, with most system variables bounded so that practically the entire system remains intact} \cite{PKundur}. Power system stability is classified based on three key power system states, they are: \textit{Rotor Angle Stability}, \textit{Frequency Stability}, and \textit{Voltage Stability}. While the rotor angle stability is considered to be a short-term effect, frequency, and voltage instabilities might have either short-term or long-term consequences.

Traditional power systems must at least abide by the $N-1$ contingency placed by NREC. This means the grid operations must not be disrupted by a single generator outage. With ever growing complexity of modern power systems, NREC raised those standards to $N-k$ where $k \geq 2$. But by increasing $k$ in this constraint, there is a chance that huge reserves of power generated might go unused. Given the stochastic nature of power demand such a constraint is not desirable for power generation. In \cite{FBouffard} introduced the loss of load probability (LOLP) criterion to address the issue of reliability-constrained market clearing. According to Bouffard LOLP is \textit{the probability that the available generation, including spinning reserve, cannot meet the system load}\cite{FBouffard}. In simple words, when demand exceeds the generation, a loss of load situation arises. Loss of load probability specifies the occurrence of such a situation in a year. LOLP is usually considered a long-term resource planning method. In \cite{2011_Varaiya} varaiya used the LOLP formulation and developed an algorithm for a risk-limiting dispatch for smart grid operation. Sj$\ddot{o}$din in \cite{2012_Sjodin} used the algorithm developed by Varaiya in developing an optimal power flow algorithm for wind-powered systems also considering the energy storage elements and fast ramp-supporting backup generators. In this work, Sj$\ddot{o}$din analyzes the risk-mitigated optimal power flow problem in a new fashion by adding storage charge or discharge dynamics to the optimal power flow problem. Adding these dynamics a finite horizon optimal control problem has been formulated.

\textit{Demand Response} is a change in power consumption by the customer to match supply accordingly. Usually, this is considered as \textit{post-contingency} step. According to FERC demand response is defined as: \textit{Changes in electric usage by end-use customers from their normal consumption patterns in response to changes in the price of electricity over time, or to incentive payments designed to induce lower electricity use at times of high wholesale market prices or when system reliability is jeopardized} \cite{2011_Murthy}. The main goal of demand response is to engage the community and consumers in managing their load usage based on pricing changes. Developments in energy storage technology disabled impromptu dispatch of power. The concept of demand response is similar to \textit{dynamic demand} mechanism. Smart grid technology enables energy communication with neighboring houses, which is extremely useful in managing power during peak hours. The difference between demand response and dynamic response is that the dynamic response mechanism senses irregularities in the grid disconnects that portion automatically and reroutes the power supply. Demand response is part of energy efficient homes which include EV charging and EM management for building.

\subsection{Assessment of Cyber Threats}
\begin{figure}[ht]
	\centering
	\includegraphics[width=0.7\textwidth]{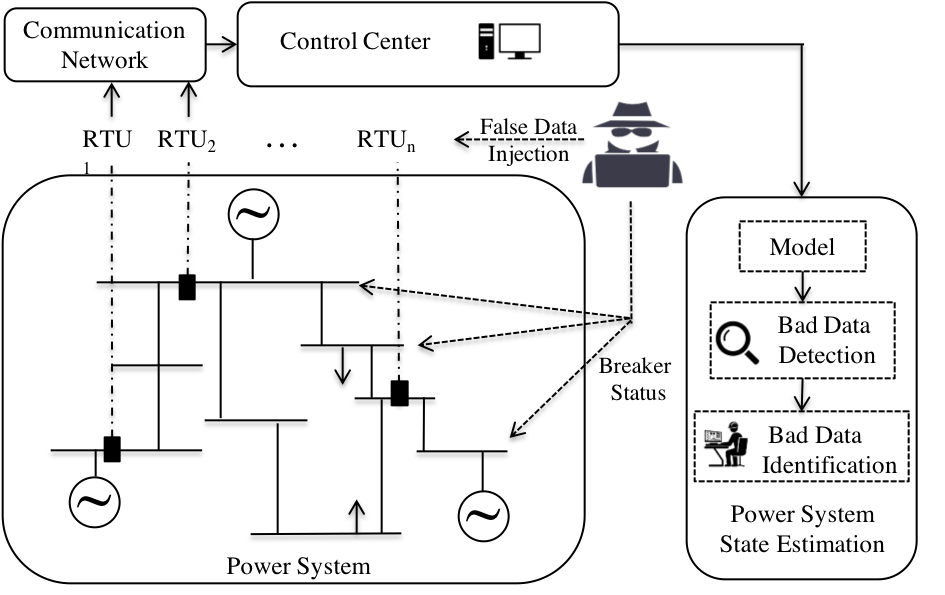}
	\caption{An illustration of CPS attacks}
    \label{CPS_Attack}
\end{figure}

\textit{Cyber security} is a critical topic and it needs to be addressed with crucial importance in \textit{Cyber Physical Systems} (CPS), especially for microgrids. It also has to be addressed from an energy dispatch perspective. With rapidly expanding smart grids more integrated communication networks are being added to the grid raising serious cyber security concerns. Most of the communication networks in microgrids have open communication channels and are prone to cyber-attacks. Cyber-attacks hinder seamless grid operation which is considered a risk.  One of the prominent methods widely used by adversaries to attack a cyber-physical system is through injecting false data. This attack mode is called \textit{False Data Injection Attack} (FDIA) shown in Figure \ref{CPS_Attack}. FDIAs tend to compromise multiple power grid readings such as sensor and phasor measurements, which include bus voltages, and real and reactive power injections at that particular bus. Adversaries may inject FDIAs between communication lines between measurements and \textit{SCADA} systems. They might also impact economic dispatch leading to a surge in operational costs and severe system outages. This might mislead control and operation centers, where the data is retained. If the adversary has complete knowledge about power grid elements and their positioning, false data can be cleverly crafted in such a way that it goes undetected by residue-based \textit{Bad Data Detectors} (BDD). In most of the realistic cases, the adversary does not have access to complete real-time knowledge of grid parameters.

Most of the existing literature on the impact of FDIAs on the power grid is under the assumption that the adversary has complete grid component knowledge. Liu in \cite{2009_Yao} studied state estimation in electric power grids under FDIAs with the assumption of complete grid knowledge. This paper presents that cleverly crafted FDIA can bypass BDD. Liu considers attacks against state estimation in DC power flow models. Even though DC power flow models are less accurate compared to AC power flow models, they are linear and more robust \cite{2006_VanHertem}. Both the attacker's perspective and the operator's perspective were considered. From the attacker's perspective both the realistic cases \textit{targeted FDIA} and \textit{arbitrary FDIA} were considered. Rahman in \cite{2012_Rahman} built on Liu's results and discusses the impact of FDIAs under incomplete information about grid parameters. The main contributions of this paper are: considering realistic FDIA attack with incomplete information, that is, the adversary does not have real-time knowledge of grid parameters. A mathematical characterization of FDIA was also provided discussing the probability of attack detection and its impact on mitigating the risk of compromising the grid. The novelty of this work compared to Liu's work lies in the introduction of \textit{vulnerability measure} which compares numerous power grid topologies. This helps in building attack-resilient power grids. Direct current-based microgrids (DCMG) have been on the rise in recent times. Still, the majority of distribution systems and grids rely on alternating current (AC). One of the main benefits of DCMG is its high efficiency and low maintenance costs. Another benefit is smooth coupling between the main grid and DCGM \cite{2020_Cheng}. Modern ship power systems, electric charging stations, etc are examples of DCMG. \textit{Consensus} based distributed control algorithm uses \textit{graph theory} approach in solving the economic dispatch problem. Graph theory acts as a basis for multi-agent consensus problems.

\subsubsection{Load Altering Attacks}
A Load Altering Attack (LAA) is a cyber-physical attack against the demand response and the demand-side management programs \cite{7131791,5976424,Lakshminarayana2022LoadAlteringAA,9838515,9393479,7723861}. Unsecured \textit{controllable loads} are altered in order to cause adverse effects in the grid. LAAs are classified into two types of attacks: Static Load Altering Attacks (SLAA) \cite{5976424} and Dynamic Load Altering Attacks (DLAA) \cite{7131791,10197626,7436350}. SLAA implies static changes in the load. This is attributed to shifting the percent of the load requirements over the actual value. DLAA are further classified into multiple LAAs namely: Open-Loop DLAA, Single-point DLAA, and Multi-point DLAAs \cite{7723861}.

\section{Notations and Preliminaries}
\subsection{Notations and Definitions}
$\mathbb{N}$, $\mathbb{R}$, and $\mathbb{R}_+$ denote the set of natural, real, and positive real numbers. $\mathcal{L}_2$ and $\mathcal{L}_{\infty}$ denote the square-integrable (measurable) and bounded signal spaces. A real matrix with $n$ rows and $m$ columns is denoted as $X \in \mathbb{R}^{n \times m}$. The identity matrix of size $n$ is denoted as $I_n$. $\lambda_{min}(X), \lambda_{max}(X)$ denote the minimum and the maximum eigenvalues of a matrix $X$. $X^\top$ denotes the transpose of a matrix $X$. \textsf{trace}$(X) \triangleq \sum_{i=1}^{n}x_{ii}$ denotes the trace of a matrix $X$. \textsf{det}$(X)$ denotes the determinant of a matrix $X$. Natural and Real scalars are denoted by lowercase alphabets (for example $x \in \mathbb{N}$ and $y \in \mathbb{R}$). The real vectors are represented by the lowercase bold alphabets (i.e. $\textbf{x} \in \mathbb{R}^{n}$). The vector of ones and zeros is denoted as $\mathbf{1}$ and $\underline{\mathbf{0}}$. 
 For any vector $\mathbf{x} \in \mathbb{R}^n$, $\|\mathbf{x}\|_2 \triangleq \sqrt{\mathbf{x}^\top\mathbf{x}}$, $\|\mathbf{x}\|_1
 \triangleq \sum_{i=1}^{n}|\mathbf{x}_i|$, and $\|\mathbf{x}\|_{\infty} \triangleq \max(|x_1|,|x_2|,\ldots,|x_n|)$ representing the 2-norm, 1-norm, and the $\infty$-norm respectively (where $|.|$ denotes absolute value). The symbol $\preceq$ denotes the component-wise inequality i.e. $\mathbf{x} \preceq \mathbf{y}$ is equivalent to $\mathbf{x}_i \leq \mathbf{y}_i$ for $i=1,2,\hdots,n$. The dot product/inner product of the two vectors $\mathbf{x} \in \mathbb{R}^n$ and $\mathbf{y} \in \mathbb{R}^n$ is denoted as $\mathbf{x}^\top \mathbf{y}$. The cross product/outer product is denoted as $\mathbf{x}^\top J \mathbf{y}$, where $J \triangleq \begin{bmatrix}
     0 & 1 \\ -1 & 0
 \end{bmatrix}$. For a function $f:\mathbb{R}^n \longrightarrow \mathbb{R}^m$, $\nabla f(x)$ denotes the gradient of the function $f$ at $x$. $\nabla^2f(x)$ denotes the hessian of the function $f$ at $x$. The Kronecker product of two matrices $X \in \mathbb{R}^{n \times n}$, $Y \in \mathbb{R}^{m \times m}$ is denoted as  $$X \otimes Y \triangleq\begin{bmatrix}
    X_{11}Y&X_{12}Y\\
    X_{21}Y&X_{22}Y\\
    \end{bmatrix} \in \mathbb{R}^{nm \times nm}.$$
    $\mathcal{B}(0,r)\triangleq\left\{\mathbf{x}\in\mathbb{R}^n|\left\|\mathbf{x}\right\|_{2} <r\right\}\subset \mathbb{R}^n$ denotes a ball of radius $r$ centered at the origin.

\subsubsection{Graph Theory}  The material in this subsection can be found in more detail in \cite{FB-LNS}.
    Consider the network with $m \in \mathbb{N}$ nodes or agents labeled by the set $\mathcal{V}=\{1,2,\hdots,m\}$, $\mathcal{E}$ is the \emph{unordered} edge such that $\mathcal{E} \subseteq \mathcal{V} \times \mathcal{V}$. The connection between the nodes is fixed $\mathcal{G}=(\mathcal{V},\mathcal{E})$. For an undirected graph, The \emph{adjacency} matrix ($A \in \mathbb{R}^{m \times m}$) is denoted as 
    $a_{ij}= \{1 \hspace{1mm}  \text{if} \hspace{1mm} (i,j)\in\mathcal{E} | 0 \hspace{1mm}  \text{otherwise}\},$ 
    i.e. if there is a path between two nodes, then the value of 1 is assigned otherwise the value in the matrix is set to 0. The \emph{degree} of a graph ($\mathcal{N}_i$) denotes the number of neighboring nodes of a given node. The degree matrix ($D \in \mathbb{R}^{m \times m}$) is the degree of a given node on the diagonal and 0's elsewhere. The \emph{Laplacian} matrix ($L=D-A$), is the difference between the degree and the adjacency matrix. The Laplacian-based weighted graph matrix is defined as follows: $W = I-\frac{1}{\tau}L$, where $\tau > \frac{1}{2}\lambda_{max}L$ is a constant, $\lambda_{max}$ is the maximum eigenvalue of $L$. 

\subsubsection{Vector Norms and Inequalities}

For any $p \geq 1$, $\norm{\mathbf{x}}_p \triangleq \bigg(\sum_{i=1}^{n}|\mathbf{x}_i|^p\bigg)^{\frac{1}{p}}$ denotes the $p$-norm of $\mathbf{x} \in \mathbb{R}^n$. For some $p > r \geq 1$ a vector $\mathbf{x} \in \mathbb{R}^n$, the topological equivalence of the norms is given as follows
\begin{align*}
    \norm{\mathbf{x}}_p \leq \norm{\mathbf{x}}_r \leq n^{(\frac{1}{r}-\frac{1}{p})}\norm{\mathbf{x}}_p.
\end{align*}
For example, the inequality mentioned above translates to
\begin{align*}
    \norm{\mathbf{x}}_{\infty} \leq \norm{\mathbf{x}}_2 \leq \norm{\mathbf{x}}_1 \leq \sqrt{n}\norm{\mathbf{x}}_2 \leq n \norm{\mathbf{x}}_{\infty}.
\end{align*}
For some $p \geq 1$, and $\mathbf{x}, \mathbf{y} \in \mathbb{R}^n$ the \textit{Triangle Inequality} is given as
\begin{align*}
    \norm{\mathbf{x} + \mathbf{y}}_p \leq \norm{\mathbf{x}}_p + \norm{\mathbf{y}}_p.
\end{align*}

\subsubsection{Signal Norms}
For a signal $\mathbf{x}(t)$, if $\norm{\mathbf{x}(t)}_p$ exists,
\begin{align*}
    \norm{\mathbf{x}(t)}_p \triangleq \bigg(\int_{0}^{\infty}|\mathbf{x}(\tau)|^p d\tau\bigg)^{\frac{1}{p}}, \hspace{2mm} p \in [1,\infty),
\end{align*}
we say that $\mathbf{x}(t) \in \mathcal{L}_p$. We say that $\mathbf{x}(t) \in \mathcal{L}_2$ (space of square integrable functions) if $\norm{\mathbf{x}(t)}_2$ exists, where
\begin{align*}
    \norm{\mathbf{x}(t)}_2 \triangleq \sqrt{\int_{0}^{\infty}|\mathbf{x}(\tau)|^2 d\tau}.
\end{align*}
The $\mathcal{L}_{\infty}$-norm (space of bounded functions) of a signal $\mathbf{x}(t)$, $\norm{\mathbf{x}(t)}_{\infty} \triangleq \sup_{t \geq 0} |\mathbf{x}(t)|$, then $\mathbf{x}(t) \in \mathcal{L}_{\infty}$ if $\norm{\mathbf{x}(t)}_{\infty}$ exists.

If $a \geq 0$, and $b \geq 0$, and for some $p, q > 1$ such that $\frac{1}{p} + \frac{1}{q} = 1$, then the following inequality holds (it is famously known as \textit{Young's Inequality} for products
\begin{align*}
    ab \leq \frac{a^p}{p} + \frac{b^q}{q}.
\end{align*}

\subsection{Optimization Review}
The material in this subsection is based on the book ``Convex Optimization" \cite{boyd2004convex} and \cite{ryu_yin_2022}.
For a function $f:X \to Y$, the domain of $f$ is denoted by $dom(f)$. A function $f:\mathbb{R}^n \to \mathbb{R}$ is \textit{convex} if $dom(f)$ is a convex set and $f(\gamma x+(1-\gamma) y) \leq \gamma f(x)+(1-\gamma)f(y)$ $\forall$ $x,y\in dom(f)$ and $0 \leq \gamma \leq 1$. The function $f$ is \textit{strictly convex} if the strict inequality holds.
 \begin{definition}
        A function $f:\mathbb{R}^n \longrightarrow \mathbb{R}^m$ is said to be $L$-Lipschitz if:
        $$\norm{f(x)-f(y)} \leq L \norm{x-y} \hspace{2mm} \forall x,y \in \mathbb{R}^n,L \in (0,\infty).$$
    \end{definition}
    \begin{definition}
        A closed, convex, and proper (CCP) function $f$ is $\mu$-strongly convex if any of the following equivalent conditions are satisfied:
        \begin{enumerate}
        \item $\langle \nabla f(x)-\nabla f(y), x-y \rangle \geq \mu \norm{x-y}^2 \forall x,y$
       \item $\nabla^2f(x) \succcurlyeq \mu I$ $\forall x$ if $f$ is twice continuously differentiable.
        \end{enumerate}\end{definition}
    \begin{definition}
        A closed, convex, and proper (CCP) function $f$ is $L$-smooth if any of the following equivalent conditions are satisfied:
        \begin{enumerate}
        \item $\langle \nabla f(x)-\nabla f(y),x-y \rangle \geq \frac{1}{L}\norm{\nabla f(x)-\nabla f(y)}^2 \forall x,y$ if $f$ is differentiable.
        \item $\nabla^2 f(x) \preceq LI$ $\forall x$ if $f$ is twice continuously differentiable. 
    \end{enumerate} \end{definition}

    \begin{definition}
        For a function $f$ which is $\mu$-strongly convex and $L$-smooth, the following inequality holds:
        $$\mu \norm{x-y}^2 \leq \norm{\nabla f(x)-\nabla f(y)}\norm{x-y} \leq L\norm{x-y}^2.$$
        $f$ is \textit{convex} $\Longleftrightarrow$ $\nabla^2f\succcurlyeq0$ $\forall$ $x\in dom(f)$.
    \end{definition}

\subsubsection{Duality}
Consider the optimization problem of the form
\begin{equation}\label{opt_problem}
    \begin{aligned}
        \Minimize_\mathbf{x} \quad & f(\mathbf{x})\\
        \SubjectTo \quad & g(\mathbf{x}) \leq 0 \\
        & h(\mathbf{x}) = 0, 
    \end{aligned}
\end{equation}
where $\mathbf{x} \in \mathbb{R}^n$ are the decision variables, $f:\mathbb{R}^n \rightarrow \mathbb{R}_+$ is the function capturing optimization objective, $g:\mathbb{R}^n \rightarrow \mathbb{R}^m$ are the inequality functions and $h:\mathbb{R}^n \rightarrow \mathbb{R}^p$ are the equality constraints. The \textit{Lagrangian} of the optimization problem in (\ref{opt_problem}) is
\begin{equation}
    \mathcal{L}(\mathbf{x},\boldsymbol{\lambda},\boldsymbol{\nu})=f(\mathbf{x})+\boldsymbol{\lambda}^\top g(\mathbf{x})+\boldsymbol{\nu}^\top h(\mathbf{x}),
\end{equation}
where $\boldsymbol{\lambda} \in \mathbb{R}^m$ and $\boldsymbol{\nu} \in \mathbb{R}^p$ are the \textit{dual variables} associated with the inequality and the equality constraints. \textit{Lagrange dual function} is defined as
$$d(\boldsymbol{\lambda},\boldsymbol{\nu}) \triangleq \inf_\mathbf{x} \mathcal{L}(\mathbf{x},\boldsymbol{\lambda},\boldsymbol{\nu}).$$

Thus, given a primal optimal $\mathbf{x}^*$, the dual problem $d(\boldsymbol{\lambda},\boldsymbol{\nu}) = \mathcal{L}(\mathbf{x}^*,\boldsymbol{\lambda},\boldsymbol{\nu})$ and explicitly expressed as
\begin{equation}
\begin{aligned}
    \Maximize \quad & d(\boldsymbol{\lambda},\boldsymbol{\nu}) \\
    \SubjectTo \quad &  \boldsymbol{\lambda} \succcurlyeq 0
\end{aligned}
\end{equation}

\subsubsection{Distributed Optimization}
Over the last two decades, there has been an increased interest in the application of distributed optimization to various systems. This has been actively deployed in multi-agent systems in signal processing and also in distributed control of robotic networks \cite{DistCtrlRobotNetw}. The objective of distributed optimization is to optimize the global objective or cost function through local nodes and nodes exchange information with neighboring nodes through the underlying aggregator network. The main motivation behind the emergence of this optimization method is to handle big data and large-scale networks more efficiently. The general form of the problem is as follows:
\begin{equation}
\begin{aligned}
    \Minimize \quad & \sum_{i}^{n}h_i(x_i) \\ 
    \SubjectTo  \quad & x_i \in \mathcal{C},
\end{aligned}
\end{equation}
where $h_i: \mathbb{R} \longrightarrow \mathbb{R}_+$ is a function representing local objective of node $i$, $\mathcal{C} \subseteq \mathbb{R}^\text{n}$ is convex closed set. A review of various existing distribution optimization techniques is given in \cite{Angelia}. In recent years a new method known as \textit{Alternating Direction Method of Multipliers} or \textit{ADMM} is gaining momentum. This was first introduced in the 1970's by Gabay, Mercier, and others. Further developments were made throughout the 1980s and 1990s.    

\subsubsection{Alternating Direction Method of Multipliers}
\textit{ADMM} \cite{boyd_admm} works on the basis of decomposed coordination, solutions at each node are coordinated to find the solution to the overall global problem. This algorithm blends dual decomposition and augmented Lagrangian method for constrained optimization. This algorithm is used to solve the following forms of optimization problems:
\begin{equation}\label{ADMM}
\begin{aligned}
    \Minimize \quad & h(\mathbf{x})+g(\mathbf{z}), \\
    \SubjectTo \quad & Q\mathbf{x}+P\mathbf{z} = \mathbf{k},
\end{aligned}
\end{equation}
where $\mathbf{x} \in \mathbb{R}^n$, $\mathbf{z} \in \mathbb{R}^m$ are variables and $Q \in \mathbb{R}^{p \times n}$, $P \in \mathbb{R}^{p \times m}$ and $\mathbf{k} \in \mathbb{R}^p$. The functions $h(\mathbf{x})$ and $g(\mathbf{z})$ are assumed to be convex. The optimal value of the problem (\ref{ADMM}) is represented as: $p^* = \inf\bigg\{h(\mathbf{x})+g(\mathbf{z})\bigg|Q\mathbf{x}+P\mathbf{z} = \mathbf{k}\bigg\}$. The augmented Lagrangian is given as follows:
\begin{equation}\label{Augmented Lagrangian}
    \mathcal{L}_{\rho}(\mathbf{x},\mathbf{z}, \boldsymbol{\lambda}) = h(\mathbf{x})+g(\mathbf{z})+\boldsymbol{\lambda}^\top(Q\mathbf{x}+P\mathbf{z}-\mathbf{k})+\frac{\rho}{2} \norm{Q\mathbf{x}+P\mathbf{z}-\mathbf{k}}_2^2,
\end{equation}
where $\boldsymbol{\lambda} \in \mathbb{R}^p$ is the associated Lagrange multiplier or dual variable and $\rho > 0$ is \textit{penalty parameter}. The variables and the dual are updated in the following sequential manner:
\begin{subequations}\label{ADMM Iterations}
\begin{align}
    \mathbf{x}^{t+1} &= \argmin_{\mathbf{x}}\mathcal{L}_{\rho}(\mathbf{x},\mathbf{z}^{t},\boldsymbol{\lambda}^{t}), \\
    \mathbf{z}^{t+1} &= \argmin_{\mathbf{z}}\mathcal{L}_{\rho}(\mathbf{x}^{t+1},\mathbf{z},\boldsymbol{\lambda}^{t}), \\
    \boldsymbol{\lambda}^{t+1} &= \boldsymbol{\lambda}^{t}+\rho \bigg({Q\mathbf{x}}^{t+1}+P\mathbf{z}^{t+1}-\mathbf{k}\bigg),
\end{align}
\end{subequations}
where the variables $\mathbf{x}$ and $\mathbf{z}$ are minimized in (\ref{ADMM Iterations}a) and (\ref{ADMM Iterations}b). The dual variable is updated in (\ref{ADMM Iterations}c). Since in this method, the variables are updated in \textit{alternating} fashion, the name of the method is \textit{alternating direction method of multipliers}.

\textit{Scaled Form} is a convenient method of representing ADMM. The residual is defined as $\mathbf{r} = Q\mathbf{x}+P\mathbf{z}-\mathbf{k} $. Combining the linear and quadratic components of (\ref{Augmented Lagrangian}) we arrive at: $$\boldsymbol{\lambda}^\top\mathbf{r}+\frac{\rho}{2}\norm{\mathbf{r}}_2^2 = \frac{\rho}{2} \norm{\mathbf{r+u}}_2^2-\frac{\rho}{2} \norm{\mathbf{u}}^2,$$ Where $\mathbf{u} = \frac{1}{\rho}\boldsymbol{\lambda}$ is scaled dual variable. Using scaled dual form the iterations are represented as:
\begin{subequations}\label{Scaled ADMM Iteration}
\begin{align}
    \mathbf{x}^{t+1} &= \argmin_{\mathbf{x}}\bigg(h(\mathbf{x})+\frac{\rho}{2} \norm{Q\mathbf{x}+P\mathbf{z}^t-\mathbf{k}+\mathbf{u}^t}_2^2\bigg), \\
    \mathbf{z}^{t+1} &= \argmin_{\mathbf{z}}\bigg({g(\mathbf{z})}+\frac{\rho}{2} \norm{{Q\mathbf{x}}^{t+1}+P\mathbf{z}-\mathbf{k}+\mathbf{u}^t}_2^2\bigg), \\
    \mathbf{u}^{t+1} &= \mathbf{u}^{t}+Q\mathbf{x}^{t+1}+P\mathbf{z}^{t+1}-\mathbf{k}.
\end{align}
\end{subequations}
The forms presented in (\ref{ADMM Iterations}) and (\ref{Scaled ADMM Iteration}) are equivalent.

\subsection{Model Predictive Control Review}
\begin{figure}[h!]
	\centering
	\includegraphics[width=0.7\textwidth]{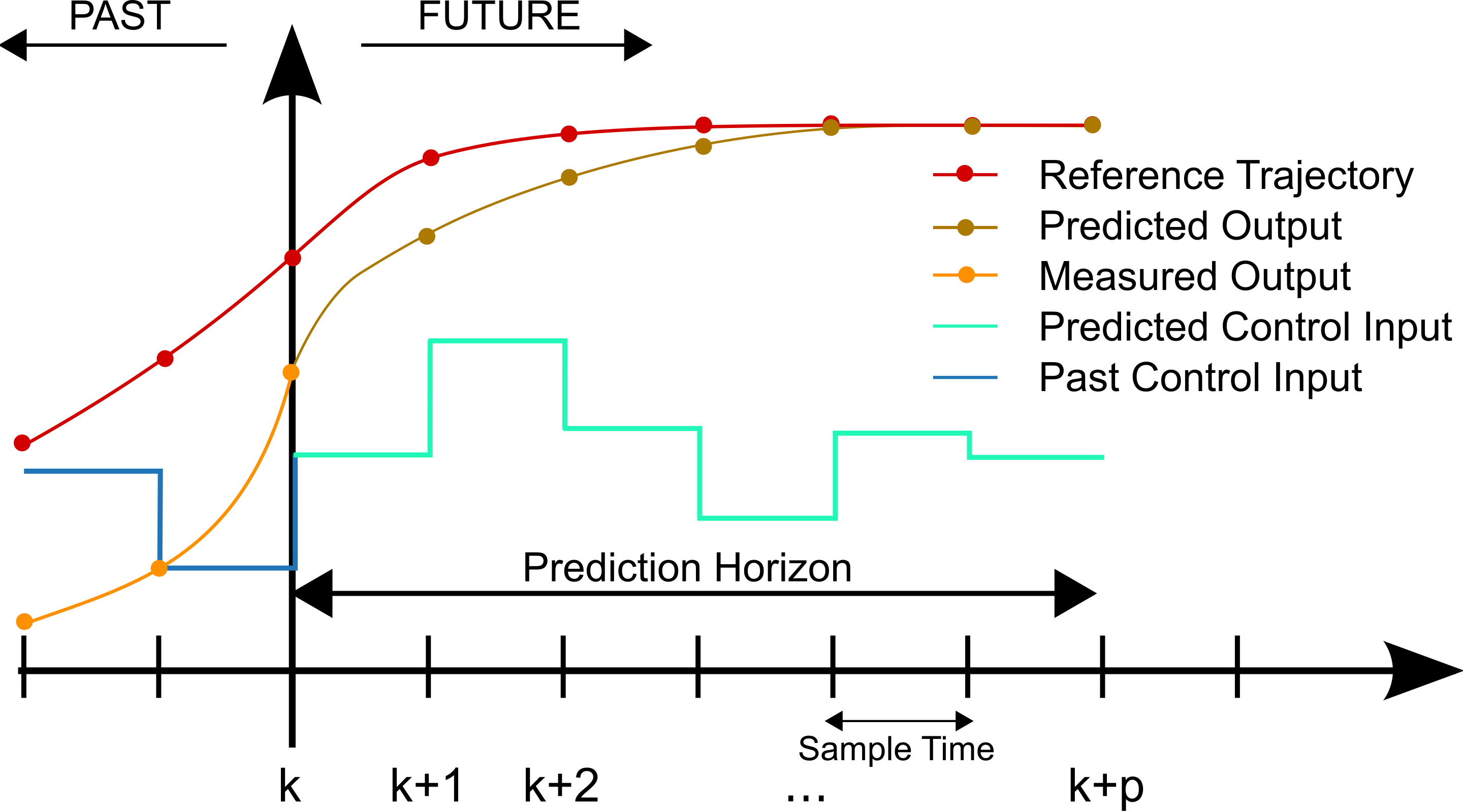}
	\caption{MPC implementation \cite{martin}.}
 \label{MPC_format}
\end{figure}

Model Predictive Control (MPC) is an optimal control technique that minimizes a cost function adhering to a set of constraints. It first originated in the field of chemical process control and slowly made its way into other fields such as power systems. MPC is based on a finite horizon optimal control. The optimization problem is solved at the current time step based on the state measurements keeping in account the future time steps honoring the set of imposed constraints on the states and the control inputs as shown in Figure \ref{MPC_format}. More detailed information about MPC can be found in \cite{2009_Rawlings}.

Consider a linear time-invariant (LTI) system
\begin{equation}\label{LTI}
\begin{aligned}
    \dot{\mathbf{x}}(t) &= A \mathbf{x}(t)+B\mathbf{u}(t), \\
    \mathbf{y}(t) &= C \mathbf{x}(t),
\end{aligned}
\end{equation}
where $\mathbf{x}(t) \in \mathbb{R}^n$ are the system states, $\mathbf{u}(t) \in \mathbb{R}^m$ denotes the system control input and $\mathbf{y}(t) \in \mathbb{R}^m$ denote the system measurable outputs. Consider a control objective where the goal is to design a control input such that the output track the desired output adhering to a set of constraints on the states and the control inputs. The MPC design is given as
\begin{equation}
    \begin{aligned}
        \Minimize \quad & \sum_{k=1}^{N}\norm{C\mathbf{x}_k-\mathbf{y}_d^k}_2^2+\sum_{k=0}^{N-1}\norm{u_k}_2^2 \\
        \SubjectTo \quad & \mathbf{x}_0 = \mathbf{x}_{meas}, \\
        & \mathbf{x}_{k+1} = A\mathbf{x}_k+B\mathbf{u}_k, \hspace{2mm} k=0,1,\hdots,N-1, \\
        &\underline{\mathbf{x}} \preceq \mathbf{x}_k \preceq \overline{\mathbf{x}}, \\
        &\underline{\mathbf{u}} \preceq \mathbf{u}_k \preceq \overline{\mathbf{u}}, \\
        &\left|\mathbf{u}_{k+1}-\mathbf{u}_k\right| \preceq r_u\mathbf{1},
    \end{aligned}
\end{equation}
where the initial value to the optimization problem $\mathbf{x}_0$ is the system state measurement $\mathbf{x}_{meas}$. $\mathbf{y}_d^k \in \mathbb{R}^m, \hspace{2mm} k=1,2,\hdots,N$ is the desired trajectory. The discretized linear dynamics are imposed as an equality constraint. $N$ is the horizon counter. $r_u$ is the rate limit of the rate of change of control input.


\section{Literature Review}

\subsection{Microgrids}
\textit{Micro grids} (MGs) were first introduced in the early 2000's by Lasseter \cite{Lasseter_2001}. MGs have been attracting attention during recent times due to their ability to provide reliable integration of \textit{distributed generation units} (DGUs) \cite{Asanso_2007}. Currently, they form a key component for distributed and decentralized power systems. They are being viewed as reliable sources of power supply to remote locations. For \textit{Smart grids}, microgrids form an important building block since all the distributed energy resources which comprise DGUs and distributed storage units (DSUs), can be utilized in a coordinated and efficient fashion. \cite{Vasilakis}. As presented in \cite{CAbbey} MGs have made significant strides towards power system resiliency and are more resilient during natural disasters like storms etc. The paper discusses the case studies of MGs during tsunamis in Japan and earthquakes in New Zealand and Chile. The main distinction of MGs from traditional distribution lines is their ability to operate in an islanded fashion. Switching from interconnected to islanded mode prompts significant imbalances, which need to be accounted for. This opens up numerous challenges in their control and operation.  In \cite{Hatziagyriou} various MG energy management frameworks have been presented. While most of the early literature proposed a hierarchical control mechanism for power sharing in microgrids, much of the focus in recent years has been diverted towards droop control \cite{Kroposki}. 

Droop control is a control mechanism used in both AC and DC systems \cite{2011_Gue}. In Existing droop control schemes for three-phase systems, frequency, and voltage are controlled for optimal power sharing. Frequency variation of even the smallest level might not be desirable in power distribution within the grid as that might lead to grid instability. Adding some power electronic components might relax these strict frequency and voltage regulations. Usually, inverters play a key role in such relaxations. The active power regulation is based on frequency control and the reactive power regulation is based on voltage control. Usually, the droop slope depends on the load demand. The DQ approach uses a DC-like control approach. The effects of DQ droop settings based on energy storage settings for inverter-based microgrids have been studied in \cite{MMBijaieh}.

\subsection{Shipboard Power Systems}
Islanded Medium Voltage DC (MVDC) Microgrids (MGs) deployed on Ship Power Systems (SPSs) provide power to various on-board equipment such as propulsion motors, hotel loads, and highly non-linear loads such as Pulsed Power Loads (PPLs). There is a significant challenge to utilize existing ramp-rate limited generators to provide balanced power to high ramp-rate loads. The situation is exacerbated by the integration of today's various power electronics equipment. This can result in an unbalanced system with degraded power quality, instability, and subsequent load shedding and reconfiguration. 

A major challenge for the operation of power electronics equipment is their Constant Power Load (CPL) behavior where the current is inversely proportional to the voltage. In this case, a negative incremental impedance can lead to instability  \cite{1976_MB,2009_Weaver}. However, fast utility load fluctuations cause more problems than the existence of CPLs \cite{2015_Cupelli}. Since the ramp-rate support of generators is limited, ESSs with high ramp-rate support capabilities are offered as a solution \cite{2019_Bijaieh}. Like DC MG systems, a hierarchical control architecture is currently utilized for power and energy management of SPSs \cite{2016_Jin}. In this case, a Power Management System (PMS) is designed to solve the underlying control allocation problem, and the Energy Management System (EMS) is utilized for resource allocation. Conventionally, the PMS exists to ensure the stability and performance of the system while addressing the generation power over-actuation, and EMS exists to feed appropriate PMS set-points to achieve various high-level objectives. This opens the ground for various optimization approaches for EMSs. 

\begin{figure}[h!]
	\centering
	\includegraphics[width=0.95\textwidth]{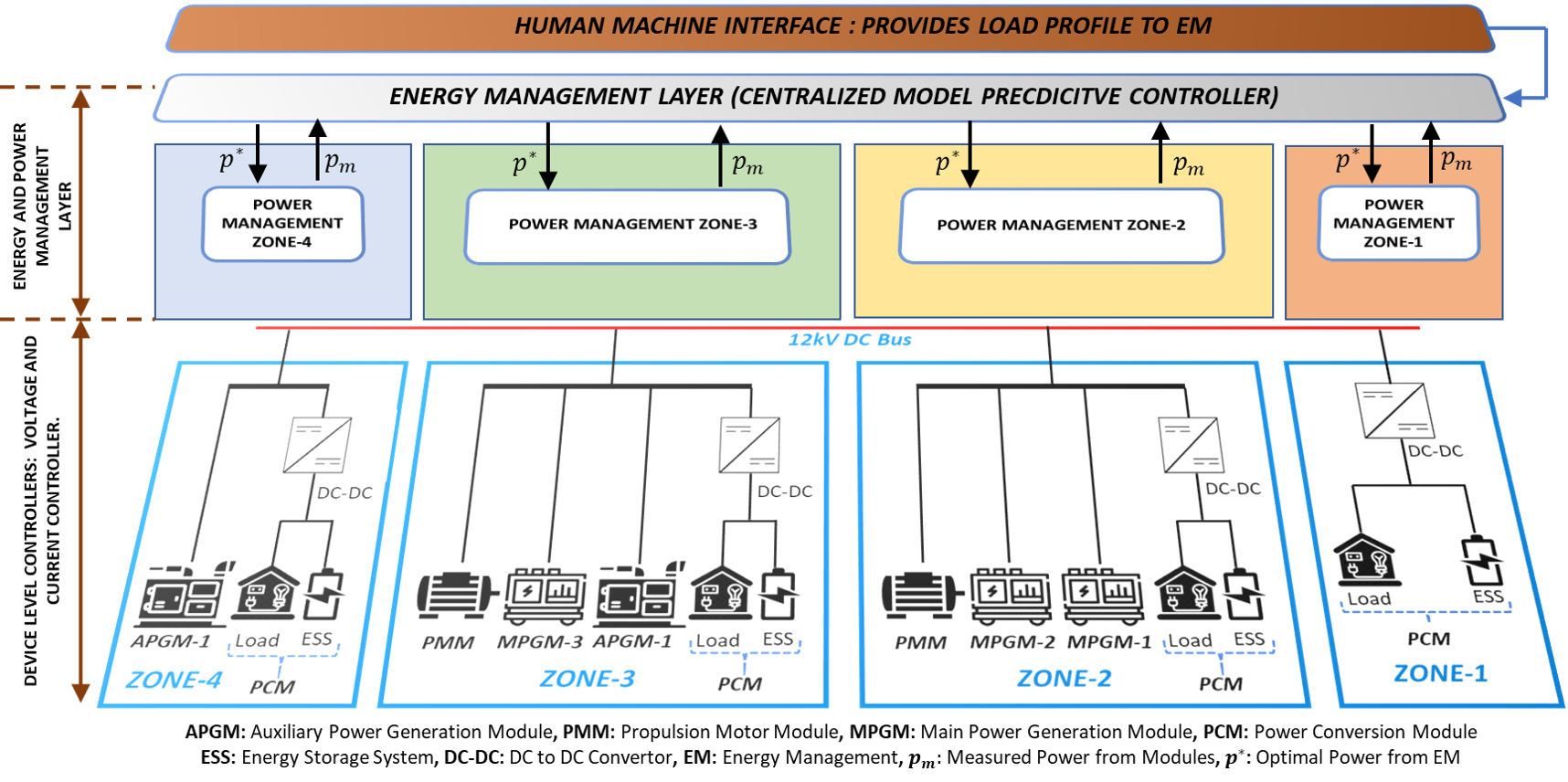}
	\caption{4-Zone notional ship power system hierarchical control structure}
 \label{SPS_hierarchical}
\end{figure}

Model Predictive Control (MPC) \cite{2009_Rawlings} is a mature technology that has been extensively used to control slow chemical processes. As the computation power improved throughout the years, it has become appropriate for control of faster systems such as power electronics. MPC utilizes the model of the system and its behavior over a prediction horizon and aims to solve an optimization problem. MPC can systematically include physical limitations such as State of Power (SOP) and ramp-rate limitations. Generally, MPC optimization includes an objective function with weighted expressions that often address a trade-off. The simplicity of addressing this trade-off compared to the complexity of conventional algorithmic approaches is a major advantage of MPC control. However, MPC problems tend to have a computation burden that might not be viable in faster operations such as in real-time implementations. Hence, control system designers should always be aware of the computation cost of model predictive controllers.

SPS is considered an islanded MVDC MG as shown in Figure \ref{SPS_hierarchical}. The notional SPS is divided into four zones. Each zone may include one or a combination of power generating, storage, and load modules. For example, zone-2 includes two main power generation modules (PGMs) each including fuel-operated generators as well as three-phase rectifiers, filters, and respective device level controllers (DLCs).   There are power conversion modules (PCMs) which include multiple power converters, energy storage systems, and AC and DC loads. High power and high ramp-rate loads are represented by propulsion motor modules  (PMMs), super loads (SLs), and AC load centers (ACLC) for DC and AC loads respectively. Zones are connected to a unified $12kV$ DC power-line which should ensure appropriate inter-zone power and energy transactions and zone and vessel-wide reconfiguration. 

\subsection{Ship System Power-Flow Model}
The power-flow model must oblige the power generated to equal the power demanded criteria. Thus, The baseline SPS power flow can be represented as:    
\begin{align} \label{SPS_P}
f(p_g,p_b,p_L)=0
\end{align}
where, $p_g$, $p_b$, and $p_L$ represent the generator sets (gen-sets), batteries, and load powers. In (\ref{SPS_P}), generators are considered to inject power in a unidirectional form while ESS can inject and absorb power bidirectionally to support the load demand. One major consideration for power distribution is the ramp-rate limitation of generation and storage elements. Critical loads may have higher power ramp-rate demand than generators. A general approach is to request power from ESSs to compensate for the high ramp-rate portion of the power demand which enables the operation of the generators within their response capabilities. This is shown to improve the survivability, stability, and quality of electrical levels of the SPS.

\subsection{Hybrid Electric Vehicles}
\textit{First} patented design for a three-wheeled vehicle can be traced back to 1886 A.D. Karl Benz applied for a patent for his design of \say{vehicle-powered by gas engine} and called it \say{Motorwagen}. In 1908 Ford created \say{Ford Model T} which went on to become the first car to be mass-produced. During the mid 20$^{th}$, the automobile industry was booming all around the world, especially in Germany. Volkswagen, Bayerische Motoren Werke (BMW), and Mercedes-Benz became popular vehicle manufacturers during that period in Germany. While Ford and Mitsubishi were leading manufacturers in the USA and Japan. The developments in the Internal Combustion Engine (ICE) technology enhanced the global market for cars. With the emergence of Japanese brands like Toyota and Honda which made cost-effective cars, more people were able to afford cars. By the early 1990s cars became the main mode of transportation in America. Engine developments led to more faster cars, but it came at the expense of fuel efficiency and an increase in Greenhouse Gas (GHSs) emissions. It is predicted that in the upcoming 50 years, the global population might reach the 10 billion mark. This is going to increase the car usage by 200$\%$. If ICE is used as a driving mechanism for these cars, there is a high possibility that we might run out of fossil fuels needed to run the engine.  

\begin{figure}[h!]
	\centering
\includegraphics[width=0.6\textwidth]{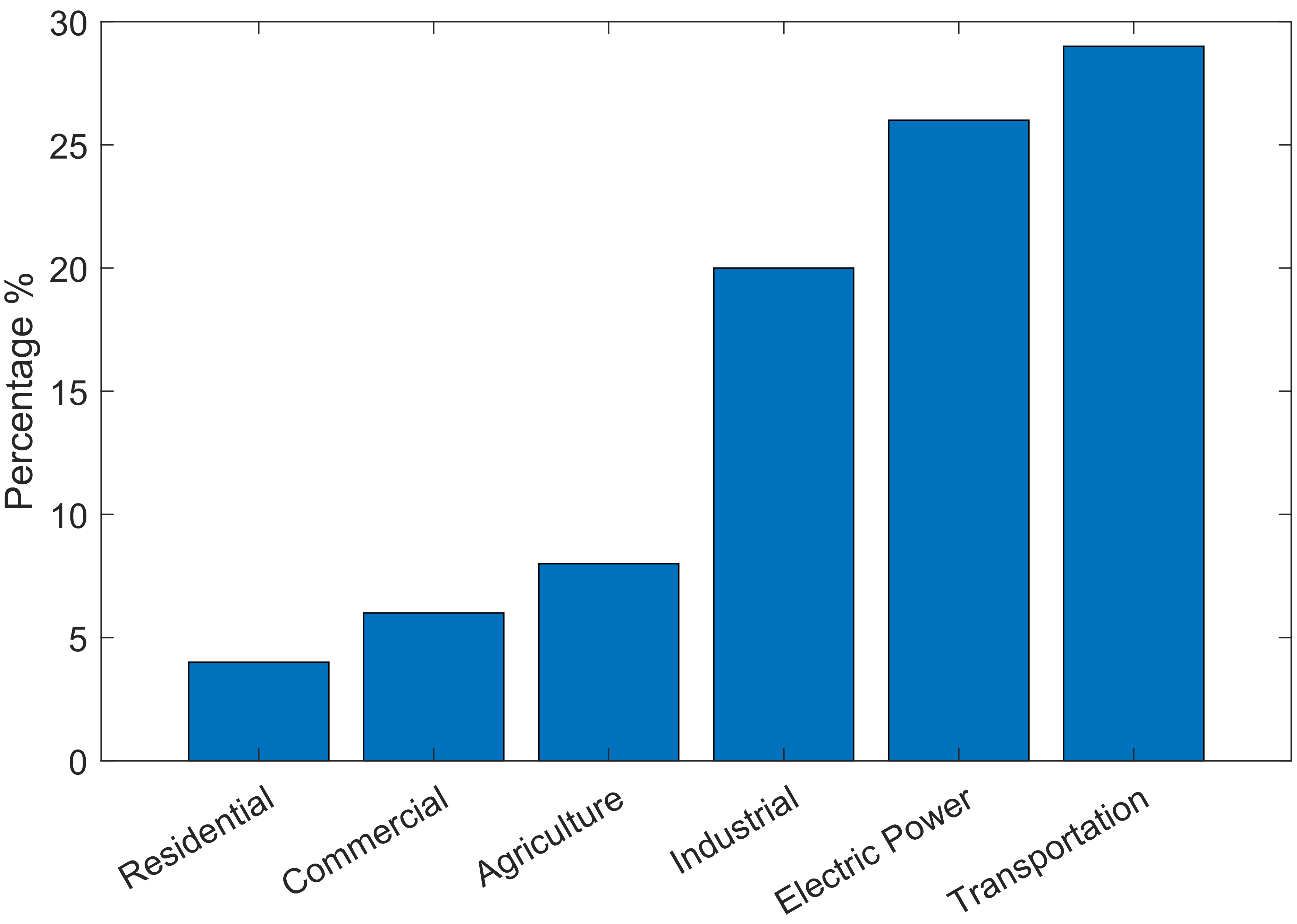}
	\caption{Graph showing greenhouse gas emissions by sector in USA}
    \label{Emissions}
\end{figure}

\textit{Transportation} sector contributes to 29$\%$ of the greenhouse gas emissions as seen in Fig \ref{Emissions}. This calls for the government's attention to designing policies that can reduce the GHG contribution from the transportation sector. The term \textit{Carbon Footprint} was coined during the mid-2000s. Since then the goal of cutting GHG emissions by manufacturing companies is being described in terms of reducing carbon footprint. United Nations Framework Convention on Climate Change (UNFCC) in 1992 based on a scientific consensus set forth set of regulations in tackling the emission of greenhouse gases. Their main goal was to stabilize the concentrations of GHG around the world. Aiming towards the goal various protocols such as Kyoto Protocol and Paris Climate Accord were introduced. Kyoto Protocol was ratified by at least 55 countries in 1997. It was at this point automobile companies started looking towards other alternatives. In 1966 the state of California started setting standards for vehicle tailpipe emissions.  More research focus was diverted towards developing a hybrid vehicle. The concept behind it was to design a drive-train that utilizes both engine and battery power, thus reducing the usage of the engine improving fuel efficiency, and mitigating emissions. The introduction of Toyota Prius in 1998 in the Japanese market and later in 2000 in the American market can be seen as the birth of HEV mass production on an assembly line. HEVs emit less GHGs compared to traditional gasoline-based cars. The emission regulations are what led to a substantial rise in HEVs. In the next section, we will review in detail manner about emergence of HEVs. 

\subsubsection{Emergence of Hybrid Electric Vehicles}

Stringent regulations on environmental emissions such as the Paris Climate Accord and also fears of extinction of oil and natural gases shortly led the automobile industry to focus on manufacturing more alternative energy-dependent, efficient, and environmentally clean vehicles as shown in Figure \ref{Emissions}. Major automotive companies have started investing heavily in developing vehicles such as Fuel cell hybrid vehicles (FHVs), Battery Electric vehicles (BEVs), and Hybrid electric vehicles (HEVs). Numerous challenges are currently being faced during the development of new technologies. FHV technology is still immature and might become dominant in upcoming years. BEV technology is gaining momentum. But due the extent to which BEVs can persist on a single charge and their long charging duration make them extremely difficult for long-distance travel. Another disadvantage is their high initial cost. But with the development of fast battery chargers and advances in the Li-ion battery model, this technology looks promising in tackling emission issues. The concept of Hybrid electric vehicles (HEV) emerges as a solution to address both problems posed above. One of the main advantages is not requiring the incorporation of an external battery source. But the most recent development is the emergence of Plug-in Hybrid electric vehicles (PHEVs). The main classification between HEVs and PHEVs is PHEVs rely more on electrical components. They can run entirely on battery power for up to 80 miles. PHEV technology is also being used in light-duty vehicles too. However further developments in PHEV technology heavily depend on the development of battery technology. Over the years various designs of HEV have been proposed, the prominent one being \cite{AFburke}.

The surge in production of HEVs in the past decade accounts for their relatively high fuel efficiency compared to traditional gasoline-based vehicles and meeting emission requirements without drastically compromising vehicle performance. HEVs can be configured in multiple ways namely: \textit{Series hybrid configuration}, \textit{Parallel hybrid configuration}, \textit{Series-Parallel hybrid configuration} and \textit{Complex hybrid configuration} \cite{Maggetto}\cite{CCChan}. The series-parallel configuration has been extensively used during recent times as the start-up operation of the vehicle is solely through battery power, thus enhancing fuel efficiency. All configurations are conglomerations of multiple components such as the Internal combustion engine (ICE) and electric machines (EMs) which include a generator and motor, whose functionality can alternate depending on the driving mode and battery pack. Usually, ultracapacitors or other additional storage elements are used in combination with battery packs to meet the specific energy and power requirements \cite{Mehrdad}. Table \ref{Char_EV_tab} shows the characteristics of the BEVs and the HEVs.

\begin{table}[ht]
\caption{Characteristics of battery electric vehicle and hybrid electric vehicle}\label{Char_EV_tab}
\begin{center}
 \begin{tabular}{||c|c|c||} 
 \hline
{\textbf{EV Subsystem Type}} &  {\textbf{Battery EV}} & {\textbf{Hybrid EV}} \\ [1ex] 
 \hline\hline
 Propulsion &  Electric Motor & Electric Motor, ICE \\ [1ex]
 \hline
 Energy Storage &  Battery, Ultra Capacitor  &  Battery, Ultra Capacitor \\[1ex]
 \hline
 Infrastructure & Grid Charging & Gasoline Stations\\[1ex]
 \hline
 Characteristics & Zero Emissions, High Efficiency & Low Emissions, High Economy\\[1ex]
 \hline
 Major Issues &  Charging Facilities & Battery Sizing \\[1ex]
 \hline
 \end{tabular}
\end{center}
\end{table}

Optimal Energy management (EM) between engine and battery plays a pivotal role in enhanced fuel efficiency of HEV, but achieving it is extremely challenging due to changing modes of operation of HEV. \cite{serra} proposes an EM strategy that takes into account the joint system i.e. the coupling between various modules and also HEV configuration. An MPC-based energy management of a power split was proposed in \cite{Borhan}. Advancements in data processing capabilities of machines have led to the use of control allocation methods in HEVs \cite{tor}\cite{Bodson}. HEVs are over-actuated in nature, utilizing both a generator and an energy storage device. Each component is integral for the peak performance of the HEV, should one fail the reliability of the vehicle drops considerably. Considering the longitudinal vehicle motion dynamics, \cite{chen} proposes adaptive energy-efficient control allocation for distributing total control effort in an over-actuated electric ground vehicle. EM techniques are studied in detail in future sections. 
\subsection{Hybrid Electric Vehicle Model}
HEV configurations play a key role in power distribution between Internal combustion engine (ICE) and Energy storage elements. Four different configurations of HEV drive-train are available \cite{Mehrdad} namely:
\begin{itemize}
    \item Series hybrid electric vehicle configuration
    \item Parallel hybrid electric vehicle configuration
    \item Series-Parallel hybrid electric vehicle configuration
    \item Complex hybrid electric vehicle configuration
\end{itemize}

In series HEV configuration, engine torque output is converted into electrical energy using a generator. This electricity can be used to propel the car forward through the motor and transmission system or charge the battery when needed. This model is assisted by three propelling devices engine, a generator, and an electric motor. This drastically reduces the efficiency of HEV in this configuration. Perfect sizing of all the components is also difficult to achieve. In this configuration engine or motor can be used in a decoupled fashion. The motor also plays a pivotal role during regenerative braking. The power generated while braking is used to charge the battery. Parallel HEV configuration enables simultaneous or parallel operation of the engine and motor to drive the wheels. The engine and Motor are coupled to a common driveshaft using dual clutch. Power required to drive the wheels can be supplied by the engine alone motor alone or both. Parallel configuration runs on two propelling devices engine and motor, so efficiency is enhanced compared to series configuration. The electric motor is used as a generator during regenerative braking to charge the battery \cite{CCChan}.

\begin{figure}[h!]
	\centering
	\includegraphics[width=0.65\textwidth]{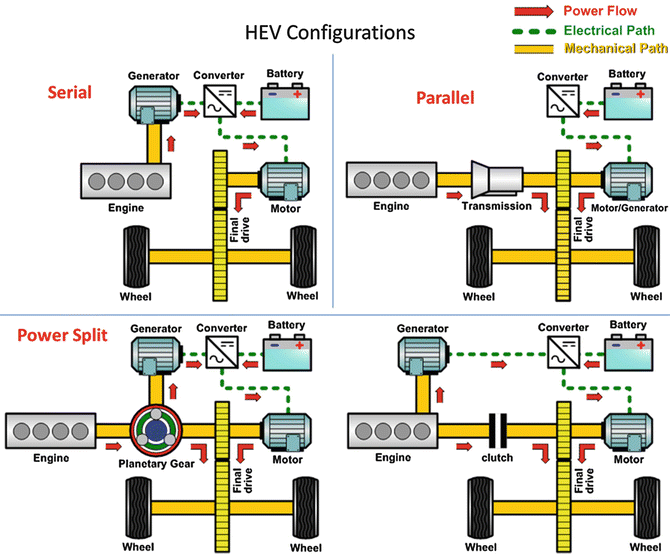}
	\caption{Hybrid electric vehicle power flow path for different configurations \cite{2013_Rizzoni}}
    \label{Power Flow}
\end{figure}

Series-Parallel HEV (SPHEV) has been a widely used configuration in HEVs in recent years. The main advantage of SPHEV is it uses only battery drive-train during startup operation. A more detailed schematic layout is discussed in the power-flow model section. In Series-Parallel or Power Split configuration the electrical and mechanical drive-trains are coupled with the help of a planetary gear system. The planetary gear system couples the engine, motor, and generator (Figure \ref{Power Flow}). The planetary gear consists of a carrier, ring, planet, and sun gears. The idea behind the operation of planetary gear was first proposed and used by Greeks during 500 B.C. The same principle was used by Italian engineer Agostino Ramelli in his \textit{bookwheel} design. 

\subsubsection{Power Flow Model}
The power flow path is distinctive of the operation mode. The main modes of operation are \textit{acceleration}, \textit{deceleration}, \textit{regenerative braking}, \textit{charging while driving} and \textit{charging while standing still}. During \textit{acceleration} the power demand is met by both battery and engine \cite{Mehrdad}. During regenerative braking mode, the mechanical drive-train is decoupled and all the power is used to charge the battery using the electrical drive-train. In HEVs the power flow function can be represented as follows:
\begin{equation}
    f(p_b,p_e,p_d) = 0,
\end{equation}
where $p_b$ is power supplied by the energy storage elements, $p_e$ is power supplied by the mechanical component engine and $p_d$ is the power demand.

\begin{table}[h!]
\caption{HEV system model design and references}
\begin{center}
 \begin{tabular}{||c|c||} 
 \hline
 {\textbf{System Modelling}} &  {\textbf{List of References}} \\ [1ex] 
 \hline\hline
 Drive-line Dynamics &  \cite{bailey},\cite{BKPowell},\cite{1997_Cikanek} \\ [1ex]
 \hline
 Internal Combustion Engine (ICE) &  \cite{syed},\cite{Heywood} \\ [1ex]
 \hline
 Electric Machine &  \cite{bailey},\cite{BKPowell} \\ [1ex]
 \hline
 Energy Storage Elements &  \cite{2021_Vedula},\cite{BKPowell},\cite{2008_Kroeze} \\ [1ex]
 \hline
 Battery Life Estimation Models & \cite{2009_Marano},\cite{2009_Serrao},\cite{2010_Filippi} \\ [1ex]
 \hline
 \end{tabular}
\end{center}
\end{table}


\subsection{Energy Management in Shipboard Power Systems}
Both new and existing warships consist of many advanced loads integrated into the ship power system. Electromagnetic railguns and electric propulsion motors are a few of the advanced loads. The ever-increasing loads in SPS have high ramp rates. The generator's ramp rate is usually not sufficient to support the high ramp rate loads. Also, the weight and size constraints limit adding more number of generators. This calls for the integration of energy storage systems (ESSs) into the SPS. The sum of all the loads must be met by the power supplied by the generators and ESSs installed on the ship.  In \cite{Gonsoulin} a centralized MPC for multiple ESS control in SPS was proposed. A high-level controller which is an EM layer acts as a control layer for ESSs control. It gives out power commands to generators and ESSs to ensure that the load demand is met. SPS power management using constrained nonlinear MPC was proposed in \cite{Stone}. The proposed controller performs a similar optimal power split between generators and ESSs as mentioned in the previous case. The end goal is to meet the load demand in real-time while abiding by the optimization constraints over a finite future horizon. The concept of PMC proposed is not just limited to power management in SPS, but also, reconfigure and add or shed loads based on predefined load priorities.

A predictive control for EM in SPS under high-power ramp-rate loads was proposed in \cite{2017_Vu_3}. This paper specifically focuses on challenges encountered in designing a controller under high ramp rate conditions. The hierarchical structure of EM, Power management (PM), and Local device control layers was proposed. The device layer regulates the output voltage of generators and the output voltage of ESSs. The PM layer includes power managers, which regulate power to maintain bus stability. The EM layer exchanges information to and from lower layers and with each other through a communication network. Particle swarm optimization-based MPC power management for SPS was proposed in \cite{Paran}. All of the above-discussed EM methods are centralized MPC-based. The table below shows the risk involved in SPS and its modeling and applied control allocation technique. The table also shows currently used approaches to address risk and the advantages of the proposed approach.

\subsection{Energy Management in Hybrid Electric Vehicles}

\begin{figure}[h!]
	\centering
	\includegraphics[width=0.6\textwidth]{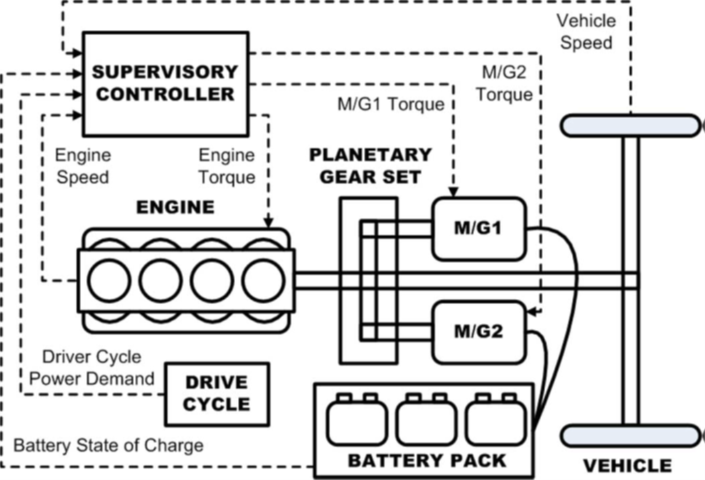}
	\caption{Hybrid electric vehicle control schematic \cite{5439900}}
    \label{HEV Control}
\end{figure}

In the previous sections, we have discussed the emergence of HEVs and HEV sub-systems design. In this section, we will focus on the \textit{Energy management} (EM) strategies in an HEV. EM strategies developed during the 1990s were mostly based on centralized optimization Figure \ref{HEV Control}. Bailey \cite{bailey} was one of the early leading people to curate all the information about HEV sub-system modeling into a single paper. In \cite{serra} authors presented a joint system model. They designed a special EM strategy taking into consideration the join system dynamics of HEV. The rule-based EM technique was widely popular during the early stages of HEV development. The idea behind this method is that the designer creates a set of rules for the system to abide by. It is often difficult to encapsulate and capture all the system rules for smooth working. Rule-based control design was also presented in this paper. The idea behind the joint system is as follows: \say{\textit{Joint system is an autonomous motion system that can drive a rotating load, delivering the demanded power output}}\cite{serra}.

Increased computational capabilities have led to a substantial increase in the use of high computation-based optimization methods such as MPC. An MPC-based EM of a power-split in HEV is proposed in \cite{Borhan} by Borhan. Since the enhancement of fuel economy is one of the main goals of HEV, this paper formulates the EM problem as a nonlinear and constrained optimal control problem that minimizes fuel input to the engine. Multiple cost functions are defined and MPC is used to solve the power-split problem between engine and battery for a given sample time at every time step. The nonlinear driveline dynamics are linearized for MPC implementation. Most of the developed and proposed EM strategies made minimizing fuel consumption their main goal. In \cite{8816673} East formulated an MPC-based optimization problem that minimizes fuel consumption and used ADMM to solve the problem. East provided convex formulations of the MPC optimization problem, which enables guaranteed convergence. Considering the real-time traffic conditions Jiao  \cite{2014_Jiao} developed an optimization problem that can manage the energy in real-time. This paper proposes a control strategy based on stochastic dynamic programming (SDP). As with other EM strategies, the ultimate goal of this strategy is to minimize the average fuel consumption and address specified battery and engine constraints while meeting the desired power demand. The EM problem is formulated as a stochastic nonlinear and constrained optimization problem. 

Over the past two decades, substantial developments have been seen in the field of HEV modeling. The components in HEV such as the Battery Model and Engine Model got complicated over the years. A new dimension in the form of managing the battery State of Charge (SOC) became a key part of the MPC optimization problem. More recent research developments focused on managing battery degradation. The battery performance degradation relies on multiple factors such as C-rate, SoC, etc \cite{2010_Filippi}. These factors play an important role in extending the operation life of a HEV. Decentralized control for PHEVs under driving uncertainty was discussed by Vaya \cite{7028989}. ADMM method was used in this paper to formulate a decentralized optimization problem. Local optimization problems are solved in parallel and the aggregator enables the flow of information between nodes to eventually find iteratively the global solution. If the formulated cost objective is convex, then a global optimum is attained. Vaya considers the driving pattern uncertain in this work.

\newpage
\section{Key Aspects Covered in this Document}
\begin{figure}[h!]
    \centering
    \includegraphics[width=0.9\textwidth]{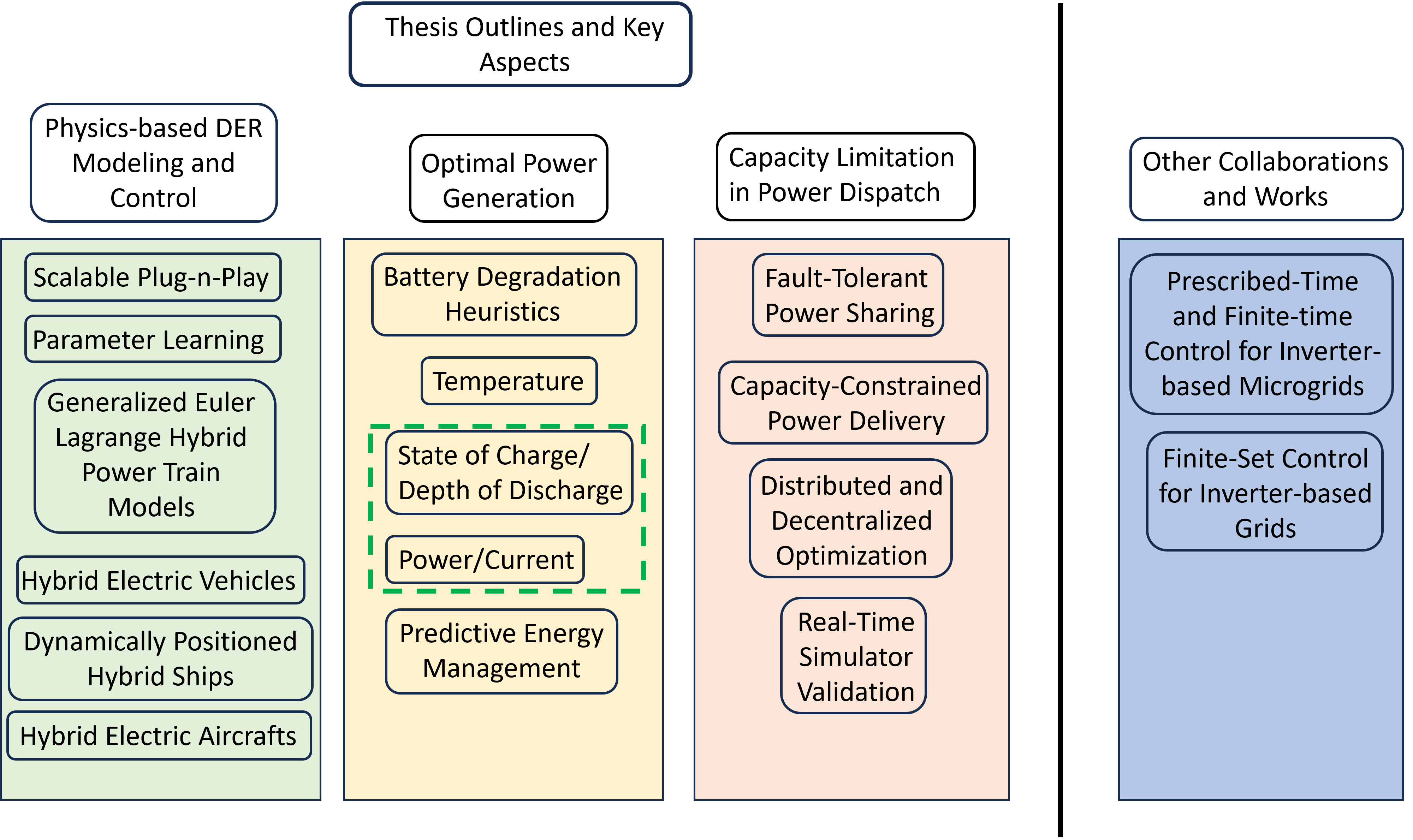}
    \caption{Outline of this document and key aspects and collaborations}
    \label{fig:thesis_outline}
\end{figure}

\section{Mathematical Modeling of DER Components}
\subsection{Motivation}
This section is in part based on \cite{9512356}. The concept of the Micro-grid (MG) as an integration of flexible clusters of energy generating, storing, and dissipating sub-systems was first presented in \cite{2001_Lasseter}. MGs are typically categorized as islanded and grid-tied based on their configuration \cite{2001_Lasseter,2007_Hatz_Microgrids}. The control structure employed for MGs is usually hierarchical and approaches such as centralized and distributed optimization \cite{2011_Gue}. The former is capable of achieving high levels of performance and optimization but may be susceptible to single points of failure and expandability limitations. The latter offers flexibility at the cost of computational complexity and increased control surface which further exposes the system to more uncertainty including cyber attacks. These phenomena, in turn, further increase the complexity of the distributed MG control.

Analogous to the hierarchical nature of control and resource allocation strategies, \cite{2011_Gue} proposed a hierarchical control architecture for AC and DC MGs that includes level zero, primary, secondary, and tertiary controllers. Currently, these control levels are widely used to provide device-level control (DLC), power management (PM) and energy management (EM) of MG systems \cite{2020_Bijaieh,2021_DSL}. The general form of this hierarchical framework enables the integration of many local (on-site) control methods such as droop-based controls, maximum power point tracking (MPPT), and many more. It is important to understand that various control methods should be carefully matched and eventually benchmarked versus appropriate metrics that are defined for specific systems. For example, a combination of droop control (as primary) and MPPT (as lower lever local) controls will require some sort of power curtailment due to the mismatch between variable generation and fixed power injection to the local point of common coupling (PCC). The interoperability of such combinations will also depend on the strictness of the constraints that are defined for the electrical levels of the system. Control specifications, appropriate metrics, and guidelines for systematic control evaluations for various control levels are provided in IEEE 2030.7 and 2030.8 standards \cite{2017_IEEEstd2030_7, 2017_IEEEstd2030_8}. 

The holistic nature of MGs and their control adds complexity to the device level and system-level control design, analysis, evaluation, and validation. From an investigator's perspective, the choice of MG models with appropriate fidelity is critical. There is always a trade-off between the fidelity of the chosen model and the required computation power to compile the model. The investigator must comprehend the critical assumptions, model constraints, the validity of the performed experiments as well as the appropriate platform for simulation.

There are numerous works dedicated to the modeling of MG systems. While works such as \cite{2014_Shafiee, 2017_Fran} present modeling, control, and stability of DC MG systems, they do not provide systematic guidelines for mathematical coupling, as well as a systematic approach for control design of the underlying mathematical sub-systems. Power-flow and energy transfer models of microgrid systems with a generalized and systematic control approach for the overall system is offered in \cite{2015_Hassell, 2020_Bijaieh} but do not provide guidelines for system-level mathematical coupling or in case the system needs to be scaled. For example, solving those problems in a symbolic framework would be very hard to manage. A solution for scalability for such problems can be found in symbolic mathematical modeling platforms that use Modelica \cite{Fritz_2011}, where, blocks of code can be packaged as subsystems coupled through the fundamental definition of across (voltage) and through (current) variables. 

In all cases, the investigator must choose appropriate tools that fit best to solve the problem at hand. There are advantages of acausal modeling such as done in Modelica over typical causal modeling that fits well with ordinary differential equation (ODE) solvers such as that used in Simulink \cite{2020_simulink}. Generally, multi-physics systems modeled in hybrid differential algebraic equations (DAEs), and models with discrete states are more suitable to be solved in modeling approaches such as Modelica. On the other hand, causal modeling is more structured in a sense that the system is decomposed into a chain of causal interacting blocks \cite{2013_Dizqah}. The problem of solving DAEs in simulation environments such as Simulink is very common. Moreover, inappropriate programming can also lead to the appearance of algebraic loops which may cause severe computation burden. Implications of the existence of algebraic loops in a model include the inability of code generation for the model, the Simulink algebraic loop solver might fail to solve, and while the Simulink algebraic loop solver is active, the simulation may run slowly \cite{1999_Shampine, 1970_rabi, 1980_more}. Common solutions to the algebraic loop problem are: to introduce unit delays to the blocks, and, to turn the DAEs into ODEs by introducing additional states to aid the solver in easier solve. Adding additional delays to the model may lead to inaccurate solutions and introducing additional states might unnecessarily increase the size of the model. Hence, providing systematic approaches for the mathematical modeling of such dynamical systems is a viable path to pursue.       

This work presents a mathematical modeling approach for MG systems. The aim is to define specific generating, transmitting and loading mathematically expressed subsystems in a modular way such that the designer would be still capable of scaling or extending the overall simulated system. While the work presented here might seem iterative in part with respect to previous efforts, the authors are compelled to use the result of this work due to its appropriateness for control system design. Moreover, based on the current pedagogy in fundamental controls in academia, students are pushed to utilize mathematical models rather than simulation packages such as  \cite{simscape} hence a modeling approach that aligns well with current pedagogy may be very useful. The characteristics of the mathematical modeling approach in this paper are: (1) Ease of modeling; where the investigator uses fundamental blocks to create ODEs. (2) Modular mathematical blocks connected according to a general connection convention. (3) An LBM that is fast to run and appropriate for control development specifically at the EM level. (4) The model does not include algebraic loops and is readily made for code generation, and subsequently real-time implementation in platforms such as Simulink Real-time. (5) The model can be run whole, or in part, on any platform that allows advanced mathematical operations. (6) The system may run at variable sub-system fidelity; reducing computation for certain applications. (7) The model and control development process is more aligned with common practices and pedagogy.   

The shortcomings associated with this type of modeling include (1) all the limitations and assumptions that are inherited from use of ideal, simplified, and average models. Integration of switching functionality in an MG setting is possible but significantly increases the complexity of the model. (2) The models are not suitable for fault analysis since the bandwidth (BW) of operation will not be typically sufficient. (3) More complexity due to the management of voltage and current signals and dealing with both KVL and KCL separately at device and system levels.                    

\subsection{Generators}
 Two generators namely a DC Generator, and an AC Generator modeling are presented in this section. The generator models presented are detailed enough to demonstrate the control approach designed for this discussion.

\subsection{DC Generator Model}
A DC Generator is modeled as a controllable voltage source in series with a \textsf{RL}-line and shunt capacitance. The mathematical model for the DC Generator is given as follows \cite{2021_Vedula}:
\begin{align}
    l_g \frac{di_g}{dt} &= -r_g i_g(t) + v_g - v_c, \\
    c_g \frac{dv_c}{dt} &= i_g - i_{g_{ref}},
\end{align}
where $i_g \in \mathbb{R}$ is the generator current, $l_g, r_g, c_g\in \mathbb{R}_+$ are the inductance, resistance, and the capacitance of the generator. The capacitance voltage/bus voltage $v_{c} \in \mathbb{R}$ and the controlled generator voltage $v_g \in \mathbb{R}$. The control objective of the DC Generator is to regulate the generator voltage $v_c$ or to track a nominal voltgae $v_{{g}_{ref}}$. Consequently, the objective is to design a control law that satisfies the following objective
\begin{align*}
\int_0^{\infty} \bigg(\underbrace{v_c(t)-v_{g_{ref}}(t)}_{\tilde{v}_g}\bigg)^2 dt < \infty.
\end{align*}
This objective can be designed in multiple ways ranging from a simple PID controller to a more complex model based nonlinear control methods. A PID controller can be designed as follows
\begin{align*}
v_g = v_c + k_P \tilde{v}_g(t)+k_I \int_0^{t}\tilde{v}_g(\tau)d\tau,
\end{align*}
where $k_P,k_I \in \mathbb{R}$ are the device level current control gains. Assuming no model knowledge an adaptive control law for parameter learning concurrently tracking the desired voltage is presented next. Consider the voltage tracking error as
\begin{align*}
    e = v_c - v_{{g}_{ref}},
\end{align*}
taking the derivative along the time variables yields,
\begin{align*}
    \dot{e} &= \dot{v}_c, \\
    c_g\dot{e} &= c_g\dot{v}_c, \\
    c_g \dot{e} &= i_g - i_{{g}_{ref}}.
\end{align*}
Now, for some $\alpha > 0$ consider the filtered error is considered as follows
\begin{align*}
    \eta = \dot{e} + \alpha e,
\end{align*}
taking the first derivative along the variables yields,
\begin{align*}
    \dot{\eta} &= \ddot{e} + \alpha \dot{e}. \\
    l_g c_g \dot{\eta} &= -r_g i_g + \alpha l_g i_g - \alpha l_g i_{{g}_{ref}} + u, 
\end{align*}
where $u \triangleq v_g - v_c$. The dynamics can be linearly parameterized into known system measurements and unknown system parameters. Consequently,
\begin{align}\label{dc_gen_dynamics}
    l_g c_g \dot{\eta} = Y\boldsymbol{\theta} + u,
\end{align}
where 
\begin{align*}
    Y = \begin{bmatrix} -i_g & \alpha(i_g - i_{{g}_{ref}}) \end{bmatrix}, \hspace{3mm} \boldsymbol{\theta} = \begin{bmatrix} r_g \\ l_g \end{bmatrix}.
\end{align*}
Now, consider the control law
\begin{align}\label{dc_gen_control_law}
    v_g = v_c - Y \boldsymbol{\hat{\theta}} - k\eta.
\end{align}
Substituting \eqref{dc_gen_control_law} in \eqref{dc_gen_dynamics} yields,
\begin{align}\label{dec_gen_closed_loop}
    l_g c_g \dot{\eta} = Y \boldsymbol{\tilde{\theta}} - k\eta,
\end{align}
where $\boldsymbol{\tilde{\theta}} = \boldsymbol{\theta} - \boldsymbol{\hat{\theta}}$. The parameter adaptation law is chosen as
\begin{align}\label{dc_gen_param_adapt}
    \boldsymbol{\dot{\hat{\theta}}} = Y^\top \eta.
\end{align}

   Now, Consider the Lyapunov function candidate as follows
    \begin{align*}
        V = \frac{l_g c_g}{2} \eta^2 + \frac{1}{2}e^2 + \frac{1}{2}\boldsymbol{\tilde{\theta}}^\top \boldsymbol{\tilde{\theta}}.
    \end{align*}
    Taking the first derivative and substituting the dynamics in \eqref{dec_gen_closed_loop} yields,
    \begin{align*}
        \dot{V} &= \eta l_g c_g \dot{\eta} + e \dot{e} - \boldsymbol{\tilde{\theta}}^\top \boldsymbol{\dot{\hat{\theta}}}, \\
        &= \eta(Y\boldsymbol{\tilde{\theta}} - k\eta) + e(\eta - \alpha e) - \boldsymbol{\tilde{\theta}}^\top \boldsymbol{\dot{\hat{\theta}}}.
    \end{align*}
    Substituting the parameter update law in \eqref{dc_gen_param_adapt} and using Young's inequality yields,
    \begin{align*}
       \dot{V} \leq -(k-\frac{1}{2})\eta^2 - (\alpha - \frac{1}{2})e^2. 
    \end{align*}
    $\dot{V}$ is negative semi-definite implying that $V \in \mathcal{L}_{\infty}$ which implies $\eta, e, \boldsymbol{\tilde{\theta}} \in \mathcal{L}_{\infty}$. Integrating $\dot{V}$ yields,
    \begin{align*}
        V(\infty) - V(0) \leq -(k-\frac{1}{2})\int_{0}^{\infty}\eta^2(\tau) d\tau - (\alpha-\frac{1}{2})\int_{0}^{\infty}e^2(\tau) d\tau.
    \end{align*}
    Implying that $\eta, e \in \mathcal{L}_2$. From Barbalat's lemma, it can be established that as $t \longrightarrow 0$, $\eta \longrightarrow 0$ and $e \longrightarrow 0$.
\begin{figure}
    \centering
    \includegraphics[width=0.65\textwidth]{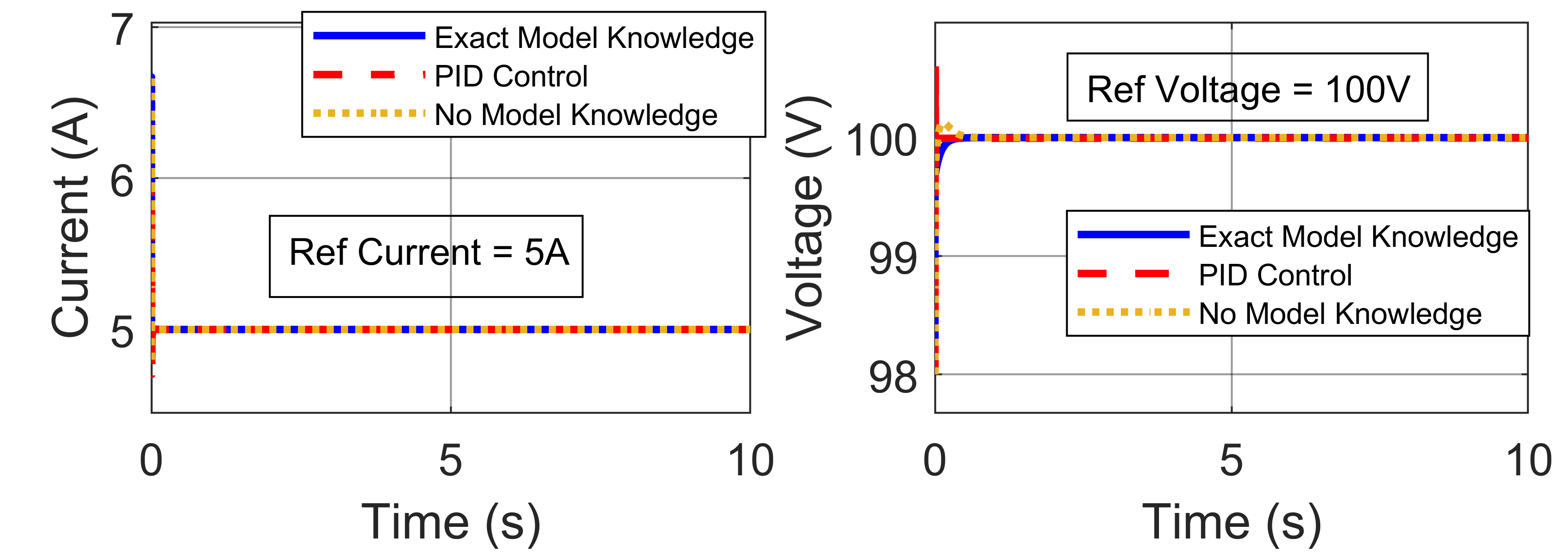}
    \caption{Current and voltage tracking for different controller designs}
    \label{fig:i_v_dcgen}
\end{figure}
The control method demonstrated above is one of the many approaches that can be employed to design the controller. Figure \ref{fig:i_v_dcgen} shows the effectiveness of different control design approaches on the current and the voltage tracking.
 
\subsection{AC Generator Model}
AC generator model consists of a synchronous generator driven by the prime mover which is powered by the gas turbine. The main goal of the generator is to regulate the grid frequency and the grid voltage around a desired nominal value. For example, in the United States (U.S.) the nominal grid frequency is 60\textsf{Hz}. The dynamics of the generator in $dq$ coordinate system are given as follows \cite{levron2018tutorial}
\begin{subequations}\label{GFrG}
    \begin{align}
     \frac{d\omega_m}{dt} &= -D\omega_m(t) + T_m(t) -T_e(t), \\
     l_s \frac{d\mathbf{i}}{dt} &= -\bigg(r_sI_2-l_s\omega_e(t) J\bigg)\mathbf{i}(t)-\mathbf{v}_c(t) + \begin{bmatrix} 0 \\ v_e(t) \end{bmatrix}, \\
     c_s\frac{d\mathbf{v}_c}{dt} &= \mathbf{i}(t) - \mathbf{i}_{L}(t) + c_s\omega_e(t) J \mathbf{v}_c(t),
    \end{align} 
\end{subequations}
where $D \in \mathbb{R}_+$ is the damping coefficient, $T_m(t) \in \mathbb{R}$ is the mechanical torque in \textsf{N-m}, $T_e(t) \in \mathbb{R}$ is the electrical torque in \textsf{N-m}. $\mathbf{i}(t), \mathbf{v}_c(t) \in \mathbb{R}^2$ are the generator current and the terminal voltage. $v_e(t) \in \mathbb{R}$ is the excitation voltage provided by the DC source and is expressed in terms of the excitation current (which is a control input) as $v_e = \omega_e l_f i_f$, $l_f$ is the field circuit inductance, and $i_f \in \mathbb{R}$ is the field circuit current. $l_s, r_s, c_s \in \mathbb{R}_+$ are the inductance, resistance and capacitance in \textsf{Henry}, \textsf{Ohm}, \textsf{Farad}. $\omega_m(t), \omega_e(t) \in \mathbb{R}$ represent the mechanical and electrical frequency and are related by the number of poles ($p \in \mathbb{N}$) of the synchronous generators $\omega_e = \omega_m p/2$. The electric torque is given as $T_e = \frac{3p}{4}l_f i_f i_q$.

\begin{figure}
    \centering
    \includegraphics[width=0.7\textwidth]{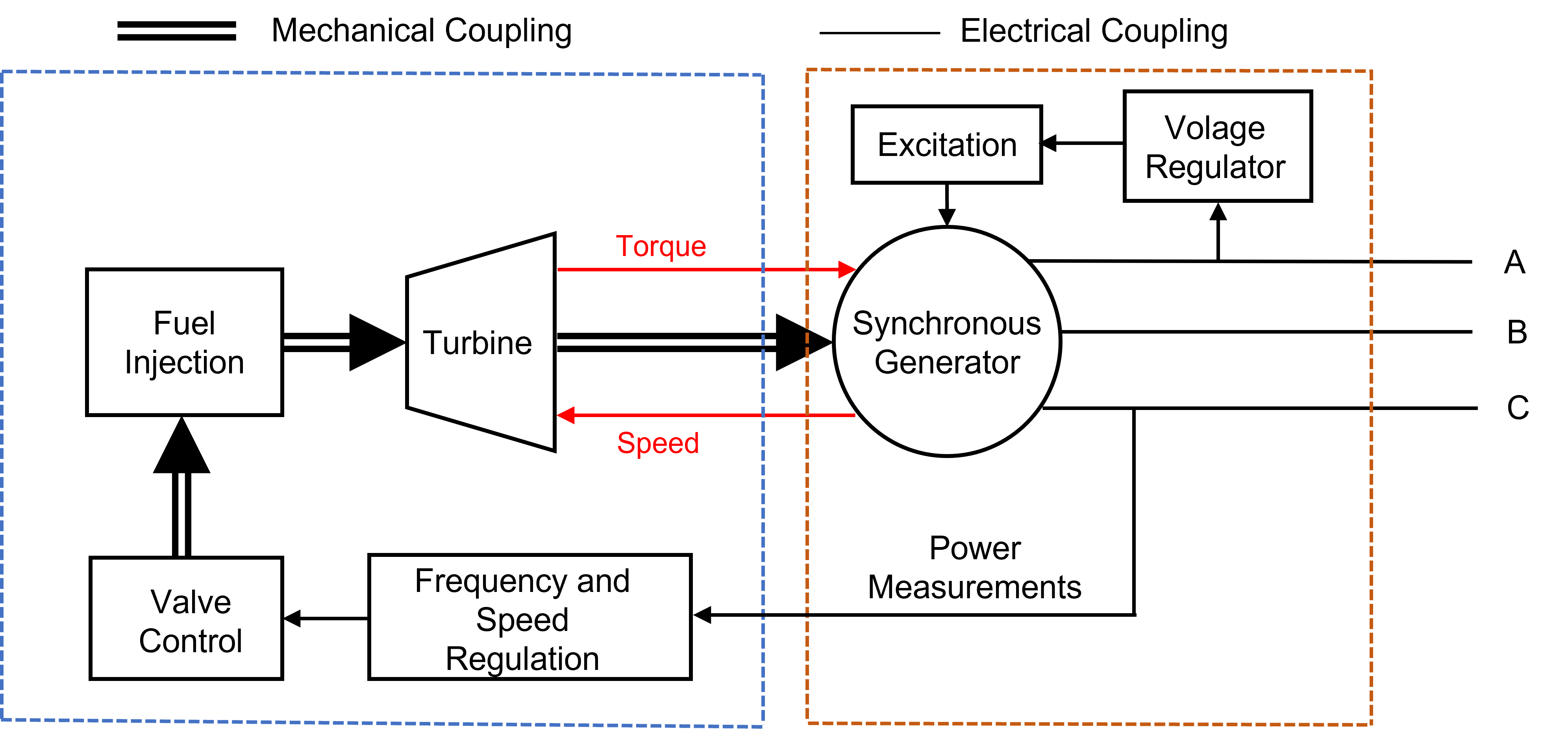}
    \caption{Schematic of effort and flow variables for AC synchronous generator}
    \label{fig:ac_sync_effortflow}
\end{figure}

The control objective is to regulate the grid frequency (for example 60\textsf{Hz} in the United States) to a nominal value while simultaneously regulating the grid voltage to a nominal value. Consequently, the speed control is designed to as follows
\begin{align*}
    T_m = k_P \tilde{\omega} + k_I \int_{0}^{\infty} \tilde{\omega}(\tau) d\tau,
\end{align*}
where $\tilde{\omega} = \omega_m - \omega_{ref}$. For instance, the US grid operates at 60\textsf{Hz}, therefore the electrical reference $\omega_e = 2\pi 60$, translating it into mechanical reference $\omega_{ref} = \frac{4\pi 60}{P}$. This sets the synchronous generator at the desired speed consequently setting the grid frequency. The other control objective is to regulate the generator's terminal voltage to a desired known value. The excitation mechanism and the $d$-axis voltage are used to regulate the capacitance/bus voltage to the nominal value. Consequently, the control law is given as
\begin{align*}
    i_f = k_P (v_{c_d} - v_{ref}) + k_I \int_{0}^{\infty} (v_{c_d}(\tau) - v_{ref})d\tau + \frac{v_{ref}}{\omega_e l_f}.
\end{align*}

\begin{figure}[t!]
    \centering
    \includegraphics[width=0.7\textwidth]{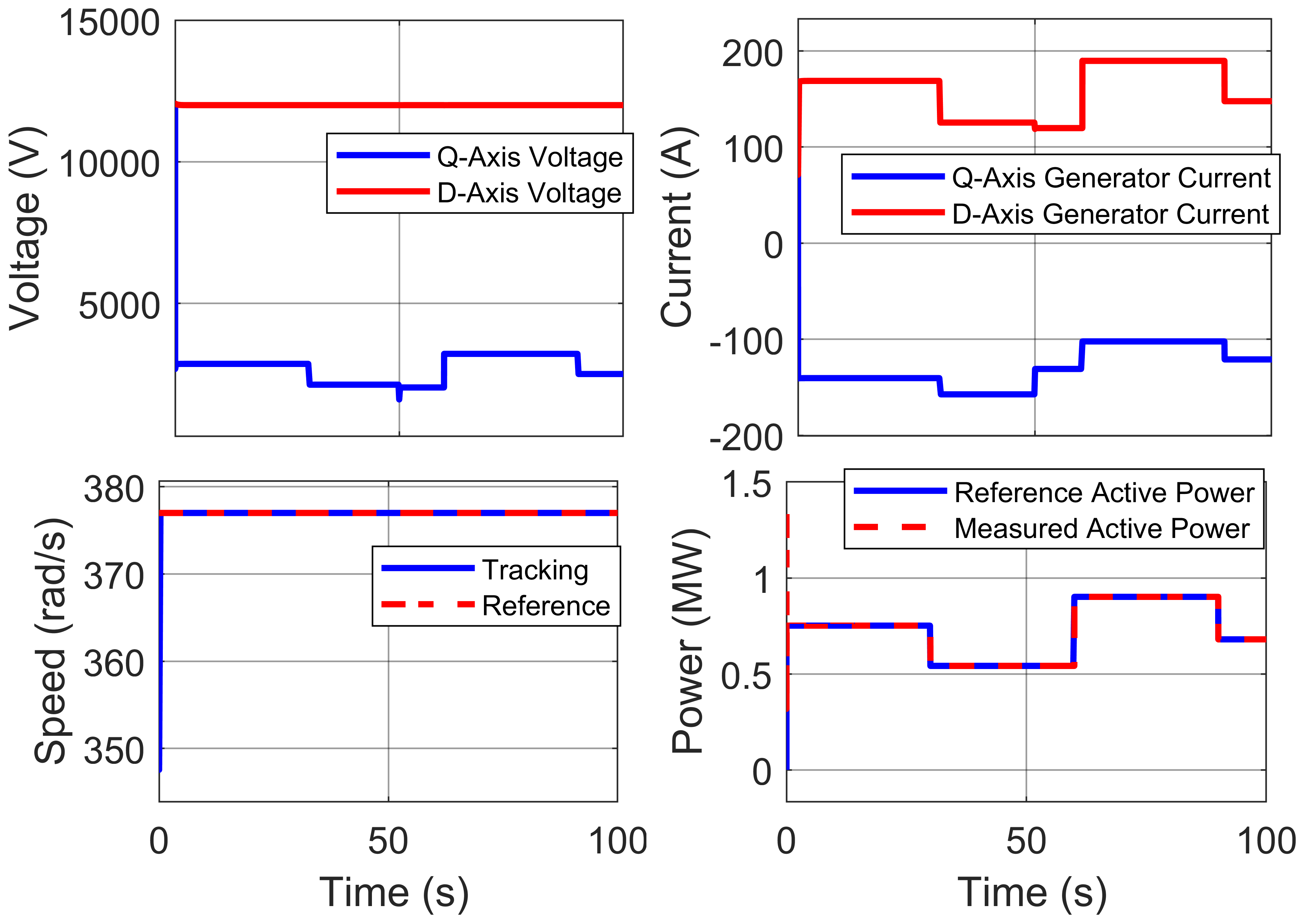}
    \caption{The voltage, speed, currents, and active power tracking for a load connected to a generator.}
    \label{fig:ac_load_gen}
\end{figure}

\subsection{Controllable Loads or Power Loads}
The \emph{controllable load} (CL) draws the required current from the grid to meet the desired active and reactive power demand. The dynamics of the CL in $dq$ coordinates is given as
\begin{equation}
    l_L\frac{d\mathbf{i}_L}{dt} = -\bigg({r_L}I-\omega_g J\bigg)\mathbf{i}_L(t) + \mathbf{v}_L(t) - \mathbf{v}_g(t),
\end{equation}
where $\mathbf{i}_L(t), \mathbf{v}_L(t), \mathbf{v}_g(t) \in \mathbb{R}^2$ are the load current, load voltage (control input) and the grid voltage. $r_L, l_L \in \mathbb{R}_+$ are the load resistance and inductance in \textsf{Ohm} and \textsf{Henry}. The objective of the power load is to draw the required current from the grid based on the demanded active and reactive power. The first step in achieving this control objective is to design the reference current given the active, reactive power profile and a grid voltage. Revisiting the definitions of the active and the reactive power
\begin{align*}
    P &= v_{g_d}i_{L_d} + v_{g_q}i_{L_q}, \\
    Q &= v_{g_d}i_{L_q} - i_{L_d}v_{L_d}.
\end{align*}
Rearranging the terms yields,
\begin{align*}
    \begin{bmatrix}
        v_{g_d} & v_{g_q} \\ -v_{g_q} & v_{g_d}
    \end{bmatrix} 
    \begin{bmatrix}
        i_{L_d} \\ i_{L_q}
    \end{bmatrix} = 
    \begin{bmatrix}
        P \\ Q
    \end{bmatrix}.
\end{align*}
Therefore, for some active and reactive power pair $(P,Q)$, the reference current is generated as follows
\begin{align*}
   \underbrace{\begin{bmatrix}
        i_{L_d} \\ i_{L_q} 
    \end{bmatrix}}_{\mathbf{i}_{{L}_{ref}}} = 
    \frac{1}{v_{g_d}^2 + v_{g_q}^2}
    \begin{bmatrix}
        v_{g_d} & -v_{g_q} \\ v_{g_q} & v_{g_d}
    \end{bmatrix}
    \begin{bmatrix}
        P \\ Q
    \end{bmatrix}.
\end{align*}

The current control for the power load is consequently given as
\begin{align*}
    \mathbf{v}_L = k_p \mathbf{\tilde{i}}_L + k_I \int_{0}^{\infty} \mathbf{\tilde{i}}_L (\tau) d\tau,
\end{align*}
where $\mathbf{\tilde{i}}_L = \mathbf{i}_L - \mathbf{i}_{{L}_{ref}}$. Figure \ref{fig:ac_load_gen} shows the performance (voltage, speed tracking), current, and active power tracking of the single load connected to a generator. The developed models can be extended to multiple generators and loads. 

\subsection{Energy Storage Elements}
The modeling of the Energy Storage Elements (ESEs) consists of the physics-based battery model, the State of Charge (SoC) and the State of Health (SoH) model.
\subsection{Model Development}
The battery model consists of multiple energy storage systems (ESSs) modeled as a current-controlled voltage source and a resistance $r_b$ coupled to the bus. The battery current $i_b$ is dictated by the controllable voltage source $v_b$. The open circuit voltage $v_{oc}$ is a function of the battery state of charge (SoC). The relation between $v_{oc}$ and SoC can be approximated using various functions ranging from lower to higher orders with the linear approximation being the fundamental one given as $v_{oc}=c_1SoC+c_2$, where $c_1,c_2$ are constants \cite{RR}\cite{Weng2013AnOM}. Since the battery voltage dynamics are fast compared to SoC, it is neglected, and, the control input $v_b \in \mathbb{R}$ is determined through
\begin{align*}
v_b = v_{bus}-\frac{p_{b_r}r_b-v_{bus}v_{oc}}{v_{bus}}, \hspace{-11mm} \\ 
i_b=\frac{v_{bus}-v_{b}-v_{oc}}{r_b},
\end{align*}
where $p_{b_r}$ is the optimal power profile from the upper-level controller.

\subsection{State of Charge and State of Health}
The dynamic model for the State of Charge (SoC) is given as follows:
\begin{equation}\label{SoC_dynamics}
    \dot{SoC} = -\frac{i_b}{Q} ,
\end{equation}
where $i_b \in \mathbb{R}$ denotes the battery current, $SoC \in \mathbb{R}$ the SoC vector, and $Q$ the total capacity of an individual battery in \textsf{ampere-hour}. The SoC dynamics in (\ref{SoC_dynamics}) is discretized using the Euler method for a time step of $T_s$ as:
\begin{equation}
     SoC_{k+1} = SoC_k-T_s \frac{i_b}{Q}.
\end{equation} 

The battery degradation model is based on the Arrhenius equation and uses the $Ah-$throughput $\displaystyle \int_{0}^{t}\left|i_b(\nu)\right|d\nu$ as a metric to evaluate the battery state of health (SoH), where $i_b(t) \in \mathbb{R}$ is the current drawn from the battery (positive in the discharge direction and vice-versa). The \emph{capacity loss} is given as \cite{SONG2018433}:
\begin{equation}
    Q_L(t) = \zeta_1 e^{\frac{-\zeta_2+TC_{r}}{RT}}\int_0^{t}|i_b(\nu)|^{\frac{1}{2}}d\nu,
\end{equation}
where $T \in \mathbb{R}_+$ is the baseline temperature (in \textsf{Kelvin}) of the battery where it is most efficient and $C_r \in \mathbb{R}_+$ is the C-rate of the battery.

\subsection{Flywheel Storage Modeling}
Flywheel storage technology consists of s spinning wheel connected to a motor that stores and releases energy through the rotational kinetic energy of the spinning wheel. During the charging phase the the motor spins the wheel leading to the storage of the energy. During the discharging phase, the wheel's rotational energy is converted back to electrical energy. The key advantage of flywheel storage as opposed to traditional energy storage is the low maintenance costs and its suitability for power-dense applications. However, the initial costs and safety concerns linger as key limitations. The basic energy stored in the flywheel (a solid disc) as rotational kinetic energy (KE) is given as
\begin{align*}
    E = \frac{1}{2}\bigg(\underbrace{\frac{1}{2}m r^2}_{I}\bigg)\omega^2, 
\end{align*}
where $r$ is the radius of the flywheel in \textsf{meters} and $m$ is the mass of the flywheel in \textsf{kg}, $\omega$ is the angular velocity of the flywheel in \textsf{rad/s} and $I$ is the flywheels moment of inertia. The dynamics of the flywheel are as follows
\begin{align}
    \frac{d\omega}{dt} = \frac{\tau}{I},
\end{align}
where $\tau$ is the applied torque to the flywheel. The power output of the flywheel is expressed as a dot product between torque and angular velocity
\begin{align*}
    p_f = \tau \omega.
\end{align*}
The state of charge for the flywheel is expressed as the current energy and the rated energy which in terms of the angular velocity and the rated angular velocity is
\begin{align*}
    SoC_f = \frac{\omega(t)^2}{\omega_{max}^2}
\end{align*}

Overall, the models for an AC power load, generator (DC and AC), and energy storage elements are presented in this chapter. The developed models will be used in subsequent chapters in modeling more detailed microgrids such as shipboard power systems and hybrid power trains.

\section{Degradation Aware Energy Management}

\subsection{Energy Storage Degradation Heuristics}
The presence of Pulse Power Loads (PPLs) in the Notional Shipboard Power System (SPS) presents a challenge in the form of meeting their high ramp rate requirements. Considering the ramp rate limitations on the generators, this might hinder the power flow in the grid. Failure to meet the ramp rate requirements might cause instability. Aggregating generators with energy storage elements usually addresses the ramp requirements while ensuring the power demand is achieved. This paper proposes an energy management strategy that adaptively splits the power demand between the generators and the batteries while simultaneously considering the battery degradation and the generator's efficient operation. Since it is challenging to incorporate the battery degradation model directly into the optimization problem due to its complex structure and the degradation time scale which is not practical for real-time implementation, two reasonable heuristics in terms of minimizing the absolute battery power and minimizing the battery state of charge are proposed and compared to manage the battery degradation. A model predictive energy management strategy is then developed to coordinate the power split considering the generator efficiency and minimizing the battery degradation based on the two heuristic approaches. The designed strategy is tested via a simulation of a lumped notional shipboard power system. The results show the impact of the battery degradation heuristics for energy management strategy in mitigating battery degradation and its health management.

\subsection{Background}
In Shipboard Power Systems (SPSs), power generation modules (PGMs) provide the power required by various loads through the direct current (DC) distribution system. Modern SPSs are equipped with advanced power loads such as rail guns, electromagnetic radars, and heavy nonlinear and pulsed power loads (PPLs), these advanced loads are high ramp rate loads \cite{Derry_1}\cite{Derry_2}. SPSs consist of two key power-supplying sources namely: Power Generation Modules and Power Conversion Modules (PCMs) and power consuming loads: Power Load Modules (PLM). PGMs are ramp rate limited which hinders in meeting the high ramp rate power required by the loads. Failure to ramp up in time to meet the load requirement might lead to system imbalances and instability. Adding additional generators is not a feasible solution considering the restrictions on the ship hull dimensions. Thus, the PCMs, which can support high ramp rates provide a solution in addressing the load requirement of high ramp rate loads \cite{Kim_2015}\cite{ESRDC_1}. The PCMs can operate bi-directionally, they can store the energy and provide immediate support when high ramp rate loads are requesting power. 

\begin{figure}[t!]
      \centering
      \includegraphics[width=0.6\textwidth]{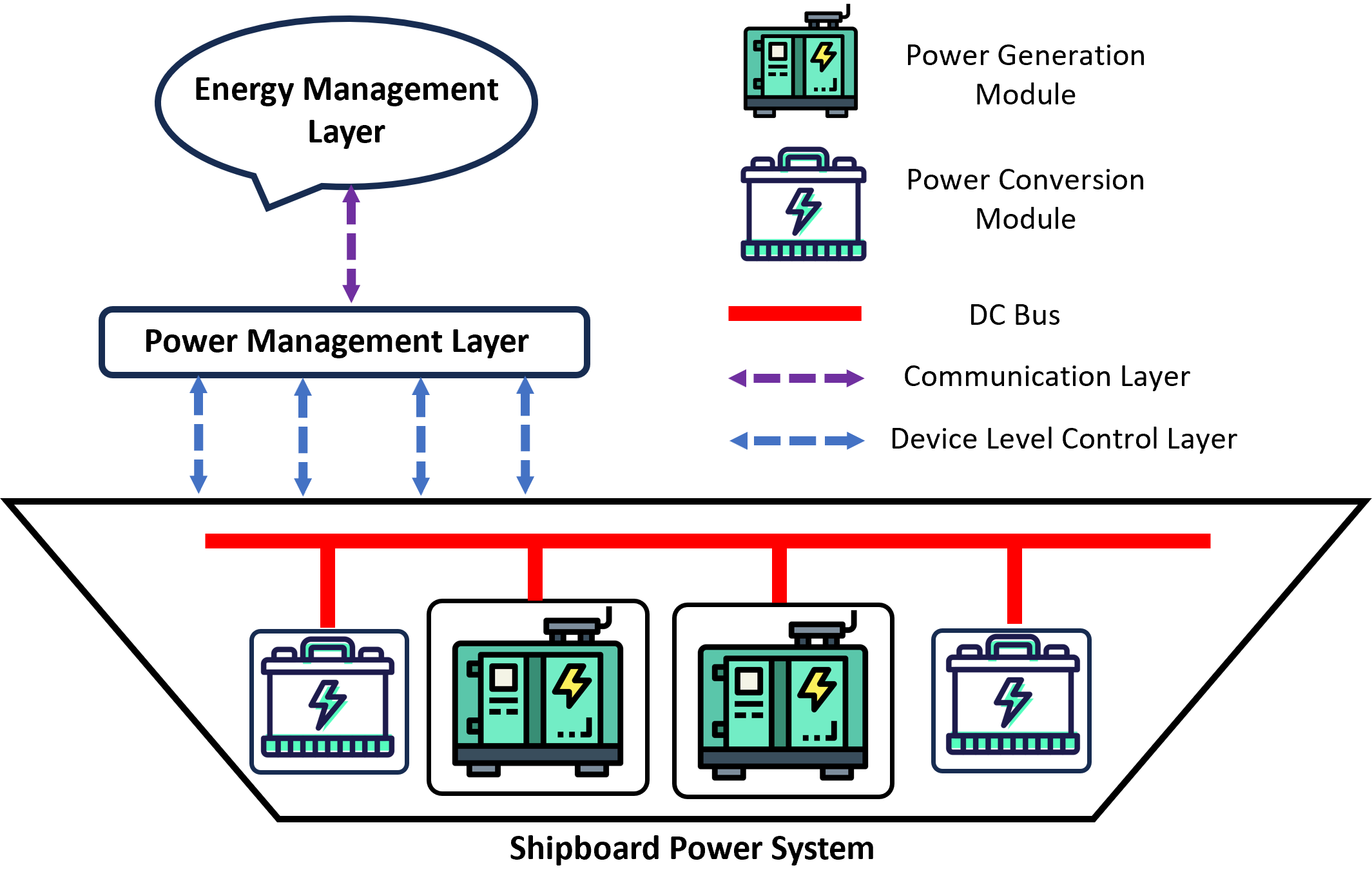}
	 \caption{An islanded notional shipboard power system}
     \label{4_Zone_SPS_2} 
\end{figure}

Existing literature considers the PGM and PCM ramp rate limitations in designing the energy management (EM) strategy \cite{2017_Vu_3,2021_Vedula,2017_Zohrabi,Zhang_2022}, which makes use of the optimal control theory (OCT) techniques such as dynamic programming (DP) and model predictive control (MPC). While limitations in the computational capabilities are the major roadblock in using the MPC over DP, it offers a better environment to incorporate constraints such as ramp rate limitations, and lower and upper power limitations on the PGMs and PCMs. Thus, due to recent advancements in the computational capabilities and its competence in handling the constraints, MPC, which originated as a control strategy for process control, made its way into the energy management of power systems.

Nonetheless, integrating PCMs into the existing SPS framework raises another important objective of managing the degradation of the PCMs. It is a well-known fact that the PCMs degrade faster than the traditional gas-turbine-based PGMs. Thus, the EM problem needs to address the \textit{PCM degradation, PGM efficient operation and the ramp rate limitations} together. Since PCM degradation is a complex process that is difficult to capture mathematically with high precision, investigators have come up with a few degradation minimization measure heuristics embedded into the cost function. Minimizing these heuristics impacts the degradation process. Table-\ref{tab:references} shows the various degradation measures used as heuristics in the existing literature. 

\begin{table}[ht]
\centering
\caption{PCM degradation measure heuristics and references}
\label{tab:references}
\resizebox{0.5\textwidth}{!}{\begin{tabular}{c|c}
\hline \hline
\textbf{Degradation Measure }    & \textbf{References}  \\
 \hline \hline
Depth of Discharge (DoD)    & \cite{Hein_2021}\cite{Steen_2021}\cite{Li_2023}\cite{Zhao_2023}   \\
State of Charge (SoC)   & \cite{Nawaz_2023}\cite{Wang_2022}\cite{Vedula_2} \\
Charge/Discharge Cycle  & \cite{Ji_2023} \\
\hline
\end{tabular}}
\end{table}

While, \cite{Hein_2021,Steen_2021,Li_2023,Zhao_2023,Nawaz_2023,Ji_2023,Wang_2022} have used different measures for capturing the degradation process, we consider the \emph{PCM power, PCM state of charge and PGM power} together in an optimization-based degradation measure heuristic framework to analyze and study the PCM degradation process. A centralized model-predictive energy management problem is formulated. The impact of the PCM degradation measurement heuristics chosen with objectives: minimizing the \textit{absolute power extracted} out of the PCM and minimizing the \textit{state of charge} of the PCM. The results compare which of the two heuristics better aids in mitigating the PCM capacity loss \cite{vedula2024batterydegradationheuristicspredictive}.

\subsection{Preliminaries: Notional Shipboard Power System}
 \begin{figure}[t!] 
	\centering
	\includegraphics[width=0.7\textwidth]{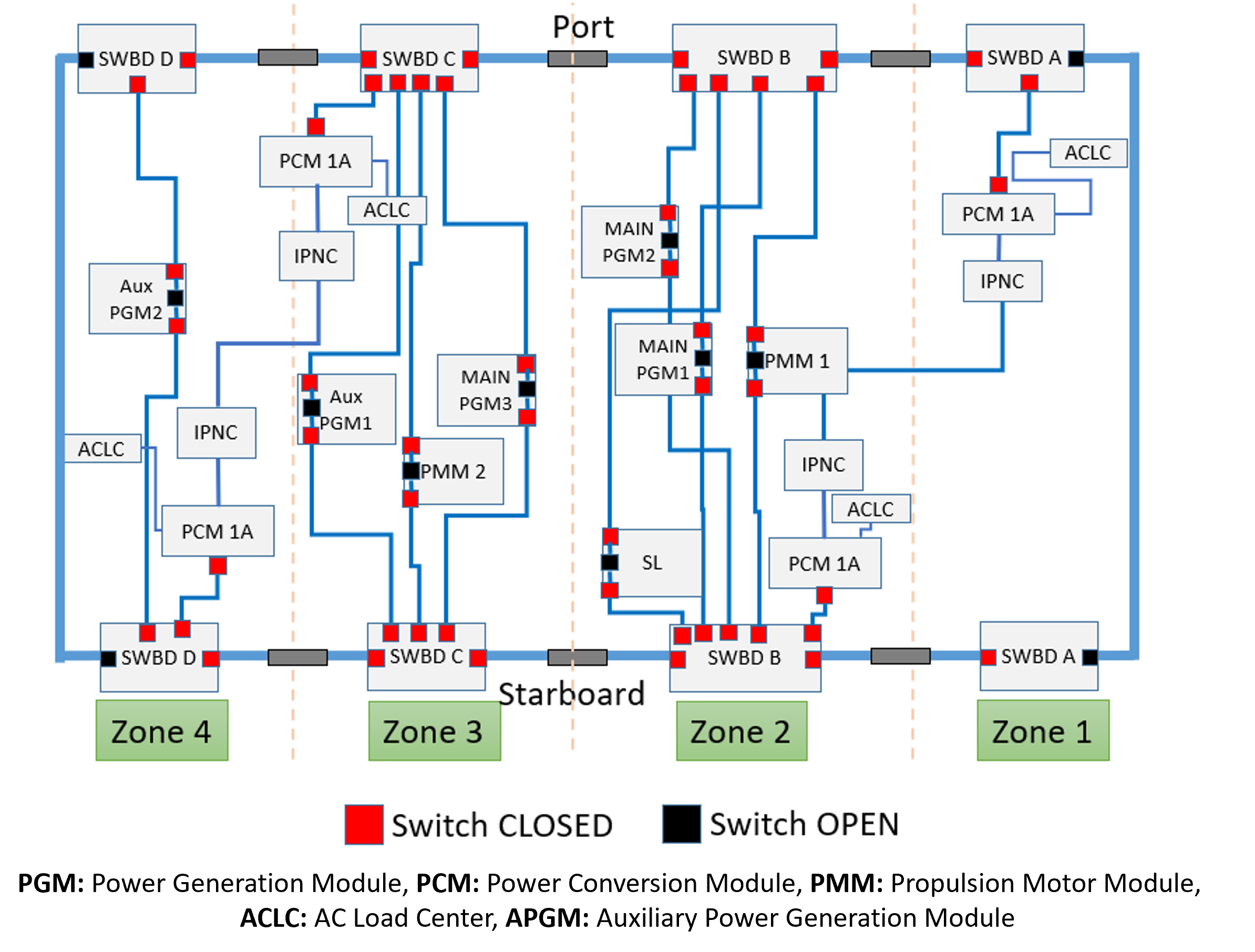} 
	\caption{Notional 4-zone shipboard power system model presented by the Office of Naval Research \cite{ESRDC_1}. 
    }
	\label{SPS_Model}
\end{figure}
The shipboard power system model proposed by the U.S. Office of Naval Research (ONR) shown in Figure \ref{SPS_Model} comprises PGMs, PCMs, and device level controllers (DLCs) for local voltage and current control. These components are categorized into a zonal structure \cite{ESRDC_1}. PGMs consist of fuel-operated gen-sets. PCMs consist of multiple energy storage elements such as ultracapacitors and battery energy storage systems. All zones are unified through a common 12\textsf{kV} DC bus. The power demand must be met under any given circumstance i.e. power supplied by all sources must match the load demand scenarios such as high ramp and uncertainty in the load to maintain the system stability. In this paper, only the main power-supplying sources as mentioned before are considered. 
There are two main PGMs and two main PCMs. The combined power supplied by these sources must meet the SPS power demand. The following notations are used in the Shipboard Power System context:
\begin{itemize}
    \item PGM: Power Generation Module. It is also referred to as a generator in some parts of the paper. The usage of PGM and generator are synonymous.
    \item PCM: Power Conversion Module. It is also referred to as an energy storage system (ESS) or battery in some parts of the paper. The usage of PCM/ESS/Battery is synonymous.
    \item The assumption is that the ESS operates bidirectionally representing charging and discharging operations. Also, we assume that the ESS can ramp up faster than the generators.
    \item PLM: Power Load Module.
\end{itemize}

The concept of MPC is based on the optimal control theory and the \textit{Bellman's principle of optimality}. \textit{``Given the optimal sequence $u^*(0),\hdots,u^*(N-1)$ and the corresponding optimal trajectory $x^*(0),\hdots,x^*(N-1)$, the sub-sequence $u^*(k),\hdots,u^*(N-1)$ is optimal for the sub-problem in the horizon $[k,N-1]$ starting from optimal state $x^*(k)$ \cite{borrelli_bemporad_morari_2017}"}. 

\subsubsection{SPS Model Development}
We consider the lumped PGM and the lumped PCM models for the model development. Thus, the total number of the components in an SPS is reduced to a single PGM, a single PCM supplying a common single PLM. 

\subsubsection{Power Generation Module Modeling}
The PGM model used in this work is modeled as a controllable voltage source in series with a \textsf{RL}-line and a shunt capacitance. The mathematical model for the individual PGM is given as follows \cite{2021_Vedula}:
\begin{align}
    l_g \frac{di_g}{dt} &= -r_g i_g(t) + v_g - v_c, \\
    c_g \frac{dv_c}{dt} &= i_g - i_{g_{ref}},
\end{align}
where $i_g \in \mathbb{R}$ is the PGM current, $l_g, r_g, c_g\in \mathbb{R}_+$ are the inductance, resistance, and the capacitance of the PGM. The capacitance voltage/bus voltage $v_{c} \in \mathbb{R}$ and the controlled generator voltage $v_g \in \mathbb{R}$. Given a power reference ${p}_{g_{ref}} \in \mathbb{R}$ from the high-level control, the current reference for the device level PGM current controller (DLC) is given as ${i}_{g_{ref}} = {p}_{g_{ref}}/{v}_c \in \mathbb{R}$. The goal of the DLC is to achieve the objective:
$$\int_0^{\infty} \bigg(\underbrace{v_c(t)-v_{g_{ref}}(t)}_{\tilde{v}_g}\bigg)^2 dt < \infty.$$
Since the design and analysis of such a DLC is not the main objective of this paper numerous linear and nonlinear methods can be employed to design it \cite{Khalil_Book}. Thus, a simple \emph{PI} controller is employed of the form:
$$v_g = v_c + K_P \tilde{v}_g(t)+K_I \int_0^{t}\tilde{v}_g(\tau)d\tau,$$
where $K_P,K_I \in \mathbb{R}$ are the device level current control gains.

\subsubsection{Power Conversion Module Modeling}
The PCM is modeled as an ESS which comprises a single or multiple hybrid energy storage systems. This might include a collection of flywheels, supercapacitors, or battery energy storage systems (BESS). The dynamics of grid-following PCM used in this work and the battery current calculations are based on \cite{2021_Vedula}. The relationship between the PCM power injected $p_{b}$ and the SoC for the individual PCM is given as follows:
\begin{equation}\label{SoC_Power_2}
    \dot{q} = -\frac{p_b}{Q_bv_c}
\end{equation}
where $q \in \mathbb{R}$ denotes the State of Charge (SoC), ${Q}_{b} \in \mathbb{R}_+$ is the total capacity of the PCM in \textsf{AHr}. The bus voltage to which the PCM is coupled is given by $v_{c} \in \mathbb{R}$. The initial SoC is denoted by $q_0 \in [0,1]$. For this work, the SoC dynamics in (\ref{SoC_Power_2}) is discretized with a sampling time of ${T_s} \in \mathbb{R}_+$ using the Euler method, and the following discrete dynamics are obtained:

\begin{equation}\label{SoC_Power_Disctrete_sps}
    q_{t+1} = q_t - \frac{T_s}{Q_b v_c}p_{b_t}.
\end{equation}

The PCM degradation model used in this work is an Arrhenius equation-based model which uses $Ah$-throughput as $\displaystyle \int_{0}^{t}\left|{i}_b(\nu)\right|d\nu$ as a metric to evaluate battery state of health (SoH). ${i}_b(t) \in \mathbb{R}$ is the current drawn from the PCM (positive while discharging and negative while charging). The PCM \textit{capacity loss} formulation is given as follows \cite{SONG2018433}.
\begin{equation}
    Q_L(t) =  e^{\frac{-\zeta_1+T_bC_{r}}{RT_b}}\int_{0}^{t}\left|{i}_b(\nu)\right|d\nu
\end{equation}
where $T_b$ is the PCM temperature in \textsf{Kelvin}, $C_r$ is the C-rate of the PCM. The capacity loss (\%) is given as follows: $$\Delta Q \% = \frac{Q_b-Q_{L}(t)}{Q_b}\times 100,$$
where $Q_{L}$ is the \emph{capacity loss} of the PCM in \textsf{ampere-hour}. \textit{This does not capture the battery end of life (EOL)}. It only captures the capacity lost during the PCM operation ($\Delta Q$).

\subsubsection{Power Load Module Modeling}
The PLM is a controllable load and is modeled as a current sink regulated through a controllable voltage. The dynamics are given as follows:
\begin{equation}
    l_L \frac{di_L}{dt} = -r_L i_L + \underbrace{v_L - v_c}_{\tilde{v}}.
\end{equation}

The goal of the PLM-DLC is to achieve the objective:
$$\int_0^{\infty} \bigg(\underbrace{i_L(t)-i_{L_{ref}}(t)}_{\tilde{i}_L}\bigg)^2 dt < \infty.$$
Given a load power $p_L$, the load current reference can be determined as $i_{L_{ref}} = p_L/v_c$. A simple \emph{PI} controller is employed of the form:
$$\tilde{v} = K_P \tilde{i}_g(t)+K_I \int_0^{t}\tilde{i}_g(\tau)d\tau,$$
where $K_{P_L},K_{I_L }\in \mathbb{R}$ are the PLM device level current control gains. 

Thus, the power flow in the SPS is given as follows:
\begin{equation}\label{SPS_Power_Flow}
    {{p}}_{g}(t)+{p}_{{b}}(t) - {p}_{{L}}(t) = 0,
\end{equation}
where ${p}_{g_i}, {p}_{b_i}, {p}_L\in\mathbb{R}$ are powers corresponding to the PGMs, the PCMs, and the PLMs respectively.

\begin{figure}[t!] 
	\centering
	\includegraphics[width=0.5\textwidth]{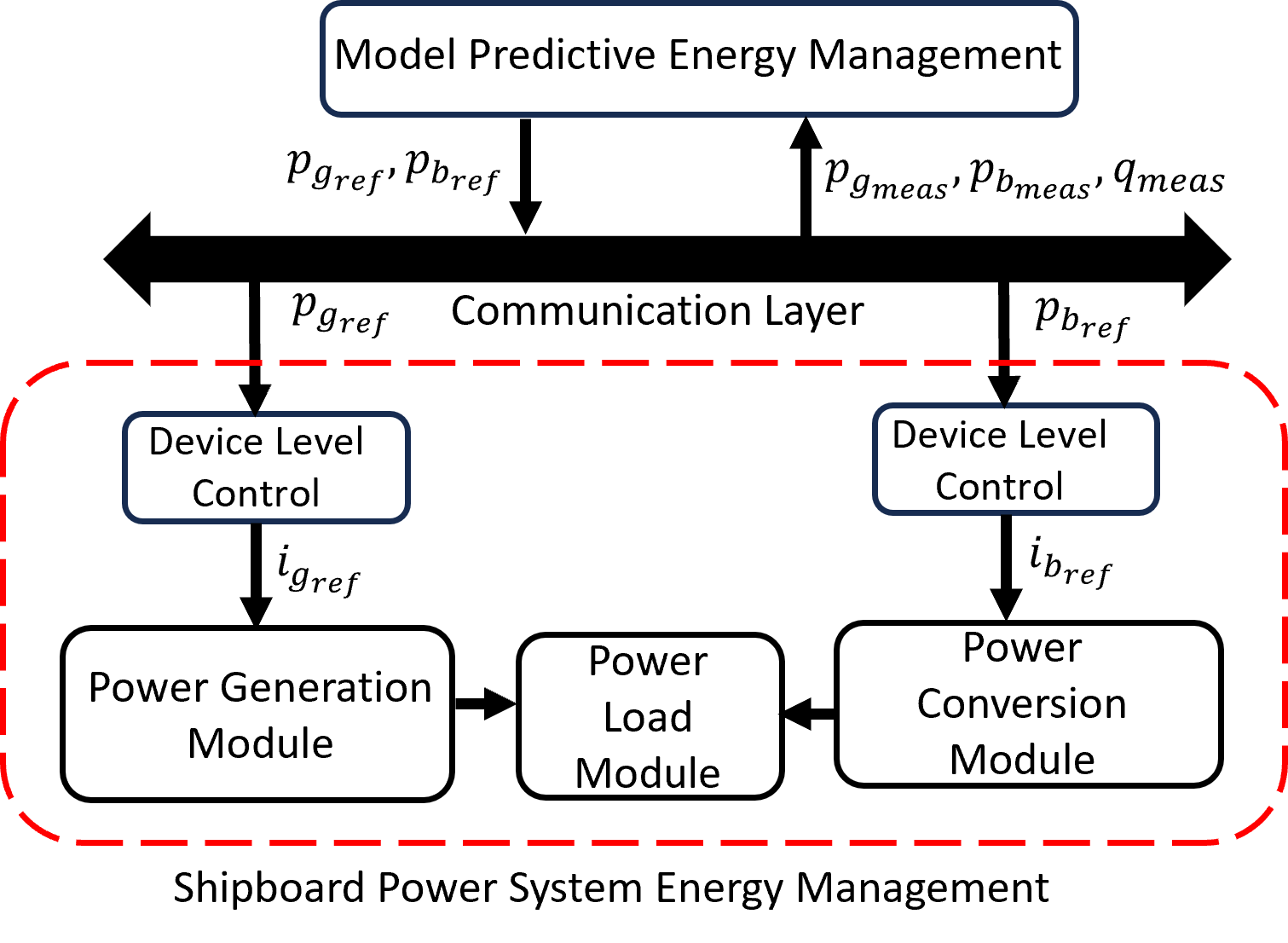} 
	\caption{Hierarchical control structure employed for the energy management in SPS}
	\label{EM_SPS_pic}
\end{figure}

\subsection{Heuristics-Based Model Predictive Energy Management}\label{sec: control development}
In this section, the optimization constraints imposed are developed followed by the model predictive energy management optimization problems considering \textit{battery power} and \textit{battery state of charge} as the PCM degradation heuristic measures.
\subsubsection{Constraints Development}
\textit{Ramp Limitations:}
Since the PGM is ramp-limited and the PCM has faster ramping capabilities compared to the PGM, the following constraints are imposed limiting the power-delivering capabilities in consequent time instances in a horizon ($H$).
\begin{subequations}\label{ramp-limitations}
    \begin{align}
        \left|\mathbf{p}_{g_k} - \mathbf{p}_{g_{k-1}}\right| \preceq r_g, \hspace{2mm} \forall k=1,\hdots,H, \\
        \left|\mathbf{p}_{b_k} - \mathbf{p}_{b_{k-1}}\right| \preceq r_b, \hspace{2mm} \forall k=1,\hdots,H,
    \end{align}
\end{subequations}
where $r_b > r_g$ and $r_g, r_b \in \mathbb{R}$ are the ramp-rate limitations of the PGM and the PCM.

\textit{Safety Constraints:} In order to accommodate for the safe operation of the PGM and the PCM, the lower and upper power limitations, SoC are imposed as
\begin{subequations}\label{safety-cons}
    \begin{align}
        \underline{\mathbf{p}}_{g} \preceq \mathbf{p}_{g_k} \preceq \overline{\mathbf{p}}_g, \hspace{2mm} \forall k=1,\hdots,H, \\
        \underline{\mathbf{p}}_{b} \preceq \mathbf{p}_{b_k} \preceq \overline{\mathbf{p}}_b, \hspace{2mm} \forall k=1,\hdots,H, \\
        \underline{\mathbf{q}} \preceq \mathbf{q}_k \preceq \overline{\mathbf{q}}, \hspace{2mm} \forall k=1,\hdots,H.
    \end{align}
\end{subequations}

\textit{SoC Constraint:} The discretized SoC dynamics in (\ref{SoC_Power_Disctrete_sps}) are imposed as an equality constraint to be satisfied at every instant of the horizon:
\begin{equation}\label{soc-balance}
    \mathbf{q}_{k+1} = \mathbf{q}_k - \frac{T_s}{Q_b v_c}\mathbf{p}_{b_k}, \hspace{2mm} \forall k=1,\hdots,H.
\end{equation}

\textit{Power Balance Constraint:} The power supply and the demand must be satisfied in a grid at any given time instant. To that extent, the following constraint is imposed based on (\ref{SPS_Power_Flow}). It is assumed that the load power measured is held constant for the length of the optimization horizon interval.
\begin{equation}\label{power-balance}
    \mathbf{p}_{g_k} + \mathbf{p}_{b_k} - p_L\mathbf{1} = 0, \hspace{2mm} \forall k=1,\hdots,H.
\end{equation}

\subsubsection{Heuristics-Based Energy Management Problem}
The objective/cost of the energy management problem consists of three optimization variables: PGM power, PCM power, and SoC. The cost is formulated to minimize the PGM power injection, PCM power injection, and PCM SoC. Since the device level control is designed and is tasked to track the power set points provided by the MPC layer, the MPC problem and the imposed constraints are given as
\begin{equation}\label{MPEM_prb}
    \begin{aligned}
        \Minimize_{\mathbf{p}_g,\mathbf{p}_b,\mathbf{q}} \quad & \frac{\beta}{2}\norm{\mathbf{p}_{g}-\mathbf{p}_{g_r}}_2^2 + \frac{\gamma_p}{2}\norm{\mathbf{p}_{b}}_2^2 + \frac{\gamma_q}{2}\norm{\mathbf{q}-\mathbf{q}_0}_2^2 \\
        \SubjectTo \quad & (\ref{ramp-limitations})-(\ref{power-balance}),
    \end{aligned}
\end{equation}
where $\mathbf{p}_{g_r} \in \mathbb{R}^H$ is the known desired operating point of the PGM. $\mathbf{q}_0 \in \mathbb{R}^H$ is the initial SoC of the PCM. The objective is to maintain the PGM around a desired operating point and to maintain the SoC as close as possible to the initial SoC. The constraint set is a \textit{polytope} comprised of an affine equality constraint and an affine inequality constraints. The objective function is strongly convex. Based on this fact the optimization problem in (\ref{MPEM_prb}) is feasible, the optimization problem is convex whose global optimal solution can be found in polynomial time by existing algorithms \cite{boyd_vandenberghe_2004}.

\subsection{Numerical Simulation}\label{sec: simulation}
\begin{figure}[t!] 
	\centering
	\includegraphics[width=0.6\textwidth]{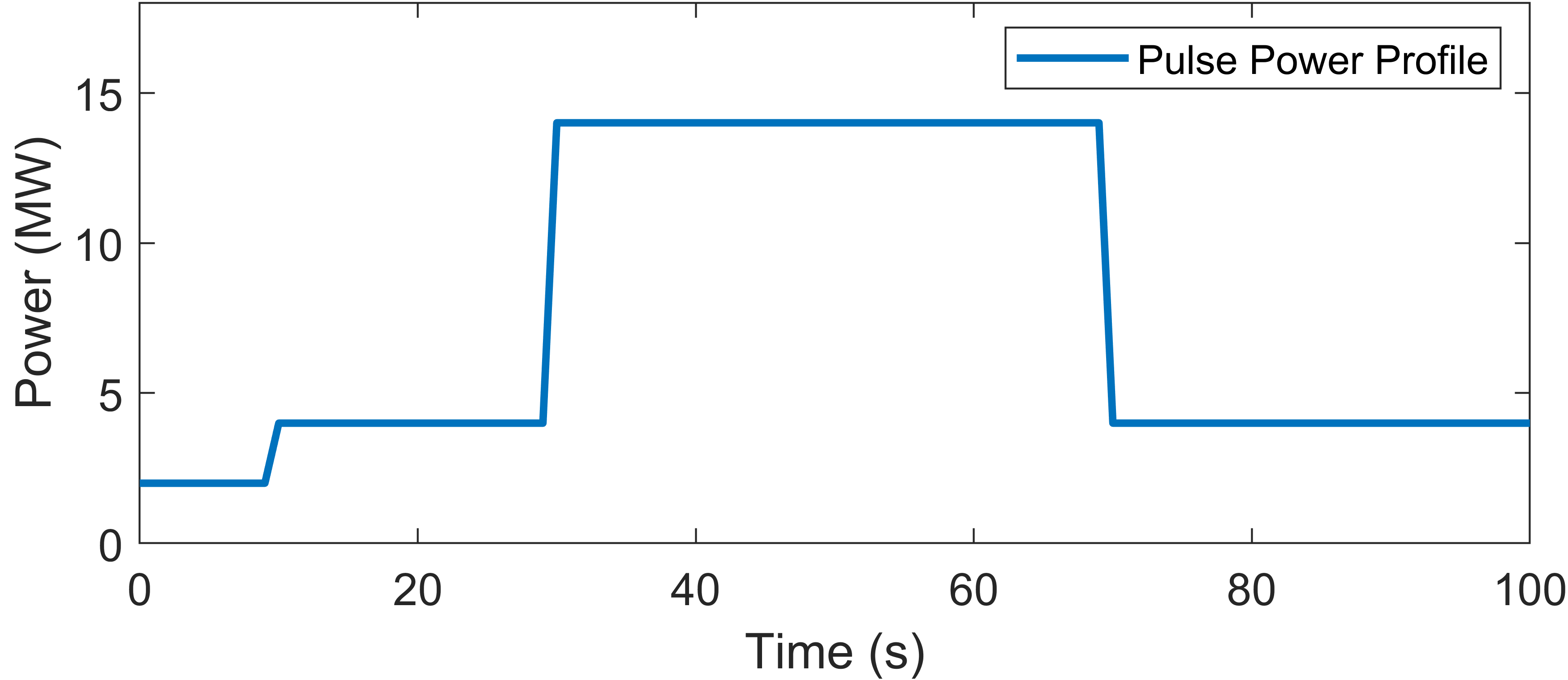} 
	\caption{Power profile}
	\label{Power_Profile}
\end{figure}
A lumped model of the SPS consisting of a single PGM, a single PCM, and a single PLM is considered for simulation purposes. The PGM and PCM combined to supply the PLM. The developed model and the optimization problems are implemented on the Digital Storm desktop with the processor Intel(R) Core(TM) i9-14900K 3.20GHz and 64GB of RAM.  The simulation time-step was chosen to be 1\textsf{ms}. The time taken by the algorithm to converge was observed to be 0.04\textsf{seconds} and about 50 iterations for the optimal values to be within the tolerance limits. Thus, the MPC and model rate transition in Simulink simulation was set at 1\textsf{second} i.e. the MPC updates every second. The MPC was implemented using YALMIP \cite{lofberg}. The power profile is designed based on the PLM ratings and an occurrence of a pulse is considered from 20-70sec. The simulation parameters are shown in Table \ref{tab:rated_2}. Figure \ref{Power_Profile} shows the designed PLM load profile. Three simulation scenarios are considered based on the optimization problem (\ref{MPEM_prb}).
\begin{itemize}
    \item \textit{No PCM Heuristic (Scenario 1):} $\beta = 1$, $\gamma_p=0$, and $\gamma_q = 0$.
    \item \textit{Power Minimization Heuristic (Scenario 2):} $\beta = 1$, $\gamma_p=1000$, and $\gamma_q = 0$.
    \item \textit{SoC Minimization Heuristic (Scenario 3):} $\beta = 1$, $\gamma_p=0$, and $\gamma_q = 1000$.
\end{itemize}
The choice of the penalty parameters $\beta$, $\gamma_p$, and $\gamma_q$ is chosen based on the observations from the simulations undertaken with different penalty values. The chosen values for the penalty best capture the results of the optimization problem for different scenarios. 

\begin{figure}[h!] 
	\centering
	\includegraphics[width=0.6\textwidth]{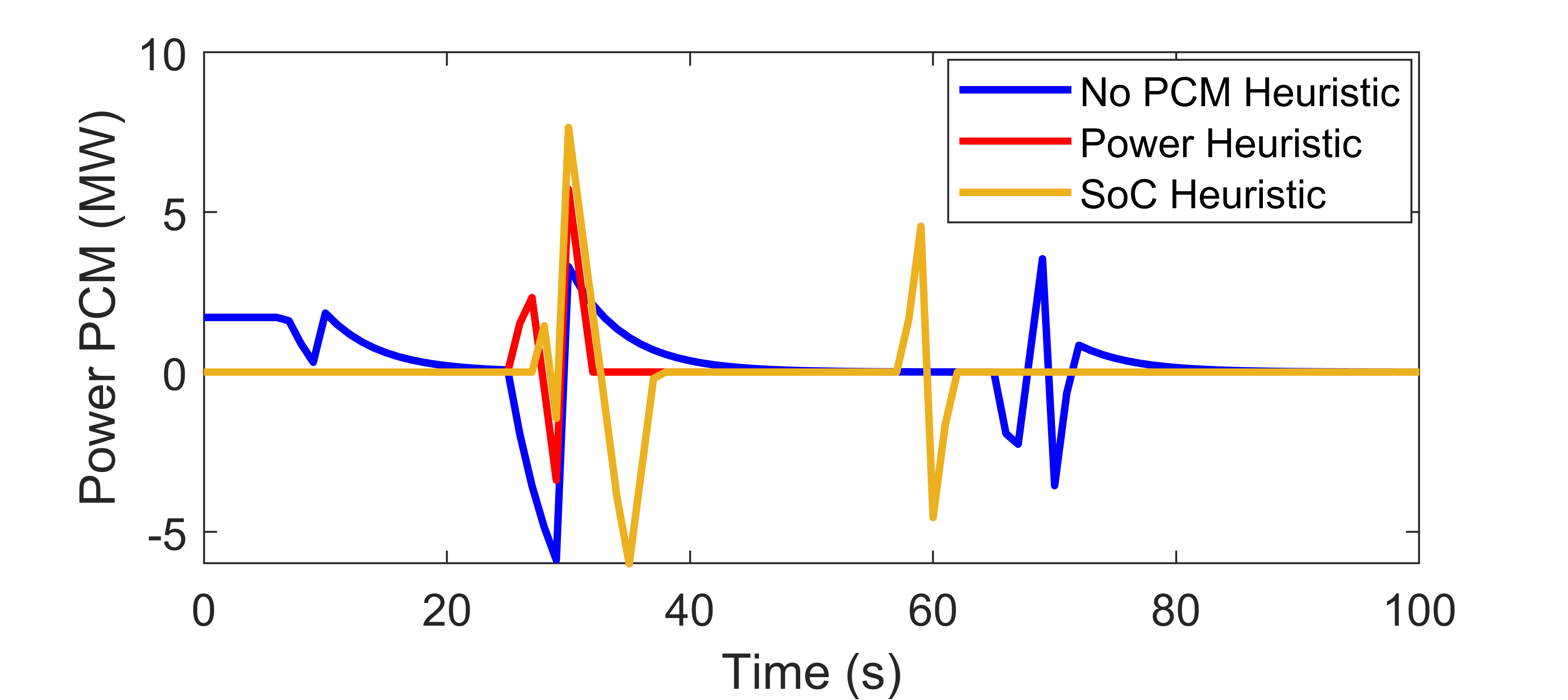} 
	\caption{PCM powers for different heuristic scenarios}
	\label{Power_PCM}
\end{figure}

Figure \ref{Power_PCM} shows the PCM powers with different scenarios implemented. It can be observed that Scenarios 2 and 3 display almost similar trends in terms of power injected by the PCM in response to the requested load profile. Scenario 2 and 3 which penalizes the power extracted and the SoC deviation from the PCM provides the most conservative response to the requested load. 

\begin{table}[ht]
\centering
\caption{Rated values and simulation parameters \cite{ESRDC_1}}
\label{tab:rated_2}
\resizebox{0.5\textwidth}{!}{\begin{tabular}{c|c|c}
\hline \hline
\textbf{Parameter}  & \textbf{Parameter}  & \textbf{Parameter}  \\
\textbf{Description} & \textbf{Notation}   & \textbf{Value}  \\  \hline \hline
PGM Upper Power Limit   & $\overline{p}_g$   & 28 MW \\ 
PGM Lower Power Limit   & $\underline{p}_g$   & 0.2 MW \\ 
PCM Upper Power Limit & $\overline{p}_b$ & 10 MW \\
PCM Lower Power Limit & $\underline{p}_b$ & -10 MW \\
SoC Lower Limit & $\underline{q}_b$ & 0.7 \\ 
SoC Upper Limit & $\overline{q}_b$ & 0.8 \\
Horizon & $H$ & 5 sec \\ 
PCM Capacity (Li-ion) & $Q_b$   & 20 AHr \\ 
PCM Ramp-rate (Li-ion) & $r_b$   & $\overline{p}_b$ \\
PGM Ramp-rate  & $r_g$   & 0.1$\overline{p}_g$ \\\hline
\end{tabular}}
\end{table}

\begin{figure}[h!] 
	\centering
	\includegraphics[width=0.6\textwidth]{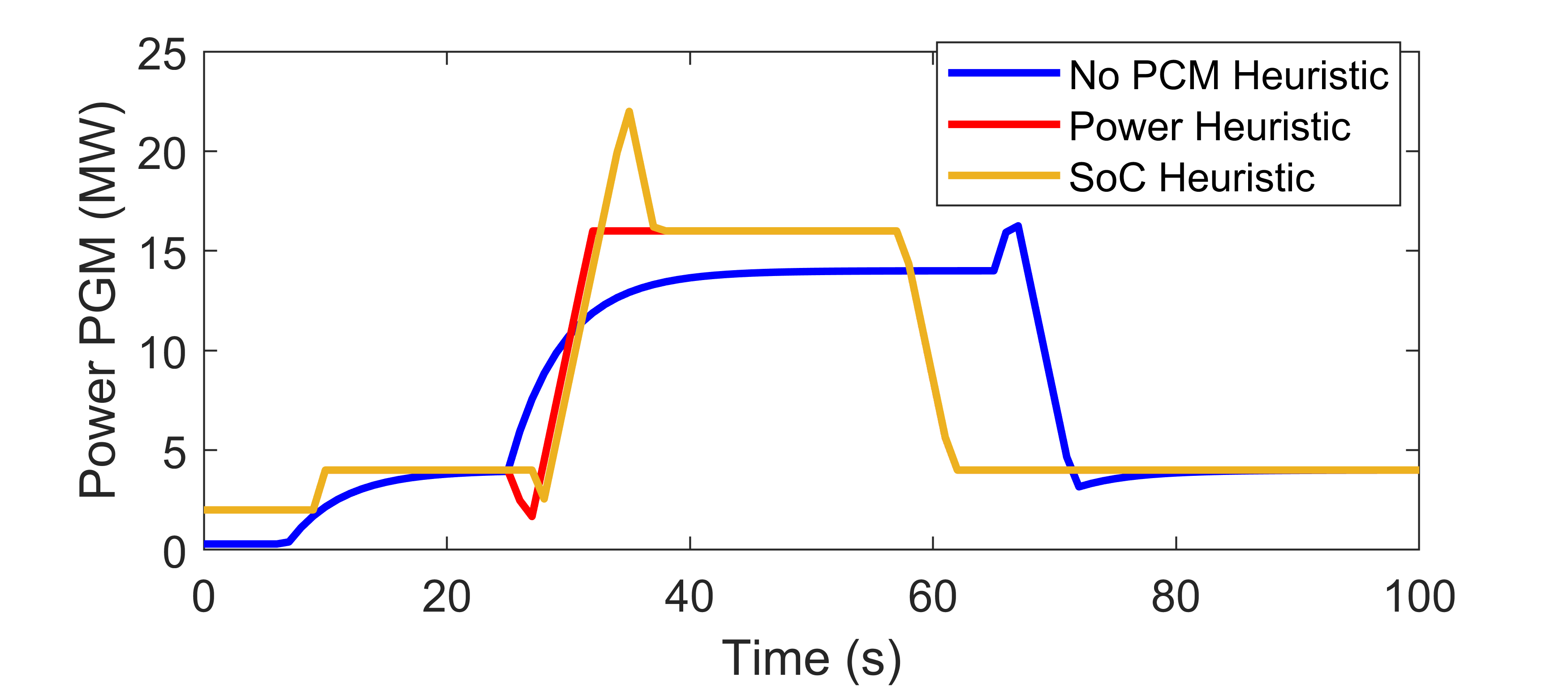} 
	\caption{PGM powers for different heuristic scenarios}
	\label{Power_PGM}
\end{figure}

The impact of the optimization scenarios and the PCM minimization heuristics on the PGM can be seen in Figure \ref{Power_PGM}. The PGM is pushed to its ramping and operational limits when the PCM power and the SoC are minimized i.e. Scenario2 and Scenario 3 as compared to Scenarios 1. The PGM acts in response to maintain the power balance constraint.

\begin{figure}[h!] 
	\centering
	\includegraphics[width=0.6\textwidth]{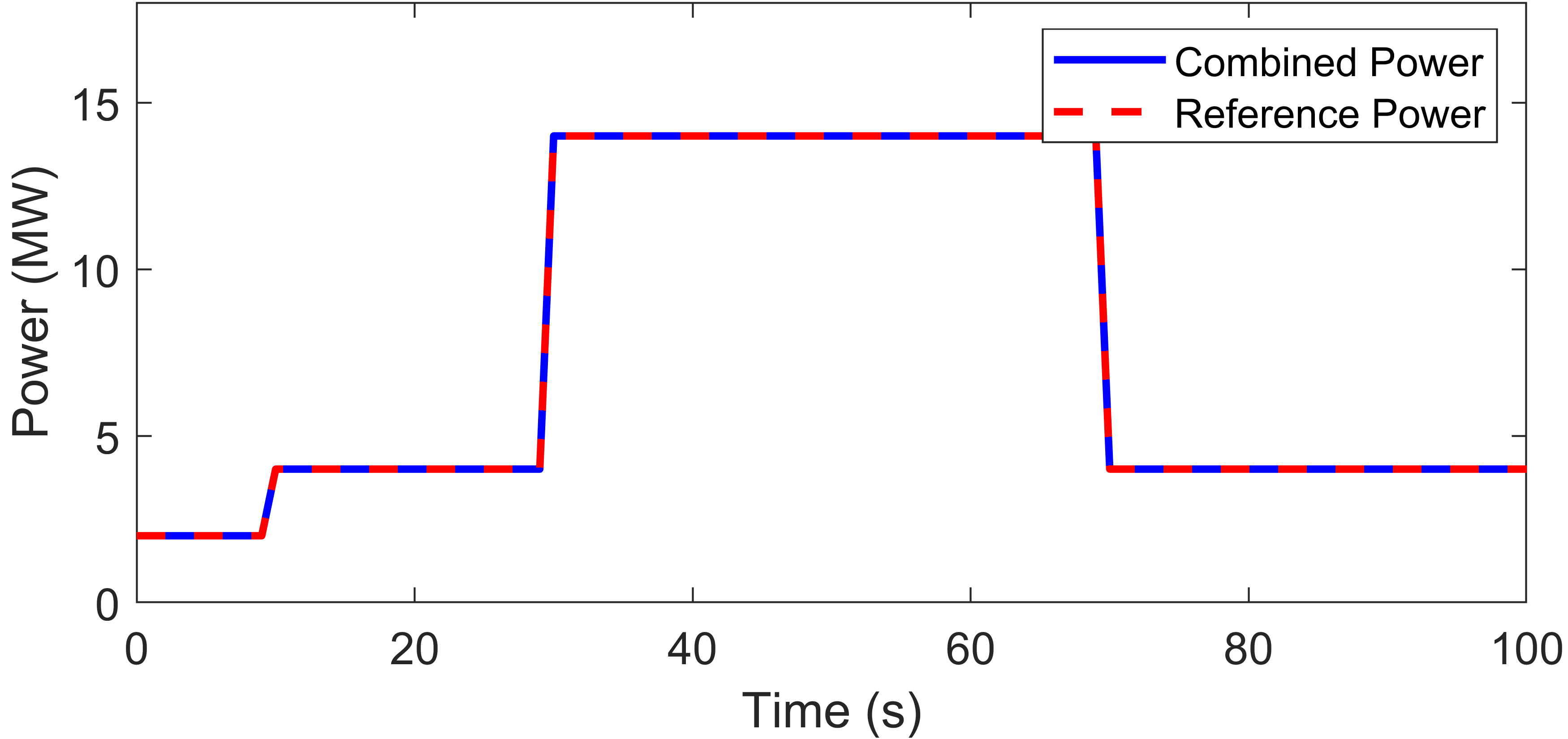} 
	\caption{Power tracking response for scenario 2}
	\label{Power_Tracking}
\end{figure}

The combined powers injected by the PGM and the PCM in response to the designed pulse power profile are shown in Figure \ref{Power_Tracking}. In order to avoid redundancy, the tracking response is plotted for Scenario 2. A similar response in tracking was also observed for other optimization scenarios.

\begin{figure}[h!] 
	\centering
	\includegraphics[width=0.6\textwidth]{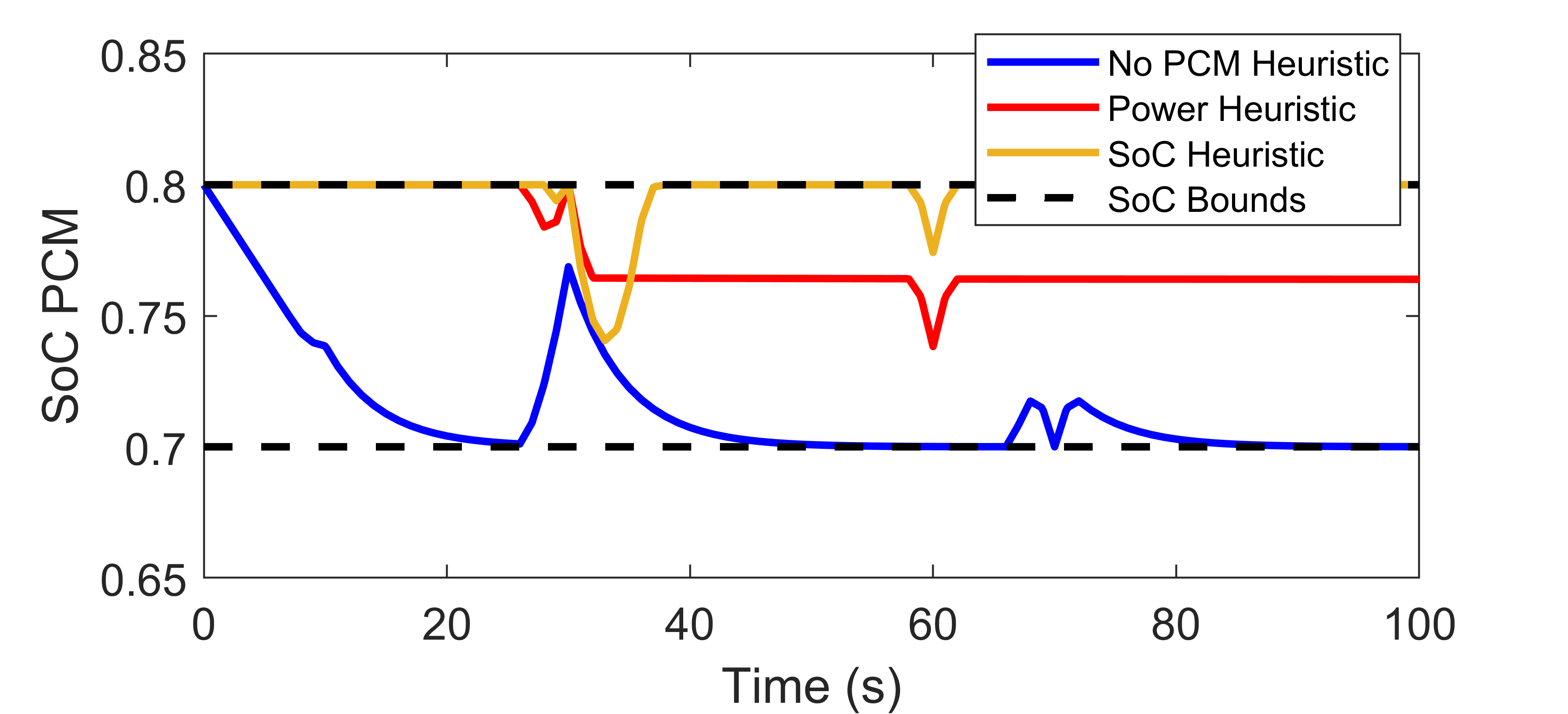} 
	\caption{SoC of the PCM for different heuristic scenarios}
	\label{PCM_SOC}
\end{figure}

Figure \ref{PCM_SOC} shows the State of Charge response for the proposed optimization scenarios. It is observed that all the scenarios obey the imposed upper and lower constraints on the SoC. Scenario-3 whose objective is to minimize the SoC, tries to keep the SoC as close as possible to the initial SoC. The SoC that is closest to the initial SoC can be observed during Scenario 3. 

\begin{figure}[h!] 
	\centering
	\includegraphics[width=0.6\textwidth]{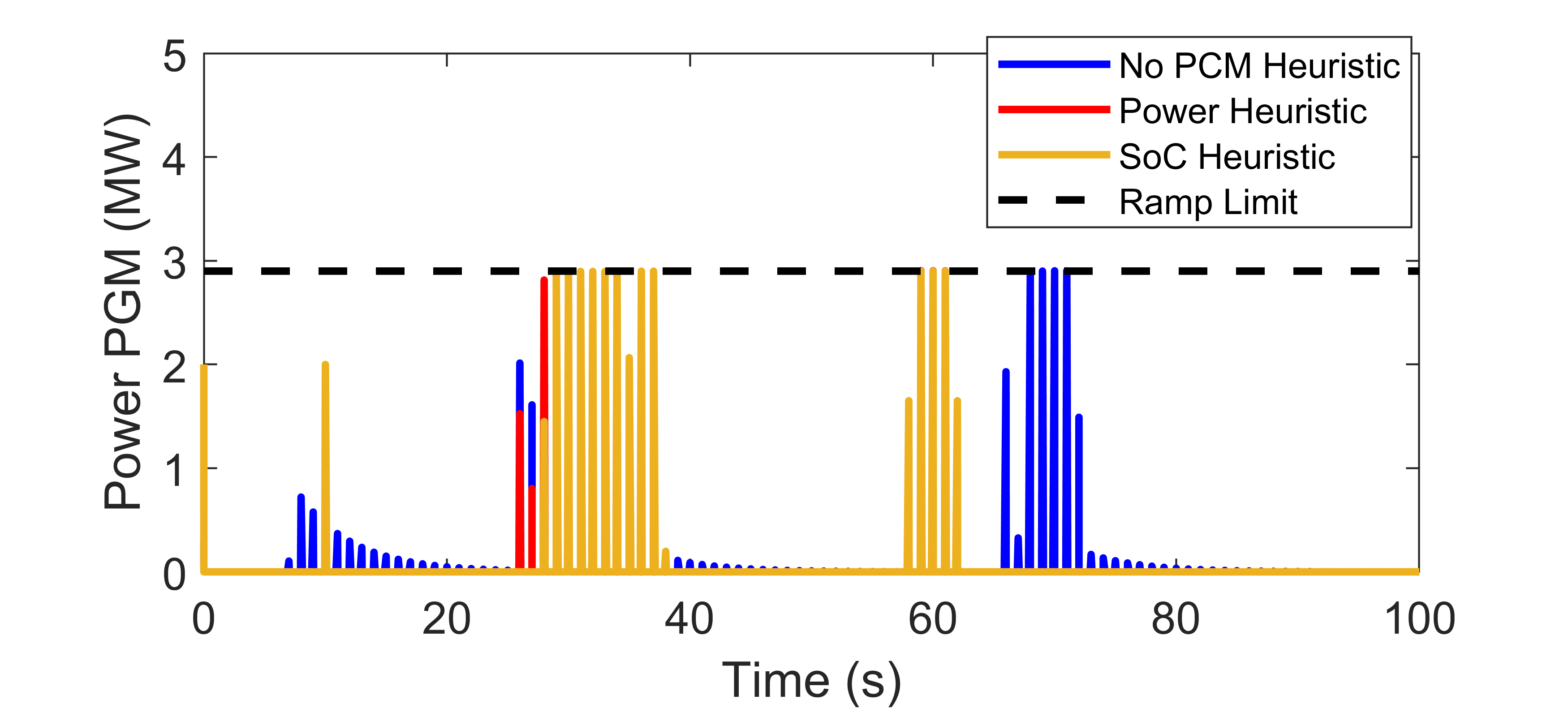} 
	\caption{PGM ramp limits: PGM power for every horizon step}
	\label{PGM_ramp}
\end{figure}

\begin{figure}[t!] 
	\centering
	\includegraphics[width=0.6\textwidth]{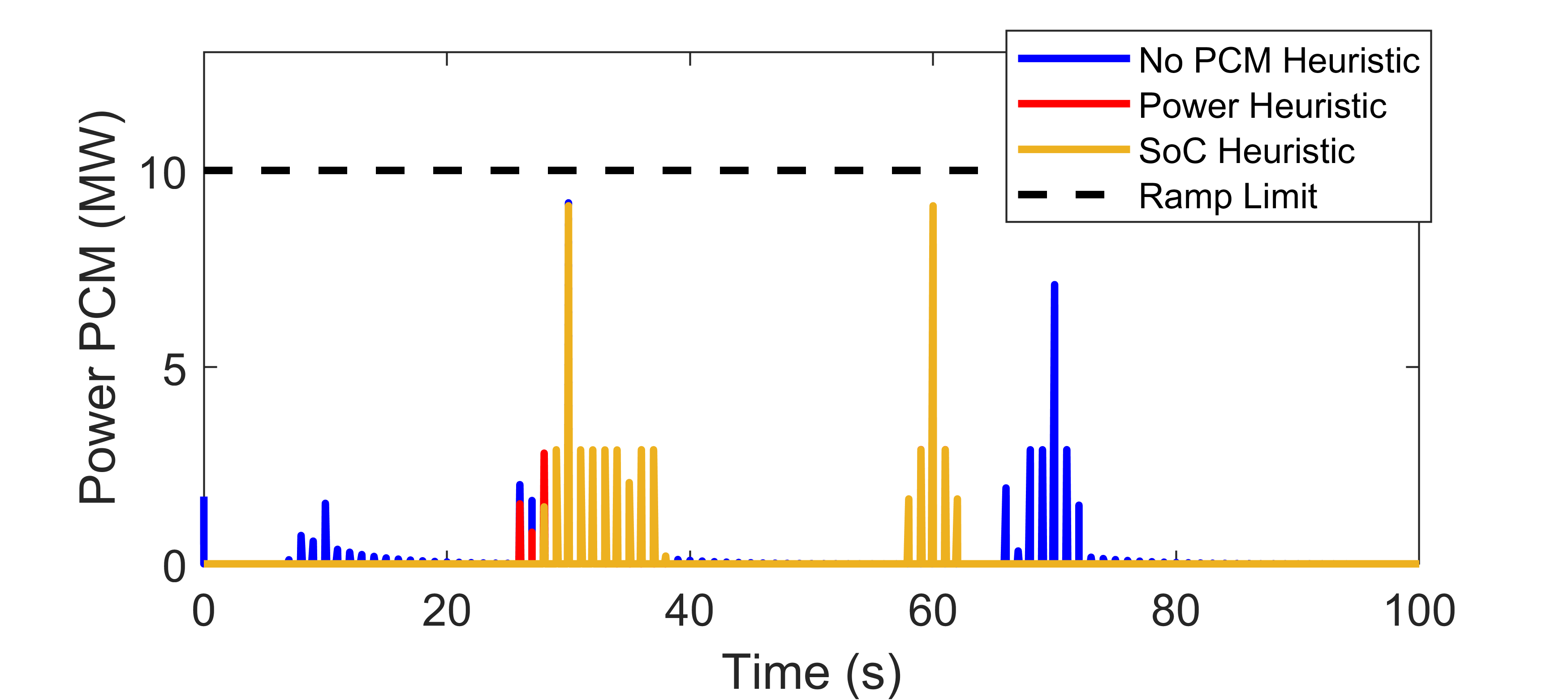} 
	\caption{PCM ramp limits: PCM power for every horizon step}
	\label{PCM_ramp}
\end{figure}

The ramp limits for the PGM and the PCM are shown in Figure \ref{PGM_ramp} and Figure \ref{PCM_ramp}. It can be seen that Scenario 2 and Scenario 3 pushes PGM to the maximum ramp limit, while for the PCM the same scenario has the lowest ramping. This indicates a more conservative PCM usage at the expense of aggressive PGM usage in order to satisfy the power balance constraint. Another key observation during the ramping of PGM is around 70sec, where both Scenarios 2 and 3 exhibit similar PGM power injections. 

\begin{figure}[t!] 
	\centering
	\includegraphics[width=0.6\textwidth]{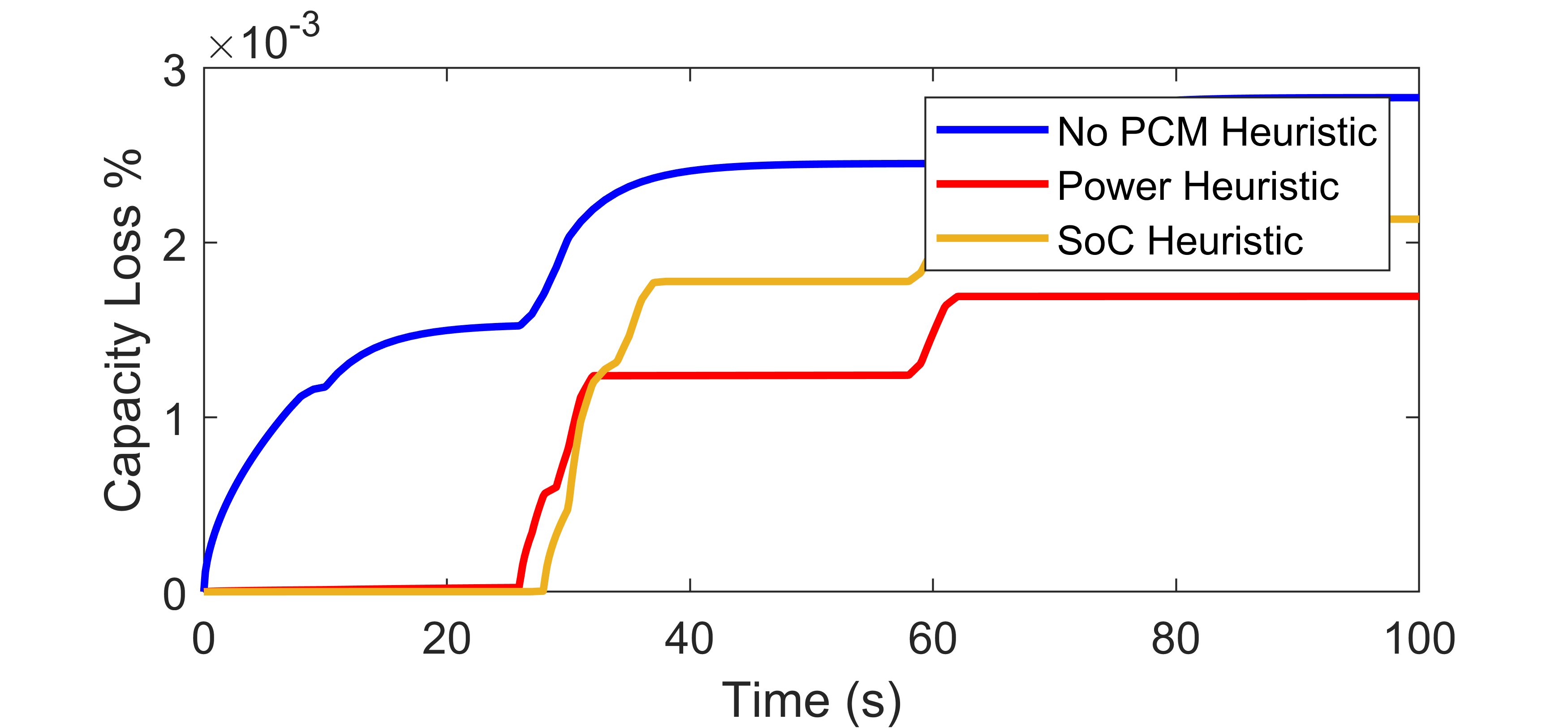} 
	\caption{Capacity loss $\%$ for different heuristic scenarios}
	\label{Cap_Loss}
\end{figure}

Finally, Figure \ref{Cap_Loss} shows the \textit{Capacity Loss $\%$} of the PCM during its operation for different scenarios. It can be seen that Scenario 2 has a minimal impact on the PCM capacity loss for the operation duration, especially during the occurrence of pulse power. 

Thus, from the results, and testing different Scenarios it can be deduced that the PCM degradation minimization comes at the expense of pushing the PGM to its operational limits. However, this trade-off can be tuned by adjusting the penalties $\beta,\gamma_p$, and $\gamma_q$.

\subsubsection{Conclusion}\label{sec: conclusion}
This section presents an MPC-based battery (PCM)-degradation-aware energy management strategy for a shipboard power system under high-ramp-rate pulse power load situations. The individual power generation module, power conversion module ramp rate limitations, the PGM-rated condition, and multiple PCM degradation heuristics in terms of PCM power and PCM SoC are considered in the problem formulation while allocating the power in the SPS. The model-based PCM degradation measures based on the PCM power and PCM state of charge are used to capture the PCM usage and thus minimize it. The effect of the PCM power minimization and PCM state of charge has been observed on the PCM capacity loss. Finally, a numerical example has been given to show the effectiveness of the proposed EM method and study which heuristics have a greater impact in mitigating PCM degradation.

\subsection{Application to Shipboard Power Systems}
Integration of loads such as pulse power loads (PPLs) into notional shipboard power systems presents a challenge in the form of meeting their high ramp rates. Ramp rate limitations on the generators can hinder power flow in the shipboard power systems. Failure to meet the ramp rate requirements may cause instability. Aggregating generators with energy storage elements usually addresses the ramp rate requirements while ensuring the power demand is achieved. This approach proposes an energy management strategy that adaptively splits the power demand between the generators and the batteries while simultaneously considering the battery degradation and the generator's efficient operation. Since it is challenging to incorporate the battery degradation model directly into the optimization problem due to its complex structure and the degradation time scale which is not practical for real-time implementation, a reasonable heuristic in terms of minimizing the absolute battery power to manage the battery degradation is proposed. A distributed model predictive energy management strategy is then developed to coordinate the power split maximizing the generator efficiency and minimizing the battery degradation. The designed strategy is tested through two numerical case studies on a single generator, battery, and load model and a notional shipboard power system on a real-time performance target machine. The results show the effect of the designed energy management strategy in mitigating battery degradation and its health management.

\subsection{Introduction}
Conceptualization of Micro-grids (MGs) as an amalgamation of fossil-fuel-driven micro-turbines, renewable energy-based generating systems, and energy storage systems (ESSs) incorporated into a distribution generation environment, driving high demand for the power was first presented in \cite{Lasseter_2001}. MGs are categorized into \emph{grid-connected} and \emph{islanded}, based on the configuration and their mode of operation \cite{Asanso_2007}. One such example of an islanded operation is Shipboard Power Systems (SPSs). In shipboard power systems, power generation modules (PGMs) and power conversion modules (PCMs) provide the power required by various loads through the direct current (DC) distribution system. Modern SPSs are equipped with advanced power loads such as rail guns, electromagnetic radars, and heavy nonlinear and pulsed power loads (PPLs) \cite{Derry_1} \cite{Derry_2}. These advanced loads are high-ramp-rate loads \cite{Kim_2015}. 

There is a challenge to utilize the existing PGMs, which are ramp-rate limited, to meet the power demanded by the high-ramp-rate loads. Failure to ramp up in appropriate time to meet the load requirements can lead to system imbalances and instability. Furthermore, adding additional generators is not a plausible solution considering the ship hull dimension restrictions. This imposes a need for high-ramp-rate power generation elements to be integrated into SPSs. Energy storage systems such as batteries can support higher ramp rates. Integrating energy storage systems (ESSs) into the existing notional SPS framework addresses the high ramp requirements \cite{ESRDC_1}. The ESSs can store the energy during normal periods and when needed can provide immediate ramp-up support. Several investigators have addressed the above challenge \emph{without} considering battery degradation \cite{2017_Vu_2,2017_Vu_3,2021_Vedula,2017_Zohrabi,Zhang_2022,Seenumani,reddy,Richard} where they have used dynamic programming (DP) and model predictive control (MPC) to address the aforementioned challenge. MPC is a renowned concept used in control literature whose applications include industry (see Survey \cite{Qin2003ASO}). MPC problem solves an objective function subject to a set of constraints that can be actuator limitations in physical systems or rate constraints (\cite{osti_206444,Wang_Boyd,Anubi_2015}). The main \emph{advantage} of the MPC over DP is that it enables incorporation of the constraints such as inequality and box constraints capturing the upper and lower limitations and ramp-rates, thus accounting for actuator limitations \cite{Mayne_2000}. 

\begin{figure}[t!]
      \centering
      \includegraphics[width=0.6\textwidth]{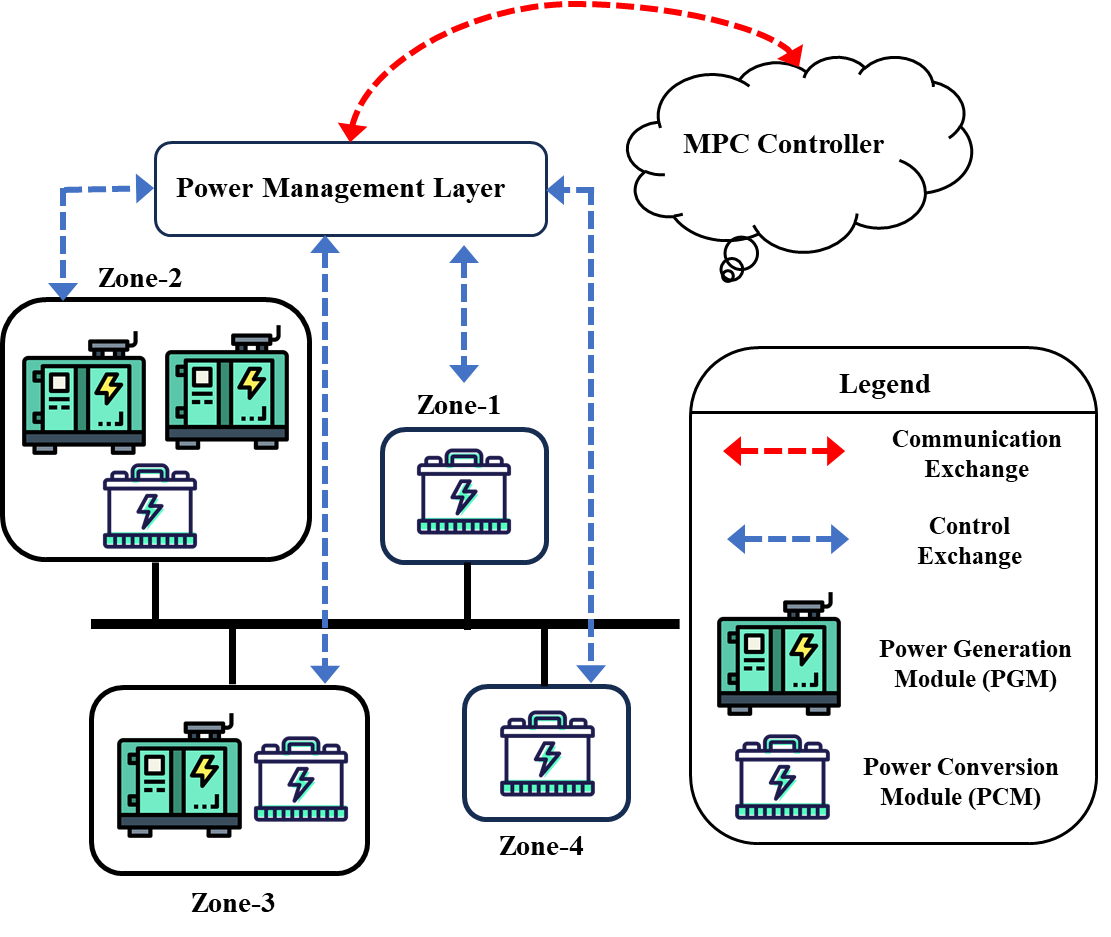}
	 \caption{A Notional Shipboard Power System with the Zonal Architecture and Control Hierarchy.}
     \label{4_Zone_SPS} 
\end{figure}

Nevertheless, batteries degrade faster than generators. This fact imposes another challenge to consider both \emph{ramp rate requirements} and \emph{battery degradation} together. Some investigators have considered this challenge while addressing the EM problem for the MGs \cite{Hein_2021,Steen_2021,Li_2023,Zhao_2023,Nawaz_2023,Ji_2023,Wang_2022,10590856}. In \cite{Hein_2021}, the authors use DP to solve the EM problem while considering the depth of discharge (DoD) as a measure for capturing battery degradation. In \cite{Steen_2021}, an energy market pricing-based EM problem for an MG, considering the battery DoD as a degradation measure, is presented. Market-energy-pricing-based approach for MG with fast Li-ion (standing for Lithium-ion)-based batteries, considering DoD as a measure for capturing the battery degradation, has been presented in \cite{Li_2023}. DoD has been also used as a battery degradation measure to solve the EM problem in \cite{Zhao_2023} where the authors proposed a neural-network-based method to address the EM problem in the MGs. In \cite{Nawaz_2023}, the authors have presented a state-of-charge (SoC)-based EM scheme for MGs; the SoC is regulated around a reference value to minimize battery degradation. Monitoring the battery charge/discharge cycles and using them to determine the battery degradation with an application to the MGs has been investigated in \cite{Ji_2023}. The work of \cite{Hein_2021,Steen_2021,Li_2023,Zhao_2023,Nawaz_2023,Ji_2023} consider battery aging in \emph{different measures} such as DoD, SoC, and charge/discharge based cycle count. In \cite{Wang_2022}, the measures used to minimize the battery usage are \emph{battery power} and \emph{battery SoC}, but the authors in \cite{Wang_2022} did \emph{not} consider the cost associated with the generators and its effect on determining the battery degradation. 

The main contributions are:
\begin{itemize}
    \item A model predictive energy management strategy is designed to incorporate battery degradation in the form of minimizing its power as a heuristic (measure) to mitigate battery deterioration.
    \item Considered both the maximizing generator efficiency and minimizing battery degradation in the designed model predictive energy management problem along with the component ramp rate limitations to test the performance under the presence of pulsed power loads.
    \item A scalable distributed plug-and-play model predictive energy management strategy is developed and tested on a shipboard power system with the real component ratings provided by the Office of Naval Research (ONR).
    \item Real-Time simulation on the SPEEDGOAT performance target machine validating the proposed model predictive energy management strategy.
\end{itemize}

\subsection{Shipboard Power System Model} 

In this section, we present the model development and the device level controller structure for a  Shipboard Power System components \emph{generators/PGMs} and \emph{batteries/PCMs} supplying a common load through a unified DC bus.  

The \emph{DC Generator or the PGM} model consists of a current-controlled DC voltage source coupled with a series $RL$ impedance connected to the bus. The assumption is that the bus voltage is already regulated to a set point. Thus, it is considered to be a constant. The dynamical model PGM is given as follows
\begin{equation}\label{gen_dynamics}
 l_{g} \frac{d i_g}{dt} = -r_{g}  i_g(t)+\Delta v(t)
\end{equation}
where $i_g(t) \in \mathbb{R}$ is the PGM current and $\Delta v(t)={v}_{bus}(t)-v_g(t)$, where $v_g(t) \in \mathbb{R}$ is the controllable voltage source and $v_{bus}(t) \in \mathbb{R}$ is the bus voltage to which the PGM is coupled. The generator inductance $l_g \in \mathbb{R}_+$ is in \textsf{Henry}, $r_g \in \mathbb{R}_+$ is the PGM resistance in \textsf{Ohm}. The current injected by the PGM into the bus is dictated by controlling the voltage source.  Given an optimal power profile $p_{g_r} \in \mathbb{R}$ from the energy management control layer, the reference current $i_{g_r} \in \mathbb{R}$ (assumed to be bounded $\in \mathcal{L}_{\infty}$) is generated as $i_{g_r}=p_{g_r}/{v_{bus}}$ and the local controller input $v_g$ is determined via a closed-loop control scheme in which the device level state tracks the reference current i.e. $i_g(t) \to i_{g_r}(t)$. Thus, the DLC is designed to achieve the following objective
$$\int_{0}^{\infty}\bigg(\underbrace{i_g(t)-i_{g_r}(t)}_{\tilde{i}_{g}}\bigg)^2 dt < \infty.$$
Since the design and analysis of such a DLC is not the main objective of this paper, numerous linear and nonlinear methods can be employed to design it \cite{Khalil_Book}. Thus, we provide a simple control design along with the closed-loop stability analysis of the designed controller.

Consider the current error $$\tilde{i}_g(t) = i_g(t)-i_{g_r}(t),$$ consequently, the control law is designed to have the following form 
\begin{equation}\label{control_law}
    \Delta v(t) = K \tilde{i}_g(t) + \Delta v_0(t),
\end{equation}
where $K$ is the control gain and $\Delta v_0(t) = r_g i_{g_r}(t)$ is the feed-forward term. Substituting (\ref{control_law}) in the open-loop dynamics (\ref{gen_dynamics}) yields
\begin{equation}
    l_g \frac{d \tilde{i}_g}{dt} = -(r_g - K) \tilde{i}_g(t).
\end{equation}
Consider the following Lyapunov candidate function
\begin{align*}
V(t) &= \frac{1}{2}l_g \tilde{i}_g(t)^2, \\
\dot{V}(t) &= l_g \tilde{i}_g(t) \frac{d \tilde{i}_g}{dt}, \\
&= -(r_g - K) \tilde{i}_g(t)^2, 
\end{align*}
since $r_g \in \mathbb{R}_+$ and with the choice of $K < 0$, the lyapunov candidate function derivative $\dot{V}(t)$ is negative-semi-definite (NSD). Since $V(t) > 0$ and $\dot{V}(t) \leq 0$ it implies $V(t) \in \mathcal{L}_{\infty}$ implies  $\tilde{i}_g(t) \in \mathcal{L}_{\infty}$. From the assumption that the current reference is bounded, the PGM current $i_g(t) \in \mathcal{L}_{\infty}$. Consequently the control input $\Delta v(t) \in \mathcal{L}_{\infty}$. Consequently the current error dynamics $\frac{d\tilde{i}_g}{dt} \in \mathcal{L}_{\infty}$ implies $\tilde{i}_g(t)$ is uniformly continuous.
Consider
\begin{align*}
    \int_{0}^{\infty} \dot{V}(t)dt &= -(r_g - K) \int_{0}^{\infty} \tilde{i}_g(t)^2 dt, \\
\int_{0}^{\infty} \tilde{i}_g(t)^2 dt &\leq \frac{V(0) - V(\infty)}{ \left|r_g + K\right|}.
\end{align*}
This concludes the DLC control design and the closed-loop stability analysis including the signal chasing for the PGM.

The \emph{Battery Energy Storage System} (BESS) or the \emph{PCM} model consists of multiple ESSs modeled as a voltage-controlled current source. The static model consists of a resistance $r_b \in \mathbb{R}_+$ coupled to the bus $v_{bus} \in \mathbb{R}$ in series with a controllable voltage source $v_b \in \mathbb{R}$ and an open circuit voltage $v_{oc} \in \mathbb{R}$ that make up the battery model. The algebraic equation governing the power exchange in the battery is given as:
\begin{align*}v_{b} &= \frac{v_{bus}^2-p_{b}r_b-v_{bus}v_{oc}}{v_{bus}}, \\
i_b &= \frac{v_{bus}-v_b-v_{oc}}{r_b}.\end{align*}
The state of charge dynamics (SoC) of the battery is given as:
\begin{equation}\label{SoC}
    \frac{d SoC}{dt} = -\frac{1}{Q}i_b(t),
\end{equation}
where $Q \in \mathbb{R}_+$ is the battery capacity in \textsf{ampere-hour}. The discretized form of the SoC dynamics given in (\ref{SoC}) is:
\begin{equation}\label{SoC_Power_Discrete}
SoC_{k+1}=SoC_k-\frac{T_d} {Q}i_{b_k},\end{equation} where $T_d \in \mathbb{R}_+$ is the discretization time-step. Here, $k,k+1$ represents the discrete time intervals. Although the SoC dynamics above do not explicitly capture the physical restriction on the $SoC \in \left[0, 1\right]$. The constraint is imposed in the optimization problem of the proposed energy management strategy.

\begin{figure}[t!]
      \centering
      \includegraphics[width=0.6\textwidth]{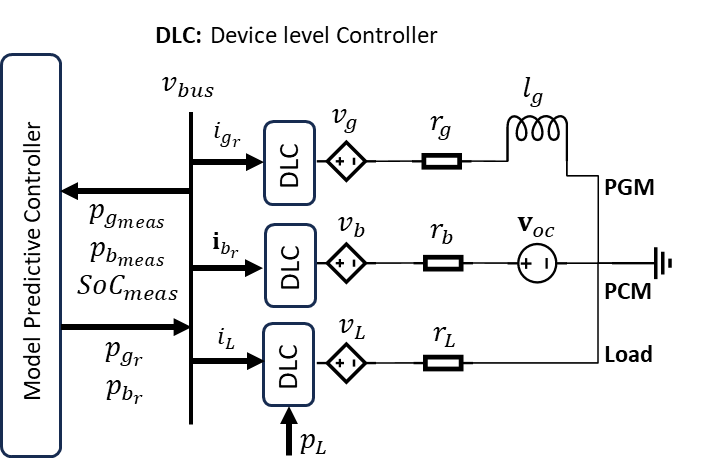}
	 \caption{Mathematical topology of the power exchange in an SPS with a circuit level schematic of PGM, PCM, and load models.}
     \label{Mathematical_Model} 
\end{figure}

The PCM degradation model used in this work is an Arrhenius equation-based model which uses $Ah$-throughput as $\displaystyle \int_{0}^{t}\left|i_b(\tau)\right|d\tau$ as a metric to evaluate battery state of health (SoH). $i_b(t) \in \mathbb{R}$ is the current drawn from the individual PCM (positive while discharging and vice-versa). The PCM \textit{capacity loss} formulation is given as follows \cite{SONG2018433} \cite{Wang_Dai}.
\begin{equation}
    Q_L(t) = \zeta_1 e^{\frac{-\zeta_2+TC_{r}}{RT}}\int_{0}^{t}\left|i_b(\tau)\right|d\tau,
\end{equation}
where $T \in \mathbb{R}_+$ is the PCM baseline operating temperature in \textsf{Kelvin}, $C_r$ is the C-rate of the PCM.  

The capacity loss (\%) is given as follows: $$\Delta Q \% = \frac{Q-Q_{L}(t)}{Q}\times 100,$$
where $Q_{L}$ is the \emph{capacity loss} of the BESS in \textsf{ampere-hour}. \textit{This does not capture the battery end of life (EOL)}. It only captures the capacity lost during the battery operation ($\Delta Q$).

The \textit{Power Load} is modeled as a static resistive load $r_L$ model which consumes the power generated by the PGMs and the PCMs. A variable voltage $v_L \in \mathbb{R}$ enables the current flow in the load based on the active power profile $p_L \in \mathbb{R}$. The assumption is that there is only one load in the model which is consuming all the generated power. Thus, The load model is given as follows:
\begin{equation}\label{power_load}
\begin{aligned}
     v_L = \frac{v_{bus}^2-p_Lr_L}{v_{bus}} \\
     i_L = \frac{v_{bus}-v_L}{r_L},
\end{aligned}
\end{equation}

Thus, the power flow in the SPS is given by the following algebraic equation:
\begin{equation}\label{PowerFlow}
    \sum_{i=1}^{n_g} \mathbf{p}_{g_i}+\sum_{j=1}^{n_b} \mathbf{p}_{b_j}-p_L = 0,
\end{equation}

where $n_g, n_b \in \mathbb{N}$ are the number of PGMs and PCMs in SPS. Figure \ref{Mathematical_Model} shows the hierarchical control structure in the SPS and the power exchange. The power and the SoC measurements are sent to the MPC layer. The MPC computes the optimal values based on the objective and set of constraints and the optimal values act as a reference for DLCs. 

\subsection{Proposed Energy Management Method}
This section presents the control development \cite{vedula2024distributedmodelpredictiveenergymanagement} and the derivation of the proposed EM strategy. First, we present the MPC problem formulation in the general form. Then, we consider the objective functions of the generator and the battery respectively, and explain the reason behind choosing the particular objectives. Finally, the distributed plug-and-play formulation for the centralized MPC problem is derived. Consider the MPC problem of the form:
\begin{equation}\label{MPC_main}
\begin{aligned}
\Minimize_{\mathbf{p}_{g_i},\mathbf{p}_{b_j}} \quad & \sum_{{i}=1}^{n_g}{C}_{g_{i}}(\mathbf{p}_{{g}_{i}})+\sum_{{j}=1}^{n_b}{C}_{b_{j}}(\mathbf{p}_{{b}_{j}})\\
\SubjectTo \quad & \sum_{{i}=1}^{n_g}\mathbf{p}_{{g}_{i}}+\sum_{{j}=1}^{n_b}\mathbf{p}_{{b}_{j}} - p_f\mathbf{1} = 0, \hspace{1mm} \forall  k,\\
& \mathbf{p}_{{g}_{{i}}} \in \mathcal{X}_g, \hspace{1mm} \forall k, \\& \mathbf{p}_{{b}_{j}}  \in \mathcal{X}_b, \hspace{1mm}\forall k,
\end{aligned}
\end{equation}
where $k \triangleq \left[\begin{array}{cccc} 1&2&\hdots&h\end{array}\right]^{\top} \in \mathbb{R}^h$ is the length of the prediction horizon, \newline $\mathbf{p}_{g_i} \triangleq \left[\begin{array}{cccc}{\mathbf{p}_{g_i}}_k&{\mathbf{p}_{g_i}}_{k+1}&\hdots&{\mathbf{p}_{g_i}}_{k+h-1}\end{array}\right]^{\top}\in\mathbb{R}^h$ is the power profile for generator $i$ over the prediction horizon of length $h$, $\textbf{p}_{b_j} \triangleq \left[\begin{array}{cccc}{\mathbf{p}_{b_j}}_k&{\mathbf{p}_{b_j}}_{k+1}&\hdots&{\mathbf{p}_{b_j}}_{k+h-1}\end{array}\right]^{\top}\in\mathbb{R}^h$ is the power profile for battery $j$ over the prediction horizon of length $h$ and ${p}_f\mathbf{1} \in\mathbb{R}^h$ is the desired total power held fixed over the prediction horizon. The generator cost ${C}_{g_i}:\mathbb{R}^{{h}}\longrightarrow\mathbb{R}_+$ could be associated with an efficiency map or operating at a desired rated power. The battery cost ${C}_{b_j}:\mathbb{R}^{{h}}\longrightarrow\mathbb{R}_+$ could be the cost associated with battery health monitoring and degradation management. ${h} \in \mathbb{N}$ represents the prediction horizon. $\mathcal{X}_g \subset \mathbb{R}^{h}$ and $\mathcal{X}_b \subset \mathbb{R}^h$ represents the set of inclusion constraints each of the PGM and PCM summarizing attributes such as upper and lower power limits, ramp rate limits, and relevant system dynamics. The optimization variables are optimized for the entire length of the horizon and the first value of the sequence is supplied/applied as a power reference (control input) to the low-level (DLC).

The PGM and the PCM costs are considered so that the objective is to maintain the PGM around a rated value ($\mathbf{p}_g^r$) which is known and to minimize the PCM power. Thus the choice for the PGM and the PCM cost functions is $$C_g(\mathbf{p}_g) \triangleq  \norm{\mathbf{p}_g-\mathbf{p}_g^r}_2^2 \hspace{2mm} \text{and} \hspace{2mm} C_b(\mathbf{p}_b) \triangleq  \norm{\mathbf{p}_b}_2^2.$$ Based on the optimization problem presented in (\ref{MPC_main}), consider the following MPC problem with the respective cost functions for the PGM and the PCM as follows:
\begin{equation} \label{Full_MPC}
\begin{aligned}
\Minimize_{\mathbf{p}_{g},\mathbf{p}_{b},\mathbf{SoC}_b} \quad & \sum_{i=1}^{n_g}\frac{\beta_i}{2}\norm{\mathbf{p}_{g_{i}}-\mathbf{p}^{r}_{{g_i}}}_2^2+\sum_{j=1}^{n_b}\frac{\gamma_j}{2}\norm{\mathbf{p}_{b_{j}}}_2^2\\
\SubjectTo \quad & \sum_{i=1}^{n_g}\mathbf{p}_{g_{i}}+\sum_{j=1}^{n_b}\mathbf{p}_{b_{j}} = p_f\mathbf{1}, \forall k, \\
\quad & \mathbf{SoC}_{b_{jk+1}} = \mathbf{SoC}_{b_{jk}}-\frac{T_d}{Qv_{bus}} \mathbf{p}_{b_{jk}}, \forall k, \\ \quad & \underline{\mathbf{p}}_{g} \preceq \mathbf{p}_{g_{i}} \preceq \overline{\mathbf{p}}_{g}, \forall k,\\
& \left|\mathbf{p}_{g_{ik}}-\mathbf{p}_{g_{ik-1}} \right| \preceq {r}_{g}\mathbf{1}, \forall k \\ \quad & \underline{\mathbf{p}}_{b} \preceq \mathbf{p}_{b_{j}} \preceq \overline{\mathbf{p}}_{b}, \forall k,\\
& \left|\mathbf{p}_{b_{jk}}-\mathbf{p}_{b_{jk-1}}\right| \preceq {r}_{b} \mathbf{1}, \forall k, \\
& \underline{\mathbf{SoC}}_b \preceq \mathbf{SoC}_{b_{jk}} \preceq \overline{\mathbf{SoC}}_b, \forall k,
\end{aligned} 
\end{equation}
where $k \triangleq \left[\begin{array}{cccc} 1&2&\hdots&h\end{array}\right]^{\top} \in \mathbb{R}^h$, $r_g \in \mathbb{R}_+$ and $r_b \in \mathbb{R}_+$ are the ramp rates of the generator and the battery respectively. $\underline{\mathbf{p}}_{g}$ and $\overline{\mathbf{p}}_{g}$ are the lower and the upper power limitations on the generator, respectively. $\underline{\mathbf{p}}_{b}$ and $\overline{\mathbf{p}}_{b}$ are the lower and the upper power limitations on the battery. $\underline{\mathbf{SoC}}_b$ and $\overline{\mathbf{SoC}}_b$ represent the lower and the upper limitations on the state of charge of the battery. $\beta_i \in \mathbb{R}_+^{n_g}$ are the tunable parameters penalizing the PGMs deviation from a known efficient operation point. $\gamma_j \in \mathbb{R}_+^{n_b}$ are the tunable parameters penalizing the battery power. The constraint set $\mathcal{X}_g$ and $\mathcal{X}_b$ are \emph{polytopes}. Based on this fact the optimization problem in (\ref{Full_MPC}) is feasible, the optimization problem is convex whose global optimal solution can be found in polynomial time by existing algorithms \cite{boyd_vandenberghe_2004}.

\subsubsection{Distributed Formulation}

The optimization problem in (\ref{MPC_main}) is broken down into the individual optimization problems solved at the respective PGMs and the PCMs coupled via an aggregator. The overall algorithm development is presented in this subsection. Considering the Lagrangian for the problem in (\ref{MPC_main}) 
\begin{align}\nonumber
\mathcal{L}(\textbf{p}_{{g}_{i}},\textbf{p}_{{b}_{j}},\boldsymbol{\lambda}) &= \sum_{{i}=1}^{n_g}{C}_{{g_i}}(\mathbf{p}_{{g}_{i}})+\sum_{{j}=1}^{n_b}{C}_{{b_j}}(\mathbf{p}_{{b}_{j}})+\mathcal{I_X}_g+\mathcal{I_X}_b
+\boldsymbol{\lambda}^\top\bigg(\sum_{{i}=1}^{n_g}\mathbf{p}_{{g}_{i}}+\sum_{{j}=1}^{n_b}\mathbf{p}_{{b}_{j}} - {p}_{f}\mathbf{1}\bigg)\\\nonumber
&= \mathcal{I_X}_g+\mathcal{I_X}_b+ \sum_{{i}=1}^{n_g}\bigg({C}_{{g_i}}(\mathbf{p}_{{g}_{i}})+\boldsymbol{\lambda}^\top\mathbf{p}_{{g}_{i}}\bigg)+\sum_{{j}=1}^{n_b}\bigg({C}_{{b_j}}(\mathbf{p}_{{b}_{j}})+ \boldsymbol{\lambda}^\top\mathbf{p}_{{b}_{j}}\bigg)\label{Lagrangian_MPC}
-p_f\boldsymbol{\lambda}^\top\textbf{1},
\end{align}
where $\boldsymbol{\lambda} \in \mathbb{R}^{h}$ is the dual variable associated with the main problem in (\ref{MPC_main}) and $\mathcal{I_X}_g$ and $\mathcal{I_X}_b$ are the indicator functions for the inclusion constraints. Given $\boldsymbol{\lambda}$, let $\mathbf{p}_{g_i}^*$ and $\mathbf{p}_{b_j}^*$ to be the solutions to the optimization problem in (\ref{MPC_main}) (\cite{boyd_vandenberghe_2004}) 
\begin{subequations}\label{Iteration_Updates_1}
\begin{align}
    \textbf{p}^*_{g_i} (\boldsymbol{\lambda}) = \argmin_{{\mathbf{p}_{{g}_{i}}} \in \mathcal{X}_g} \bigg\{\text{C}_{{g_i}}(\mathbf{p}_{{g}_{i}})+\boldsymbol{\lambda}^{\top}\textbf{p}_{{g}_{i}}\bigg\}, \\
    \textbf{p}^*_{b_j}(\boldsymbol{\lambda}) = \argmin_{{\mathbf{p}_{{b}_{j}}} \in \mathcal{X}_b} \bigg\{\text{C}_{{b_j}}(\mathbf{p}_{{b}_{j}})+\boldsymbol{\lambda}^{\top}\textbf{p}_{{b}_{i}}\bigg\}.
\end{align}
\end{subequations}
Thus, the dual problem is given as:
$$\Maximize \quad \mathcal{L}(\textbf{p}^*_{g_i},\textbf{p}^*_{b_j},\boldsymbol{\lambda}) \triangleq d(\boldsymbol{\lambda}), $$
which is solved using the gradient ascent algorithm
\begin{align}\nonumber
    \boldsymbol{\lambda}_{t+1} &= \boldsymbol{\lambda}_t +  \alpha \nabla d(\boldsymbol{\lambda}_t)\\\label{eqn:dual_ascent1}
                               &= \boldsymbol{\lambda}_t + \alpha \left(\sum\limits_{i=1}^{n_g}\mathbf{p}_{g_i}^*+\sum\limits_{j=1}^{n_b}\mathbf{p}_{b_j}^* - p_f\mathbf{1}\right).
\end{align}
for a step size $\alpha > 0$ and $t$ is the \emph{iteration counter}. Since the minimizing power profiles $\textbf{p}^*_{g_i}$ and $\textbf{p}^*_{b_j}$ depends on $\boldsymbol{\lambda}$, the dual ascent step in \eqref{eqn:dual_ascent1} becomes a recursion and challenging to evaluate. Consequently, we incorporate the primal problems in \eqref{Iteration_Updates} within a sequential minimization-maximization scheme \cite{boyd_admm} as follows:
\begin{subequations}\label{Iteration_Updates}
\begin{align}
    \mathbf{p}_{g_i}^{t+1} &= \argmin_{{\mathbf{p}_{{g}_{i}}} \in \mathcal{X}_g} \bigg\{{C}_{{i}}(\mathbf{p}_{{g}_{i}})+\boldsymbol{\lambda}^{\top}_t\mathbf{p}_{{g}_{i}}\bigg\}, \\
    \mathbf{p}_{b_j}^{t+1} &= \argmin_{{\mathbf{p}_{{b}_{j}}} \in \mathcal{X}_b} \bigg\{{C}_{{j}}(\mathbf{p}_{{b}_{j}})+\boldsymbol{\lambda}^{\top}_t\mathbf{p}_{{b}_{j}}\bigg\},\\
    \boldsymbol{\lambda}_{t+1} &= \boldsymbol{\lambda}_{t}+\alpha \bigg(\sum_{{i}=1}^{n_g}\mathbf{p}^{t+1}_{{g}_{i}}+\sum_{{j}=1}^{n_b}\mathbf{p}^{t+1}_{{b}_{j}} - p_f\mathbf{1}\bigg). 
\end{align}
\end{subequations}

\begin{algorithm}
\caption{Primal-Dual Min-Max Algorithm}\label{alg:cap}
\begin{algorithmic}[1]

\State \textbf{for} $t=1$ \textbf{do}
\State Get the current measurements for $\mathbf{p}_{g_i}$ and $\mathbf{p}_{b_j}$.
\State Initialize: $\boldsymbol{\lambda}_t$
\State  Solve the optimization problem (\ref{Iteration_Updates_1}a) and (\ref{Iteration_Updates_1}b) to obtain the primal optimal values $\mathbf{p}_{g_i}^*$ and $\mathbf{p}_{b_j}^*$ 
\State  Update dual variable $\boldsymbol{\lambda}_{t+1}$ using the primal optimal values and initial $\boldsymbol{\lambda}_t$ $$\boldsymbol{\lambda}_{t+1} = \boldsymbol{\lambda}_t + \alpha \left(\sum\limits_{i=1}^{n_g}\mathbf{p}_{g_i}^*+\sum\limits_{j=1}^{n_b}\mathbf{p}_{b_j}^* - p_f\mathbf{1}\right)$$
\State Check stopping criterion  $$\mathbf{p}_{g_i}^*+\mathbf{p}_{b_j}^*-\boldsymbol{\lambda}_{t+1} \leq \epsilon_{tol}$$
\If {condition is not satisfied} 
    \State go to step 4
\Else 
    \State update the power control sequence $\mathbf{p}_{g_i}^*$ and $\mathbf{p}_{b_j}^*$
\EndIf
\State Update the initial value $\boldsymbol{\lambda}_t=\boldsymbol{\lambda}_{t+1}$
\State \textbf{end for}
\State  $t \gets t+1$ and proceed to step 1.
\end{algorithmic}
\end{algorithm}

Next, more details on the resulting nodal problems in the above scheme with the choice of the cost functions and the inclusion constraints are given.
\subsubsection{PGM Node Optimization Problem}
From (\ref{Iteration_Updates}a) the optimization problem at the $i^{th}$ PGM node assuming the cost function $C_i(\mathbf{p}_{g_i})$ to be a quadratic function is given in (\ref{PGM_Node}). The objective of the considered cost is to maintain the PGM power at around a given rated value or a set point value at every instant of the horizon given by $\mathbf{p}^r_{g_i}$. This rated value is the operating point for that generator as prescribed by the manufacturer. The optimization problem is given as follows:
\begin{equation} \label{PGM_Node}
\begin{aligned}
\Minimize_{\mathbf{p}_{{g}_{i}}} \quad & \frac{\beta_i}{2}\norm{\mathbf{p}_{{g}_{i}}-\mathbf{p}^r_{{g_i}}}_2^2+\boldsymbol{\lambda}^{\top}_t  \mathbf{p}_{{g}_{{i}}}\\
\SubjectTo \quad & \underline{\mathbf{p}}_{g_i} \preceq \mathbf{p}_{{g_i}_k} \preceq \overline{\mathbf{p}}_{g_i}, \hspace{1mm} \forall k\\
&\left|\mathbf{p}_{{g_i}_k}-\mathbf{p}_{{{g}_{{i}}}_{k-1}}\right| \preceq {r}_{g}\mathbf{1}, \hspace{1mm} \forall k,
\end{aligned} 
\end{equation}
where ${r}_{g} \in \mathbb{R}_+$ is the ramp-rate limitation associated with generator. $\overline{\mathbf{p}}_{{{{g}_{i}}}} \in \mathbb{R}_+$ and $\underline{\mathbf{p}}_{{{{g}_{i}}}} \in  \mathbb{R}_+$ are the upper and lower limits on generator power respectively.

\begin{figure}[t!] 
	\centering
	\includegraphics[width = 0.6\textwidth]{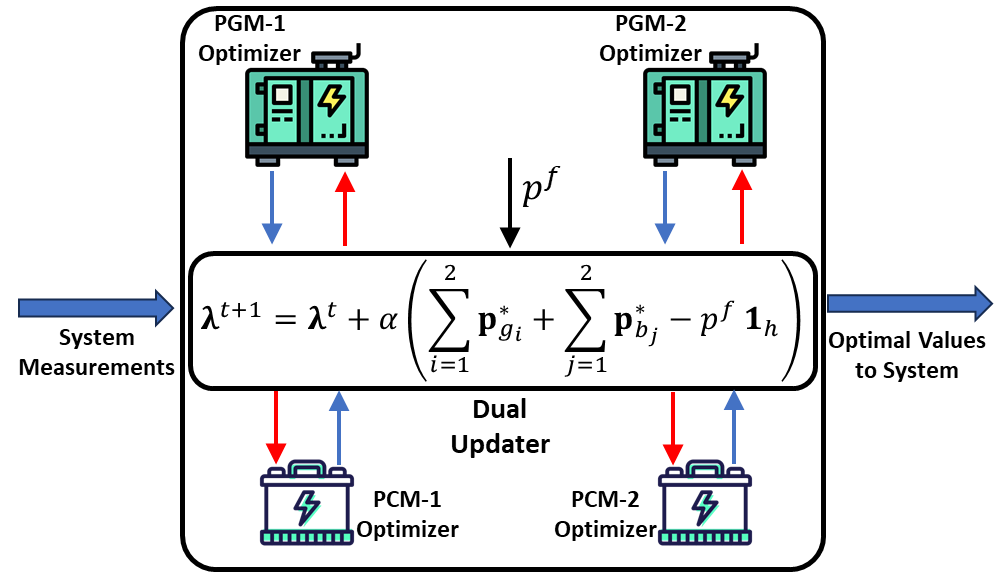} 
	\caption{Distributed dual ascent algorithm implementation schematic for numerical simulation consisting of 2-PGMs and 2-PCMs.
    }
	\label{DA_Algorithm}
\end{figure}

\subsubsection{PCM Node Optimization problem}
From (\ref{Iteration_Updates}b) the optimization problem at each battery/PCM node is given in (\ref{PCM_Node}). The discretized SoC dynamics in (\ref{SoC_Power_Discrete}) is introduced as an equality constraint. Since capturing the PCM degradation is a long process, the heuristic considered in this work is to minimize the PCM power, which translates into minimizing the PCM degradation. Thus, the cost is introduced as the weighted 2-norm of the PCM power minimization. The goal is to minimize the usage of batteries as much as possible. Thus, the optimization problem is as follows:
\begin{equation} \label{PCM_Node}
\begin{aligned}
\Minimize_{\mathbf{p}_{{b_j}},\mathbf{SoC}_{b_j}} \quad & \frac{\gamma_j}{2}\norm{\mathbf{p}_{b_j}}_2^2+\boldsymbol{\lambda}^{\top}_t  \mathbf{p}_{b_j}\\
\SubjectTo \quad &  \mathbf{SoC}_{b_{jk+1}} = \mathbf{SoC}_{b_{jk}}-\kappa \mathbf{p}_{b_{jk}}, \hspace{1mm} \forall k, \\ & \underline{\mathbf{p}}_{b_j} \preceq \mathbf{p}_{{b_j}_k} \preceq \overline{\mathbf{p}}_{b_j}, \hspace{1mm} \forall k,\\
& \underline{\mathbf{0}} \preceq \mathbf{SoC}_{{b_j}_k} \preceq \mathbf{1}, \hspace{1mm} \forall k, \\
& \left|\mathbf{p}_{{b_j}_k}-\mathbf{p}_{{b_{j}}_{k-1}}\right| \preceq {r}_{b}\mathbf{1}, \hspace{1mm}, \forall k,
\end{aligned} 
\end{equation}
where, $\displaystyle \kappa = {T_d}\slash{Q v_{bus}}$, $r_b \in \mathbb{R}_+ $  is the ramp-rate limitation associated with battery, $\overline{\mathbf{p}}_{b_j} \in \mathbb{R}_+$ and $\underline{\mathbf{p}}_{b_j} \in \mathbb{R}_-$ are the upper and lower limits on battery power. This indicates the discharging and the charging mode of the PCM. 

The solution to the optimization problems in (\ref{PGM_Node}) and (\ref{PCM_Node}) are used in calculating the dual variable update. The schematic for this aggregation and distribution setup is provided in Figure \ref{DA_Algorithm}. All the nodes present in the system send their current iterate value to the dual updater node and the information of the previous iterate of the dual variable $\boldsymbol{\lambda}_t$ is sent to all the nodes.  Thus, the dual updater, which acts as a common scheduler to all nodes implements the update law given in (\ref{Iteration_Updates}c). The convergence of the proposed algorithm under the assumed objectives and the constraint sets which are polytopes is well established in the literature \cite{boyd_vandenberghe_2004} \cite{boyd_admm}.

\subsection{Case Studies}

In this section, the developed control formulations and the model in the previous section are used to demonstrate the effectiveness of the MPC EM scheme. 
\subsubsection{Load Scenario}
\begin{figure}[h!]
      \centering
      \includegraphics[width=0.6\textwidth]{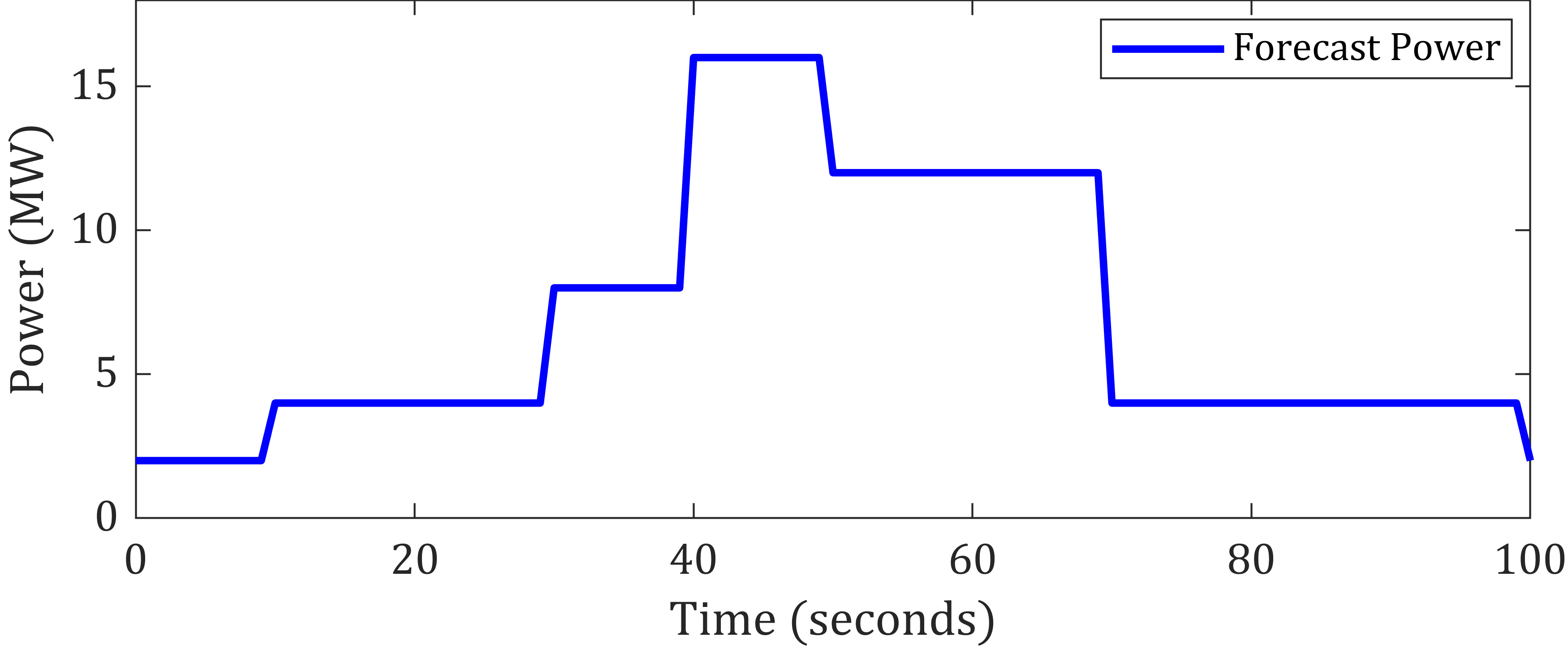}
	 \caption{The power profile used in the numerical simulations based on the load rating.}
     \label{Load Forecast} 
\end{figure}

Since the main focus of this study is to analyze the power split between the PGMs and the PCMs under the presence of the pulsed power loads, we design a pulsed load profile based on the maximum load rating of the load provided by the Office of Naval Research ONR (\cite{ESRDC_1}). The designed load is the power reference provided to the power load designed in (\ref{power_load}). Figure \ref{Load Forecast} shows the pulsed load designed to be at $16$\textsf{MW} at peak operation with three incremental pulse occurrences that are beyond the ramp capacity of the PGMs. The designed load profile is used in all of the simulation scenarios presented in this paper.

\subsubsection{Single PGM, PCM, and Load}
First, the designed model predictive energy management strategy is tested on a single PGM, PCM, and Load scenario to fully study and analyze the effect of the weights $\beta$ and $\gamma$ on the total energy delivered by the PCM, PGM, and its effect on the PCM capacity loss \%. The simulation environment used in this setup is based on MATLAB-Simulink. The simulation is run on a desktop with the following configurations Digital Storm Intel Core i9-13900K (5.7GHz Turbo) with 64GB RAM. The optimization problem presented in (\ref{Full_MPC}) with $n_g = 1$ and $n_b=1$ is implemented to study the trade-off characteristics. The MPC is implemented using YALMIP (\cite{Lofberg_2004}). The parameters presented in Table \ref{tab:rated} are based on the documentation and the notional SPS component sizing provided by the Office of Naval Research (ONRs) ESRDC (\cite{ESRDC_1}). For all the results presented, the prediction horizon is fixed at ($h=5$\textsf{seconds}). The simulation was run at the fixed step of $10^{-3}$\textsf{seconds}. The rate transition or the communication (measurement exchange) delay between the model and the MPC EM layer was set at $1$\textsf{second} to facilitate optimization convergence. 

\begin{table}[ht]
\centering
\caption{Rated values and simulation parameters for single PGM, PCM, and load Model \cite{ESRDC_1}}
\label{tab:rated}
\resizebox{0.6\columnwidth}{!}{\begin{tabular}{c|c|c}
\hline \hline
\textbf{Parameter}  & \textbf{Parameter}  & \textbf{Parameter}  \\
\textbf{Description} & \textbf{Notation}   & \textbf{Value}  \\  \hline \hline
Generator Ramp Rate    & $r_g$   & 2.9 MW/Hr   \\
ESS Ramp Rate    & $r_b$   & 10 MW/Hr  \\
Rated Bus Voltage   & $v_{bus}$   & 12 kV  \\
Generator Upper Power Limit   & $\overline{p}_g$   & 28 MW \\ 
Generator Lower Power Limit   & $\underline{p}_g$   & 0.2 MW \\ 
ESS Upper Power Limit & $\overline{p}_b$ & 10MW \\
ESS Lower Power Limit & $\underline{p}_b$ & -10MW \\
SoC Lower Limit & $\underline{SoC}_b$ & 0.4 \\ 
SoC Upper Limit & $\overline{SoC}_b$ & 0.9 \\ \hline
\end{tabular}}
\end{table}

\begin{figure}
    \centering
    \includegraphics[width=0.6\textwidth]{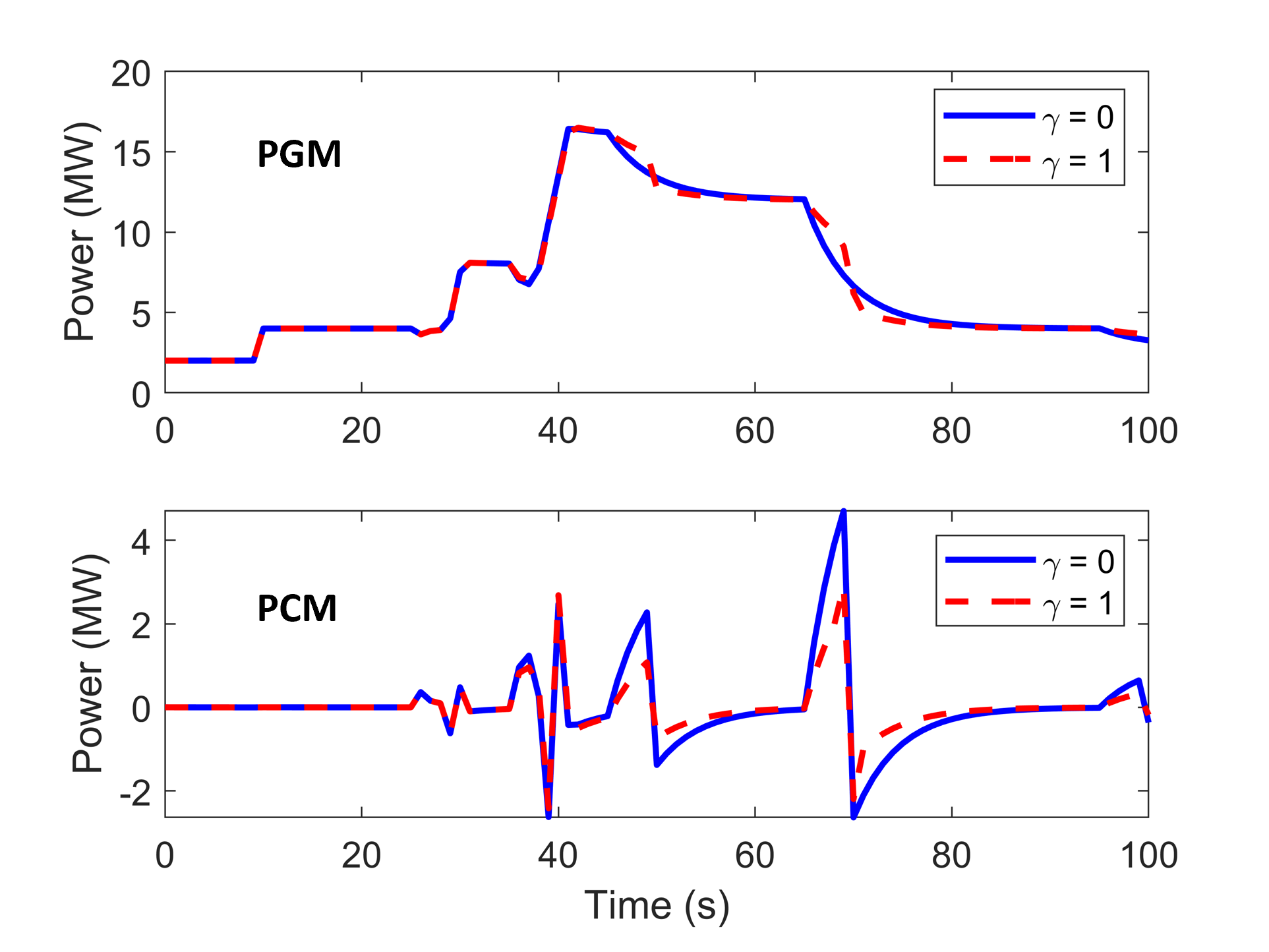}
    \caption{Power split for a single PGM and single PCM for different weights $\gamma$.}
    \label{fig:PGM_PCM_split}
\end{figure}

\begin{figure}[t!]
      \centering
      \includegraphics[width=0.6\textwidth]{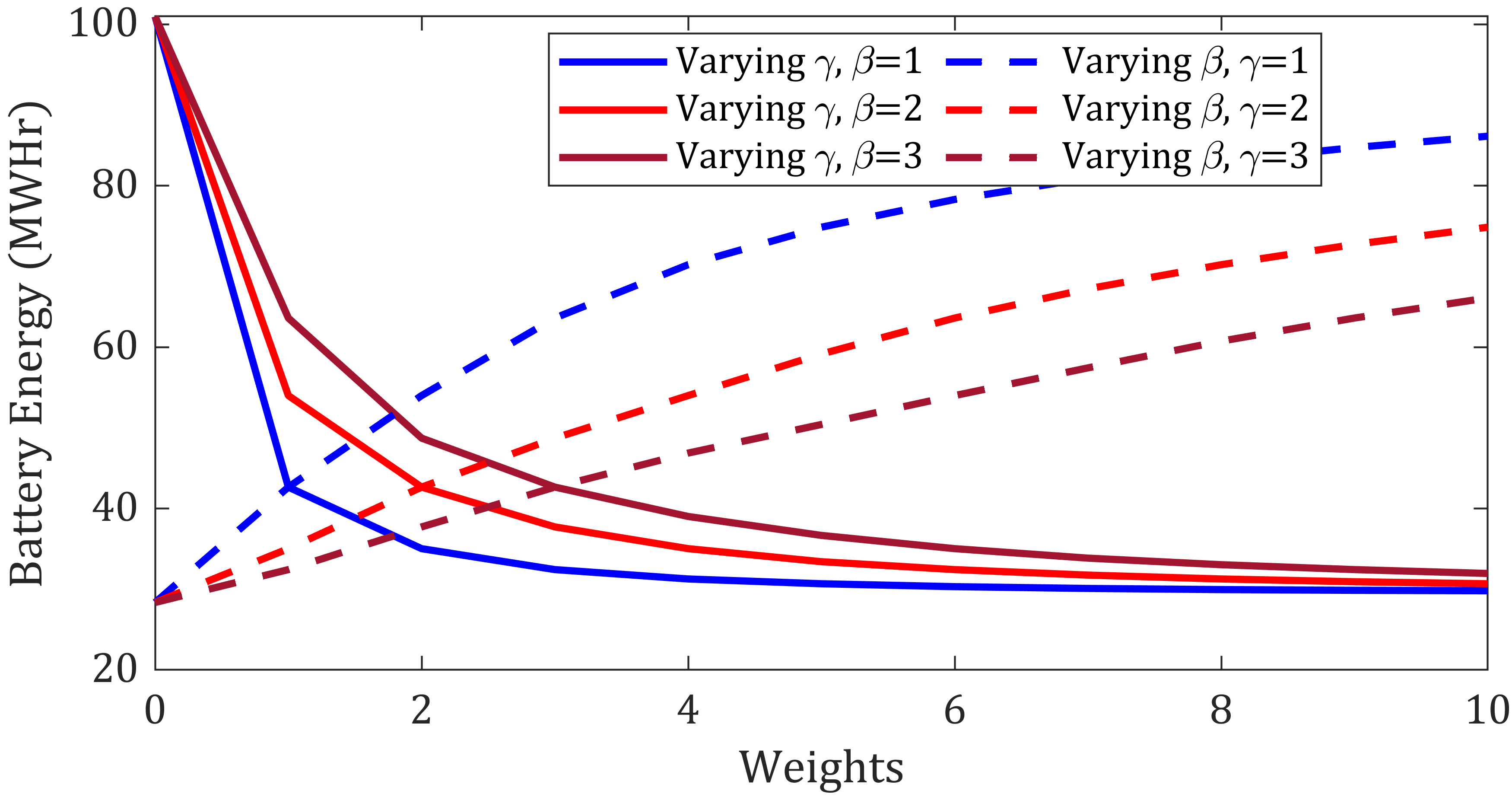}
	 \caption{The effect of different weighting of $\gamma$ and $\beta$ versus the battery energy.}
     \label{P_batt vs Weights} 
\end{figure}

Figure \ref{fig:PGM_PCM_split} shows the power split for the choice of the weight $\gamma$. Since the core focus of the work is on determining the effect of the weights $\beta$ and $\gamma$ on the generator and the battery usage and finding a trade-off between them, the results focus on the trade-off part. First, the developed EM algorithm is tested on a \emph{single PGM, PCM, and Load} model. Figure \ref{P_batt vs Weights} shows the total energy the battery utilized in \textsf{MWHr} on the y-axis versus the weight on the x-axis. The \emph{solid lines} in the figure present the effect of the battery penalty parameter $\gamma$ tuned from 0-10 against the generator penalty parameter $\beta$. As the $\gamma$ increases the total energy utilized by the battery decreases and can be seen around 30\textsf{MWHr}. This trend shifts upward for the penalty tuning of $\gamma$ from 0-10 but for the increased value of $\beta$. This shift is attributed to the power balance constraint. Similarly \emph{dotted lines} depict the effect of the generator penalty parameter $\beta$ tuned from 0-10 against the battery penalty $\gamma$. It can be seen that as the value of the $\beta$ increases, the total energy utilized by the battery increases since the emphasis is on maintaining the generator around the rated operating point. 

\begin{figure}[t!]
      \centering
      \includegraphics[width=0.6\textwidth]{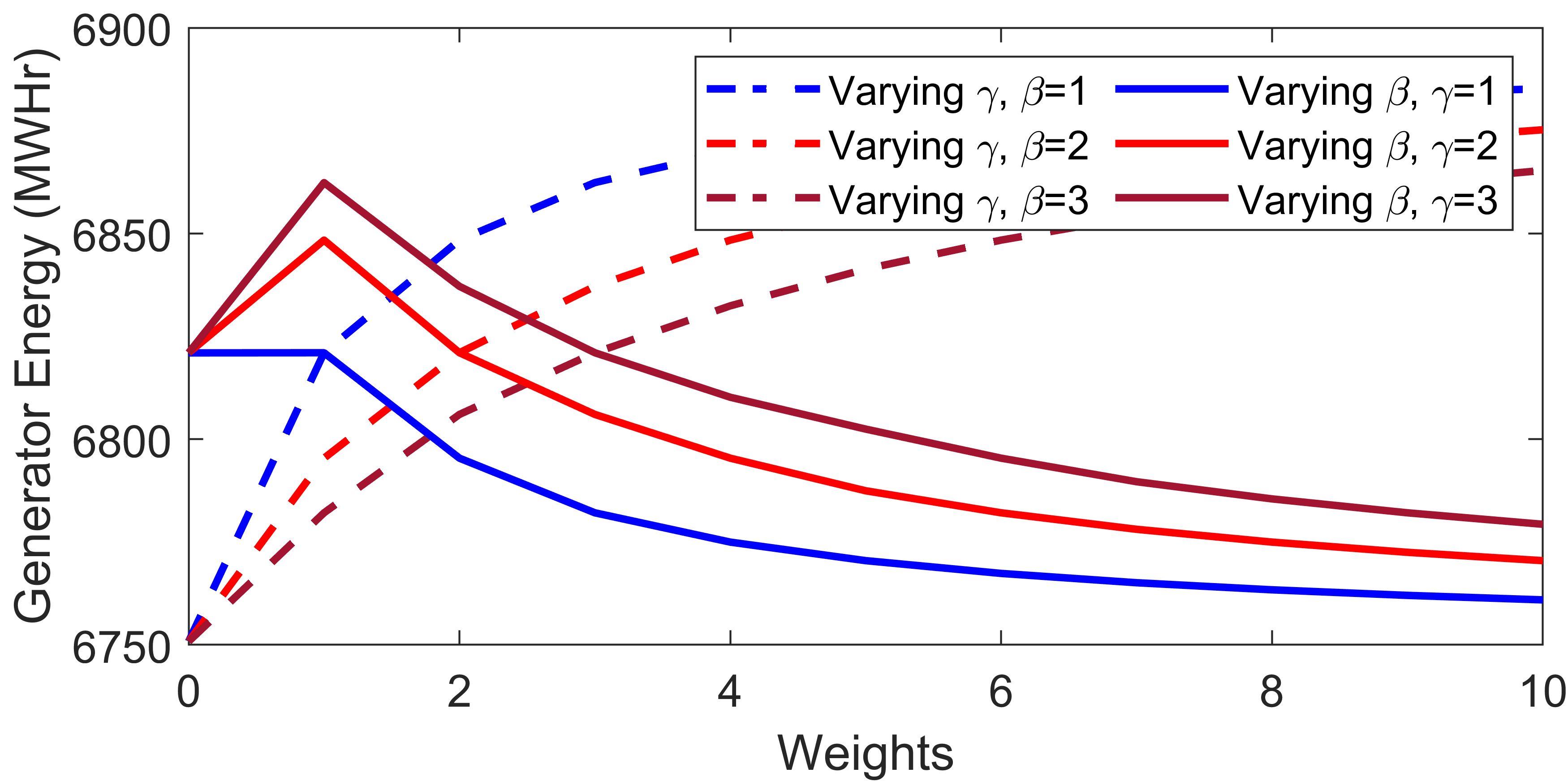}
	 \caption{The effect of different weighting $\gamma$ and $\beta$ versus the generator energy.}
  \label{Gen_Power vs Weights} 
\end{figure}

\begin{figure}[t!]
      \centering
      \includegraphics[width=0.85\textwidth]{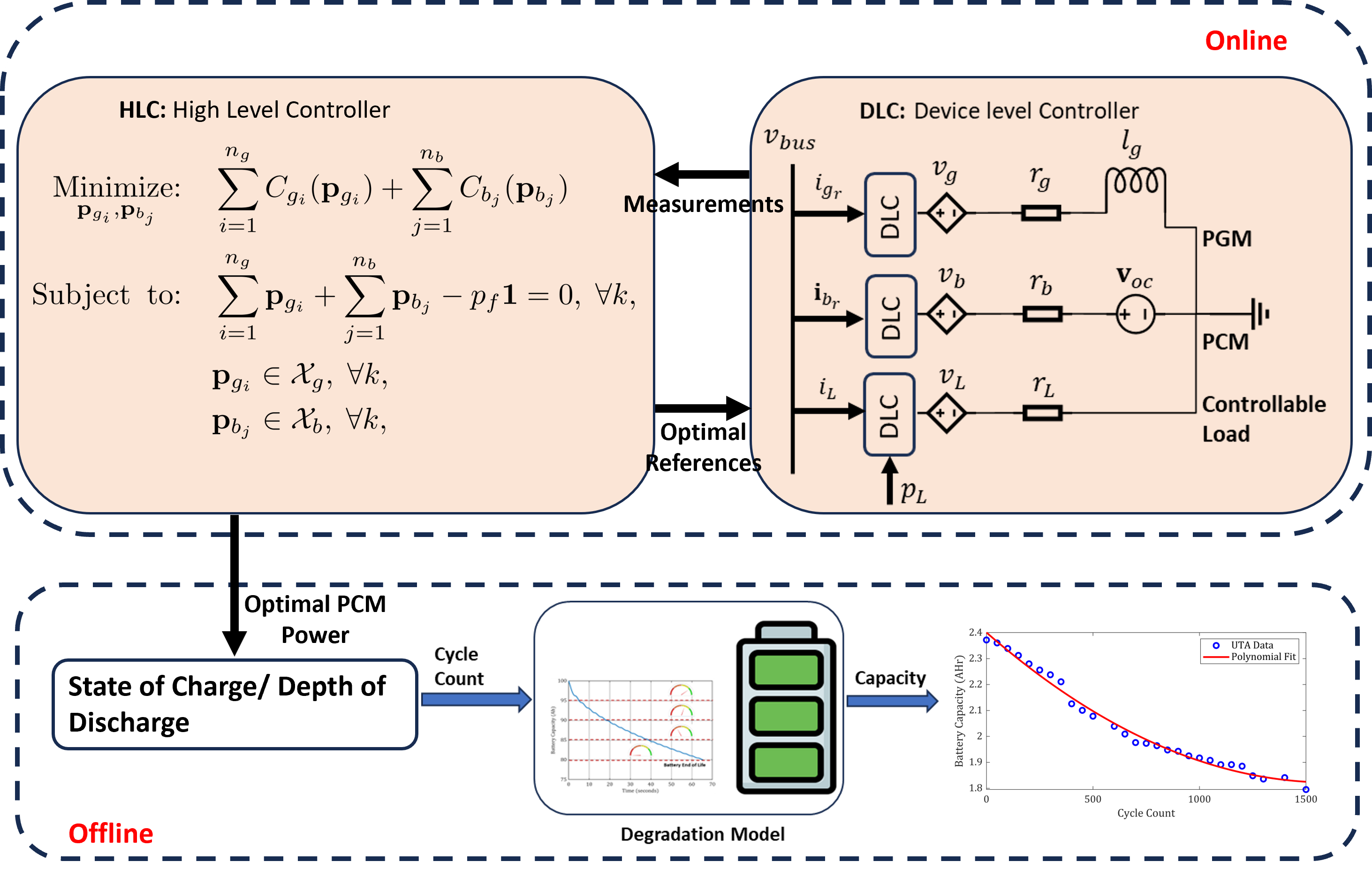}
	 \caption{Overall high-level and low-Level control implementation and offline battery degradation calculation}
     \label{Optimization Framework} 
\end{figure}

In Figure \ref{Gen_Power vs Weights} an experiment similar to the previous one is performed, and the energy of the generator is captured versus the weights $\gamma$ and $\beta$. The trends observed for the generator are opposite to the ones observed for the battery. Mathematically this is what is expected to maintain the power supply and demand equality constraint. The \emph{solid lines} depicts the effect of the weight $\gamma$ tuned from 0-10 while keeping $\beta=1$ and then at $\beta=2$ and $3$. When $\gamma=0$, the objective functions goal is minimizing $\frac{\beta}{2}\norm{\mathbf{p}_g-\mathbf{p}_g^r}_2^2$. 
It can be seen that in this case, the energy of the generator is around 6830\textsf{MWHr}. But as the value of the weight $\gamma$ goes to 10, the objective shifts to minimizing the battery usage, and thus, the generator energy goes up to maintain the power balance. Similarly, the \emph{dotted lines} show the effect of weight $\beta$ swept from 0-10 on the generator energy when the weight $\gamma=1$, $\gamma=2$, and $\gamma=3$. 

\begin{figure}[b!]
      \centering
      \includegraphics[width=0.6\textwidth]{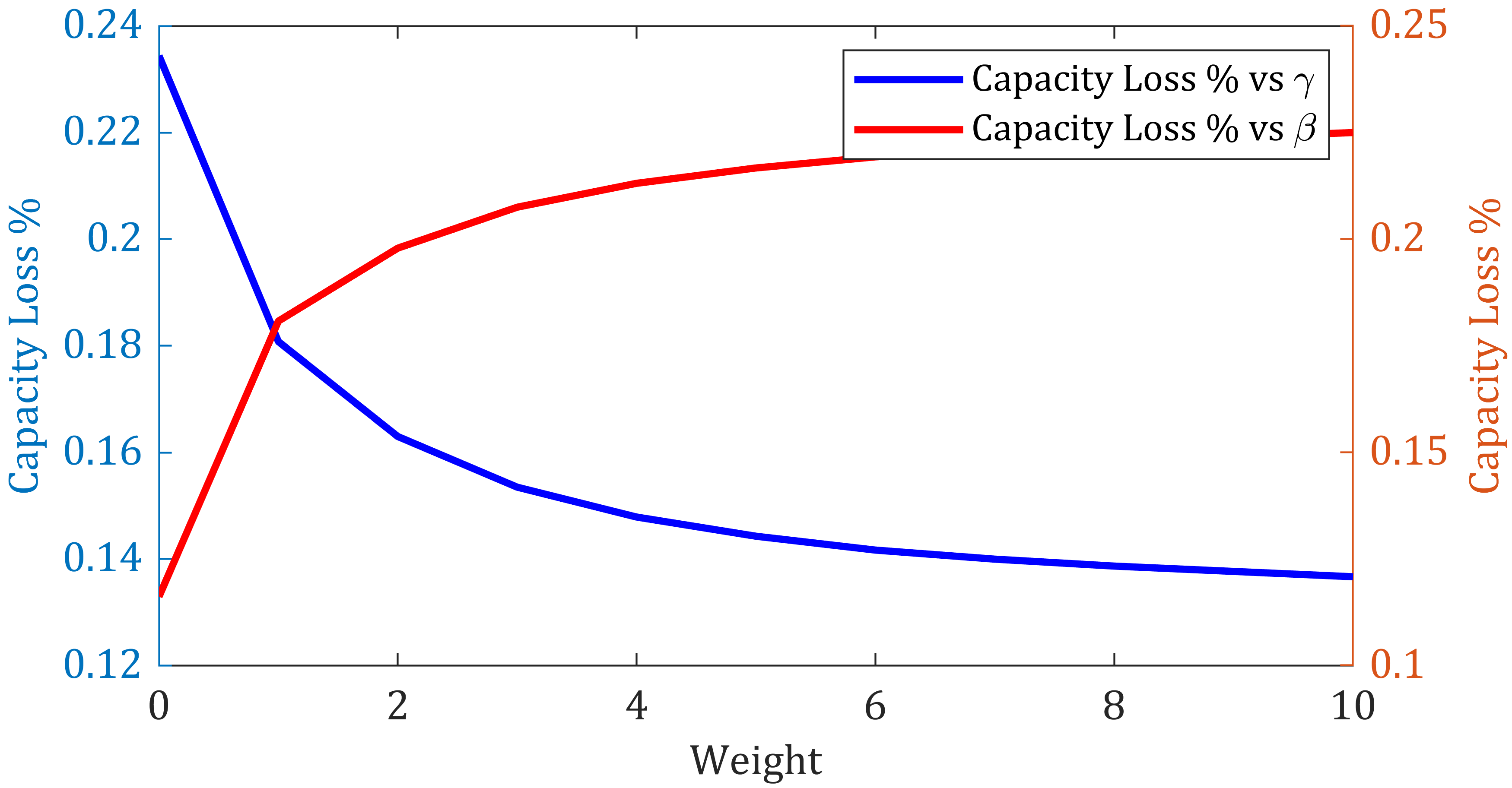}
	 \caption{Battery capacity loss $\%$ vs the weight $\gamma$.}
     \label{Cap_loss vs Beta} 
\end{figure}

Figure \ref{Cap_loss vs Beta} shows the effect of the weight $\gamma$ incremented from 0-10 while keeping the weight $\beta=1$ on the battery capacity loss \%. It can be seen that, as the emphasis on the weight $\gamma$ increases the capacity loss percent decreases. In summary, the total capacity loss decreases with $\gamma$ but increases with $\beta$. From the simulated results, we can conclude that a designer can tune the weights $\beta$ and $\gamma$ to obtain the optimal performance of the generator and the battery, respectively.  

\subsubsection{4-Zone Shipboard Power System}
The shipboard power system model proposed by the U.S. Office of Naval Research (ONR) comprises of PGMs, PCMs, and device level controllers (DLCs) for local voltage and current control. These components are categorized into a zonal structure \cite{ESRDC_1}. PGMs consist of fuel-operated gen-sets. PCMs comprise multiple energy storage elements such as ultracapacitors and battery energy storage systems. All zones are unified through a common 12\textsf{kV} direct current (DC) bus. The power demand must be met under any given circumstance i.e. power supplied by all sources must match the load demand scenarios such as high ramp and uncertainty in the load to maintain the system stability. In this paper, only the main power-supplying sources as mentioned before are considered. There are two main PGMs and two main PCMs. The combined power supplied by these sources must meet the SPS power demand. 

\begin{figure}[h!]
      \centering
      \includegraphics[width=0.6\textwidth]{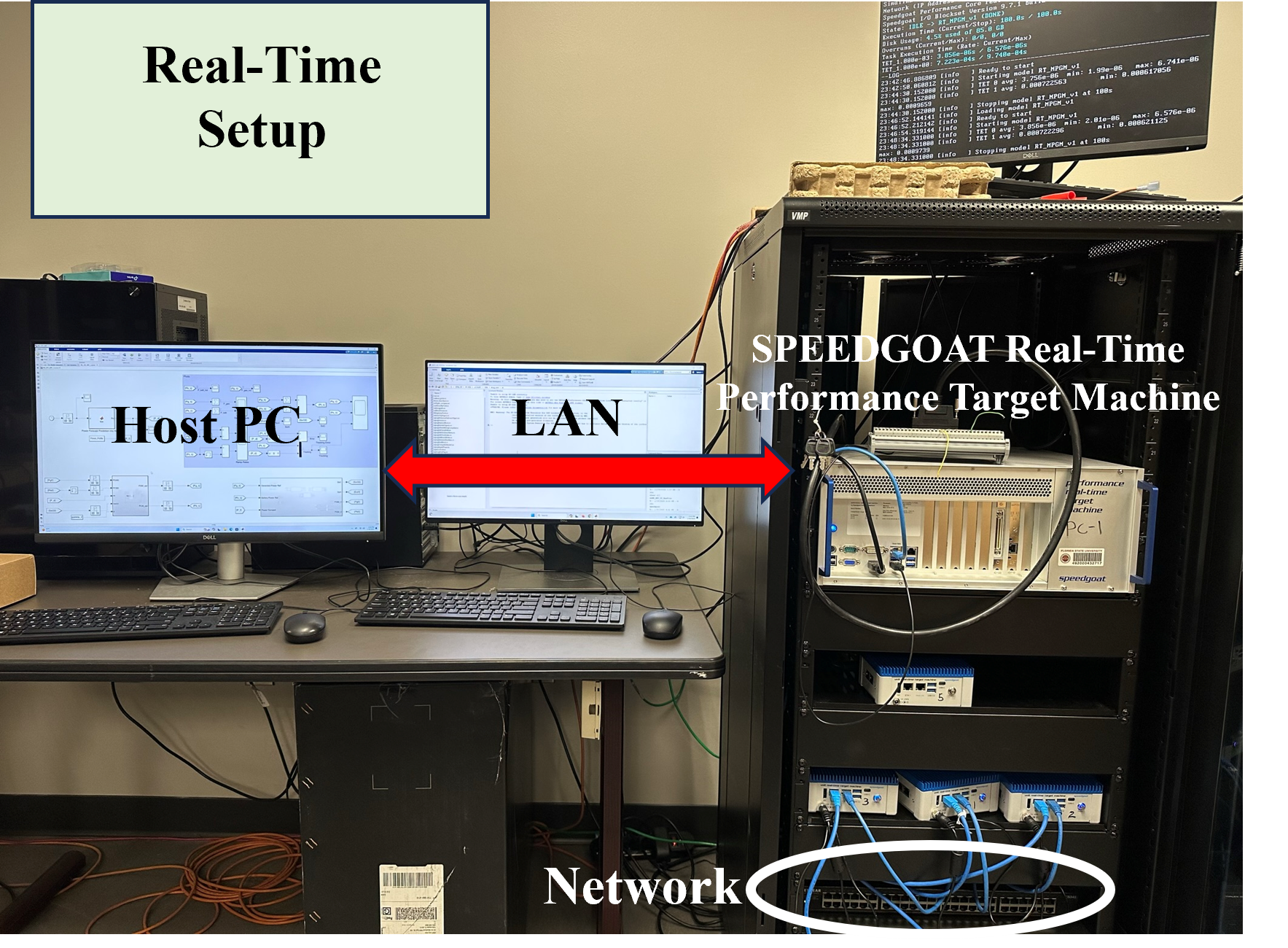}
	 \caption{The real-time setup showing the connection between the host PC and the real-time performance target machine.}
     \label{SIL_setup} 
\end{figure}

\begin{figure}
    \centering
    \includegraphics[width=0.9\textwidth]{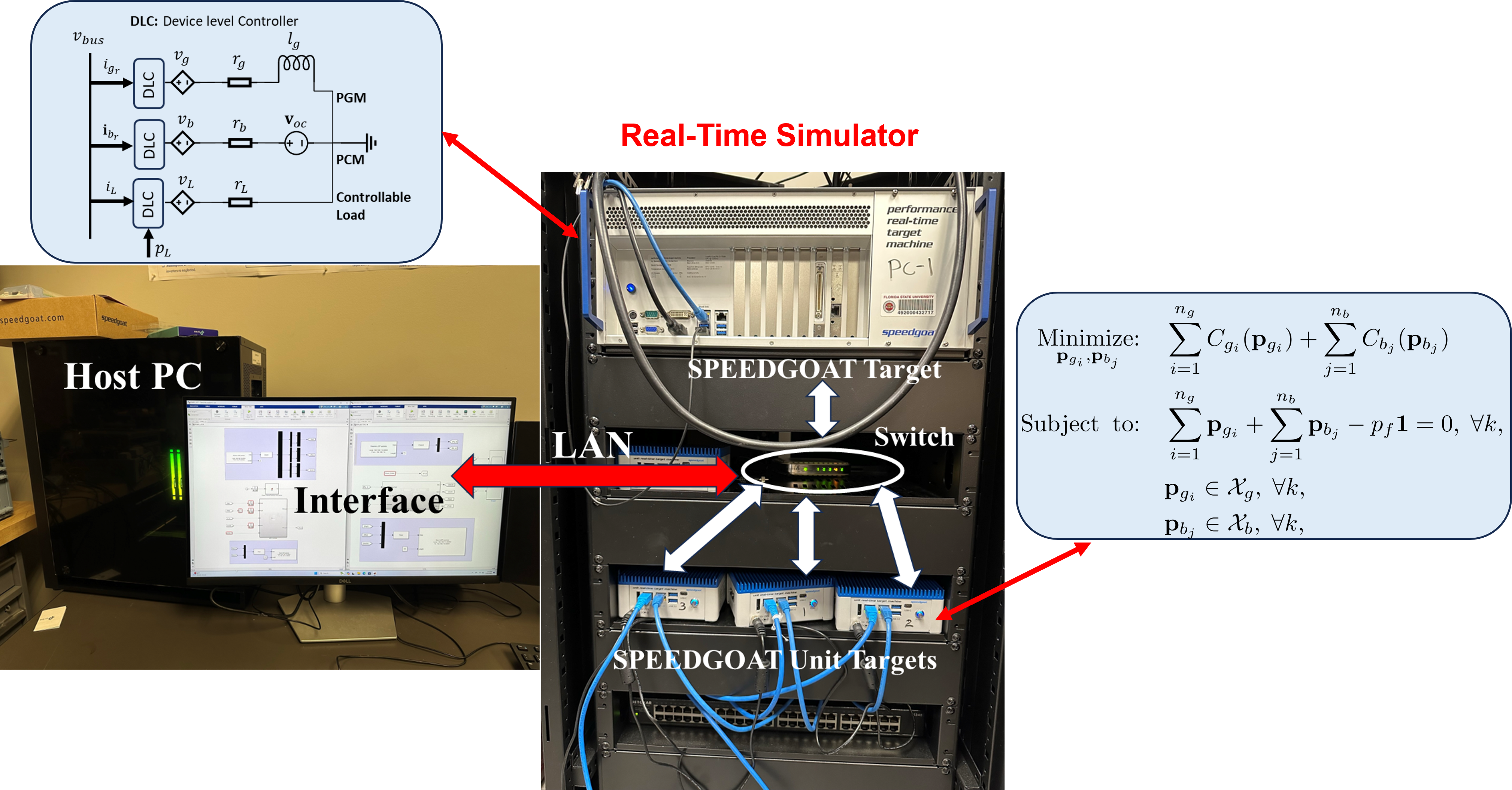}
    \caption{Software-in-Loop setup for multiple PGMs and PCMs test case}
    \label{fig:sil_multiple}
\end{figure}

In this case study, we demonstrate the effectiveness of the developed distributed EM controller by testing it via a \emph{real-time} simulation on an SPS consisting of 2-PGMs and 2-PCMs. Figure \ref{SIL_setup} shows the simulation setup consisting of a host computer and a SPEEDGOAT performance target machine with real-time simulation capabilities connected via a LAN connection (shown in Figure \ref{SIL_setup}, and Figure \ref{fig:sil_multiple}).  The optimization problems in  (\ref{PGM_Node}) and (\ref{PCM_Node}) are implemented at the individual PGMs and PCMs. The simulation time-step was chosen to be 1\textsf{ms}. The time taken by the algorithm to converge was observed to be 0.04\textsf{seconds} and about 30 iterations for the optimal values to be within the tolerance limits. Thus, the MPC and model rate transition in Simulink simulation was set at 1\textsf{second}. The simulation parameters and the model component ratings used in the simulation are provided in Table-\ref{tab:rated} and Table-\ref{tab:rated_4_Zone}.

\begin{table}[t!]
\centering
\caption{Rated values and simulation parameters for 4-zone SPS model \cite{ESRDC_1}}
\label{tab:rated_4_Zone}
\resizebox{0.6\columnwidth}{!}{\begin{tabular}{c|c}
\hline \hline
\textbf{Parameter/Simulation Description}    & \textbf{Parameter Value}  \\
 \hline \hline
Rated Power PGM-1    & 29 MW   \\
Rated Power PGM-2    & 29 MW \\
Desired Operating Power PGMs ($\mathbf{p}_{g_i}^r$) & 15 MW \\
Rated Power PCM-1    & $\pm$ 10.64 MW \\
Rated Power PCM-2    & $\pm$ 10.64 MW \\
Rated Bus Voltage ($v_{bus}$) & 12 kV  \\
Capacity of PCMs ($Q_T$)   & 20 AHr \\
\hline
PGMs weight ($\beta_i$) & 1 \\
PCM-1 weight ($\gamma_1$) & 1 \\
Gradient Step ($\alpha$)  & 0.1 \\
 \hline
\end{tabular}}
\end{table}

\begin{figure}[b!]
      \centering
      \includegraphics[width=0.6\textwidth]{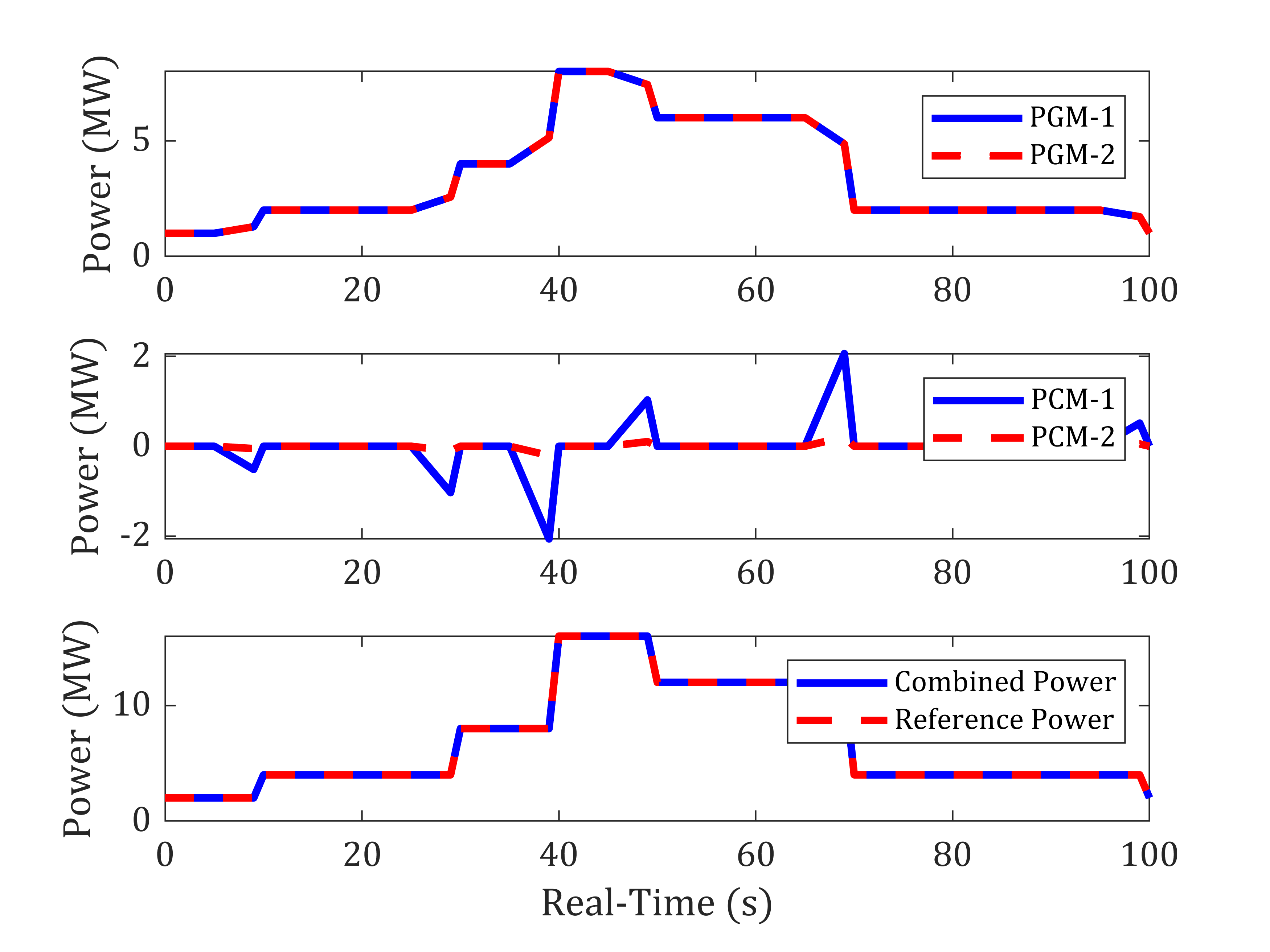}
	 \caption{Real-Time simulation on SPEEDGOAT target machine showing the battery degradation aware power split for the designed model predictive energy management.}
  \label{RT_Power_split} 
\end{figure}

\begin{figure}[h!] 
	\centering
	\includegraphics[width = 0.6\textwidth]{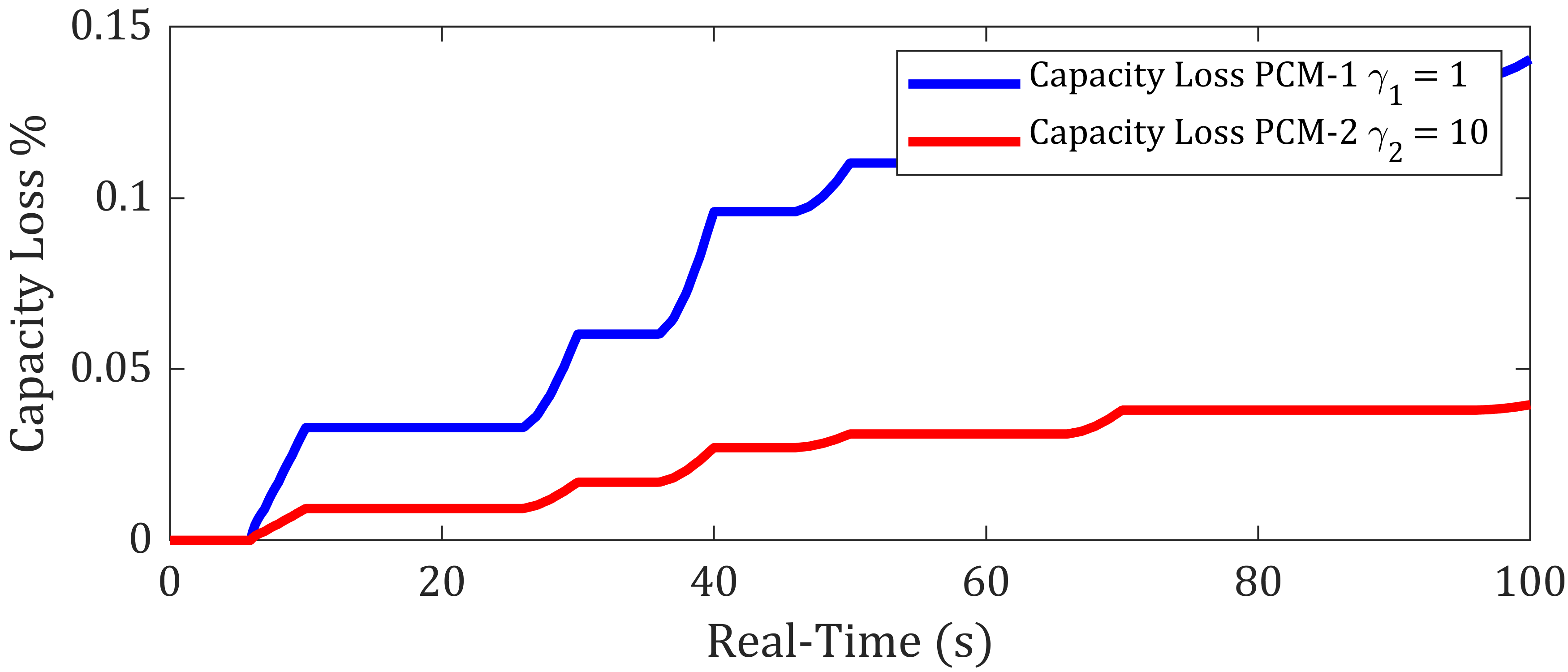} 
	\caption{Real-Time simulation on SPEEDGOAT target machine showing the capacity loss \% of the PCM-1 and PCM-2 for $\gamma_1$ and $\gamma_2$.
    }
	\label{RT_SoHs}
\end{figure}

\begin{figure}[h!] 
	\centering
	\includegraphics[width = 0.6\textwidth]{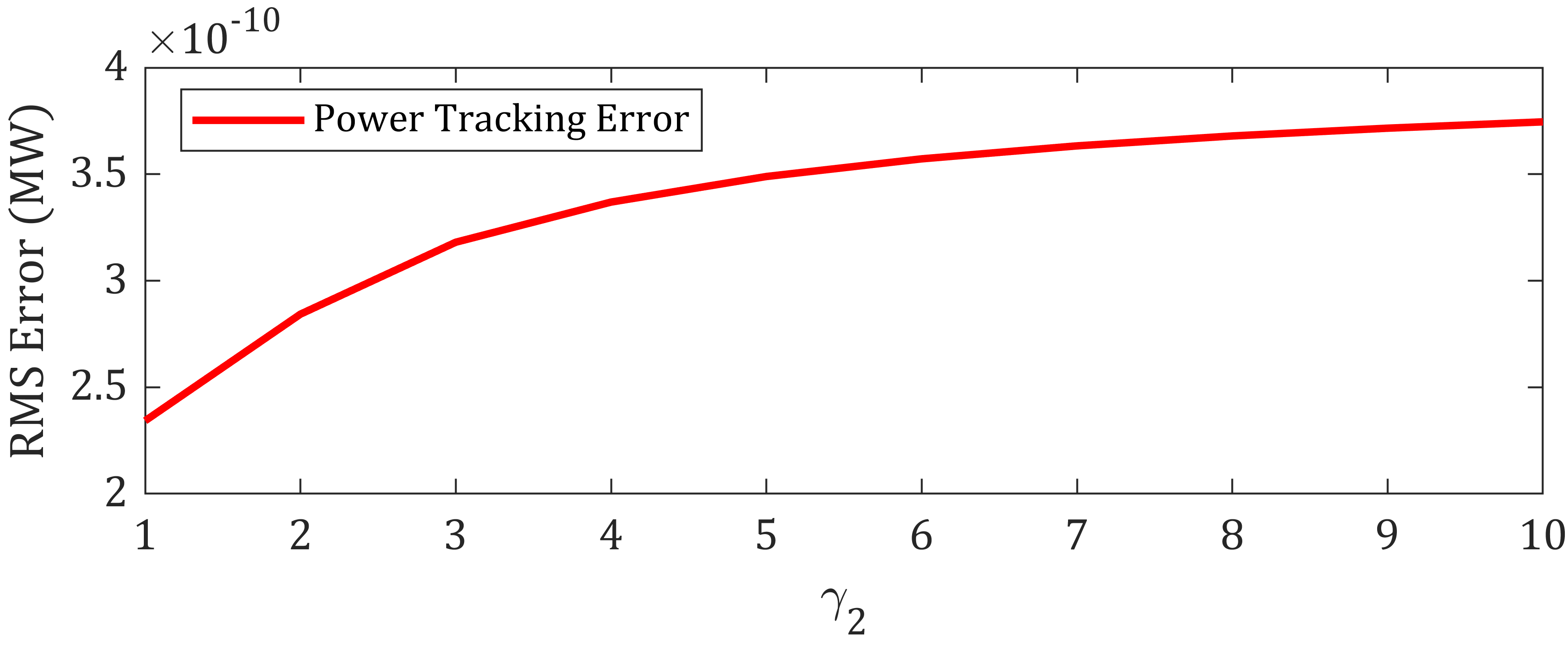} 
	\caption{The root mean square (RMS) Power Tracking Error and the impact of the weight $\gamma_2$.
    }
	\label{RMS_Power_Error_vs_Gamma}
\end{figure}

\begin{figure}[h!] 
	\centering
	\includegraphics[width = 0.6\textwidth]{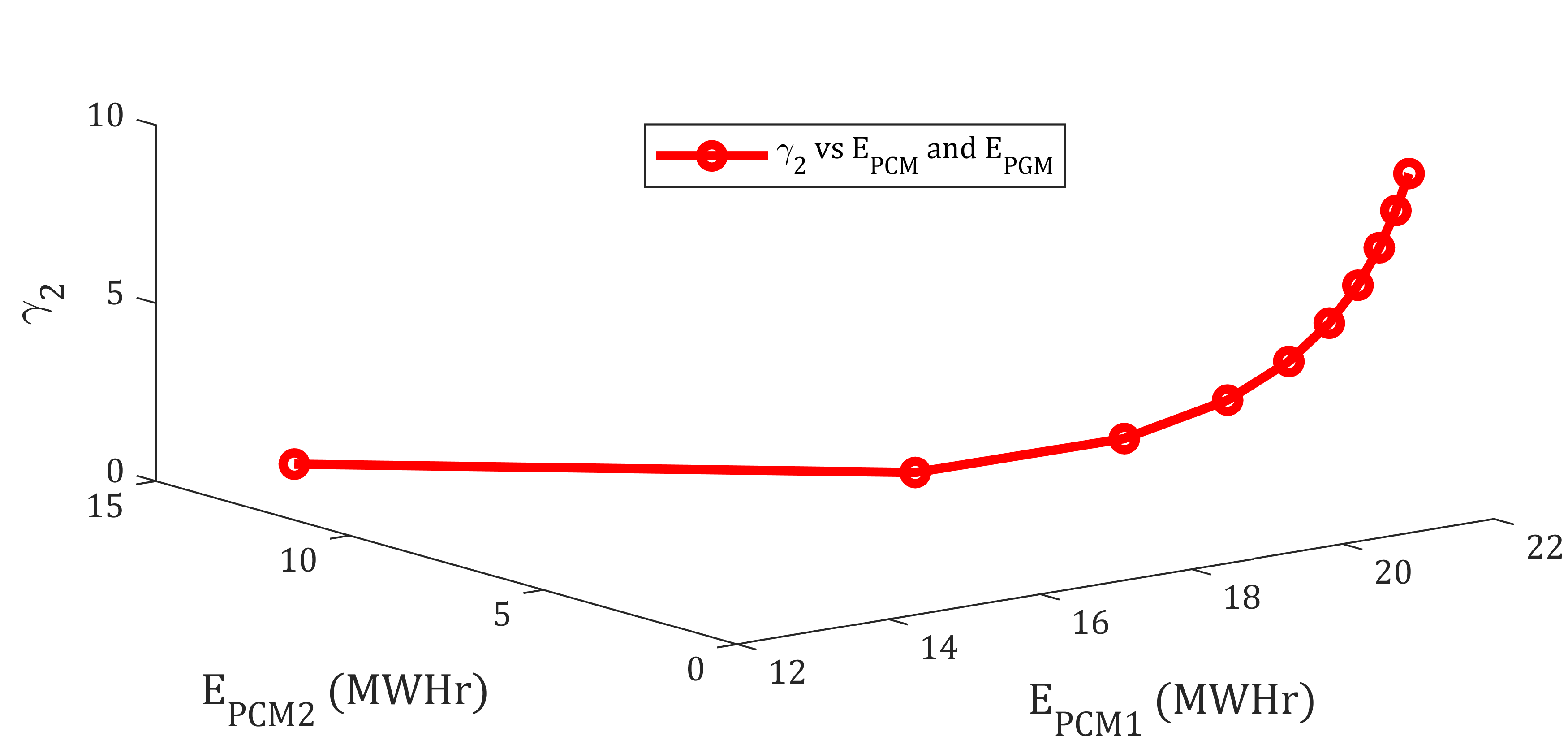} 
	\caption{The impact of the weight $\gamma_2$ keeping $\gamma_1 = 1$ on the energy of the PCM-1 and PCM-2.
    }
	\label{Gamma_vs_PCMs}
\end{figure}

Figure \ref{RT_Power_split} shows the power split between the PGMs, PCMs, and the tracked load simulated on a SPEEDGOAT real-time performance target machine. Figure \ref{RT_SoHs} shows the effect of the weight $\gamma_1$ and $\gamma_2$ on the PCM capacity loss \%. To study the effect of different $\gamma_2$ on the PCM-2 capacity loss the value of $\gamma_2$ is tuned from 0-10. Figure \ref{RMS_Power_Error_vs_Gamma} shows the RMS power tracking error for each simulation for different values of $\gamma_2$ keeping $\gamma_1,\beta_1,\beta_2$ as constants. It can be seen that the error is in the order of $10^{-10}$ which is a numerical error and is almost negligible. Since the main focus of this work is on studying the effect of the weights $\beta_i$ and $\gamma_j$ on the EM, we present the results emphasizing the trade-off between the PGMs and the PCMs. The following assumptions are made on the simulation setup to facilitate the result comparisons. The prediction horizon ($h$) is fixed at 5-time steps (\textsf{seconds}) for all the simulation scenarios. The weights on the PGMs $\beta_i$ are assumed to be constant for all the simulation cases. The weight of the PCM-1 ($\gamma_1 = 1$) is fixed for all the scenarios. The weight associated with PCM-2 ($\gamma_2$) is tuned from $[1-10]$ and the observations are presented capturing the impact of the weighting of the PCMs on the total energy supplied by the PCMs. The ramp rates for the PGMs ($r_g$) are assumed to be $10\%$ of the rated power of PGMs and the ramp rates for the PCMs ($r_b$) are assumed to be $95\%$ of the rated power of PCMs.

\begin{figure}[h!] 
	\centering
	\includegraphics[width = 0.6\textwidth]{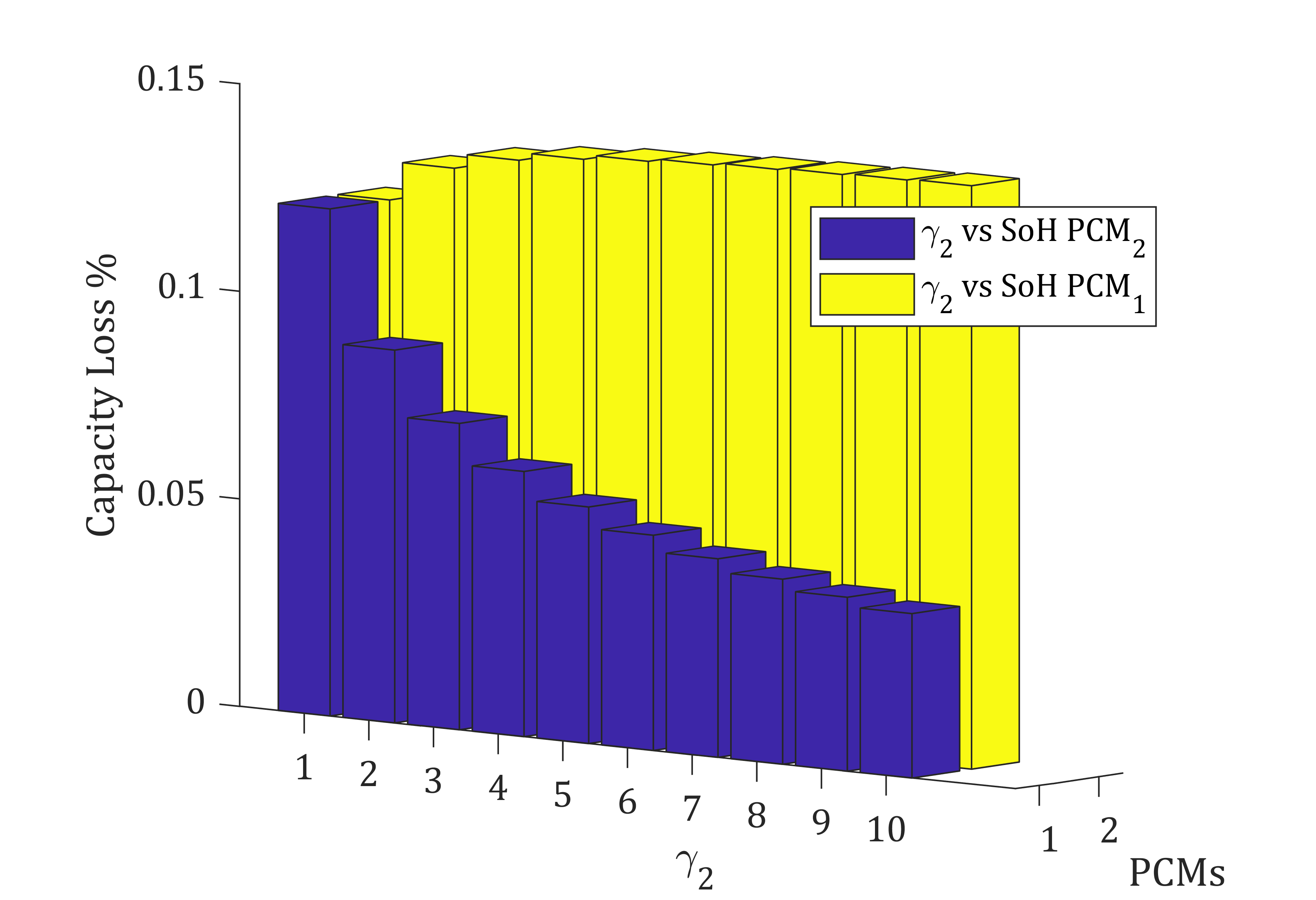} 
	\caption{The impact of the weight $\gamma_2$ on the SoH of the PCM-2.
    }
	\label{Gamma_vs_SoH}
\end{figure}

Figure \ref{Gamma_vs_PCMs} shows the impact of the weight $\gamma_2$ keeping the weight $\gamma_1 = 1$ on the energies of the PCMs in \textsf{MWHr}. The energy of the PCM-1 is plotted on the X-axis, the energy of the PCM-2 is plotted on the Y-axis, and the weight $\gamma_2$ is plotted on the Z-axis. It can be seen that, as $\gamma_2$ increases the energy of the PCM-2 decreases, but to meet the other constraints, the energy of the PCM-1 increases. Figure \ref{Gamma_vs_SoH} shows the effect of the weight $\gamma_2$ keeping $\gamma_1 = 1$ on the capacity loss \% of the PCMs. The X-axis represents the number of PCMs. The Y-axis represents the weight $\gamma_2$ and the Z-axis represents the capacity loss \% of the PCMs. It can be seen that as $\gamma_2$ increases the capacity loss \% of PCM-2 decreases at the expense of the increase in the capacity loss \% of PCM-1. Thus, the effect of the weighting $\gamma_2$ can be observed. 

\begin{figure}[h!] 
	\centering
	\includegraphics[width = 0.6\textwidth]{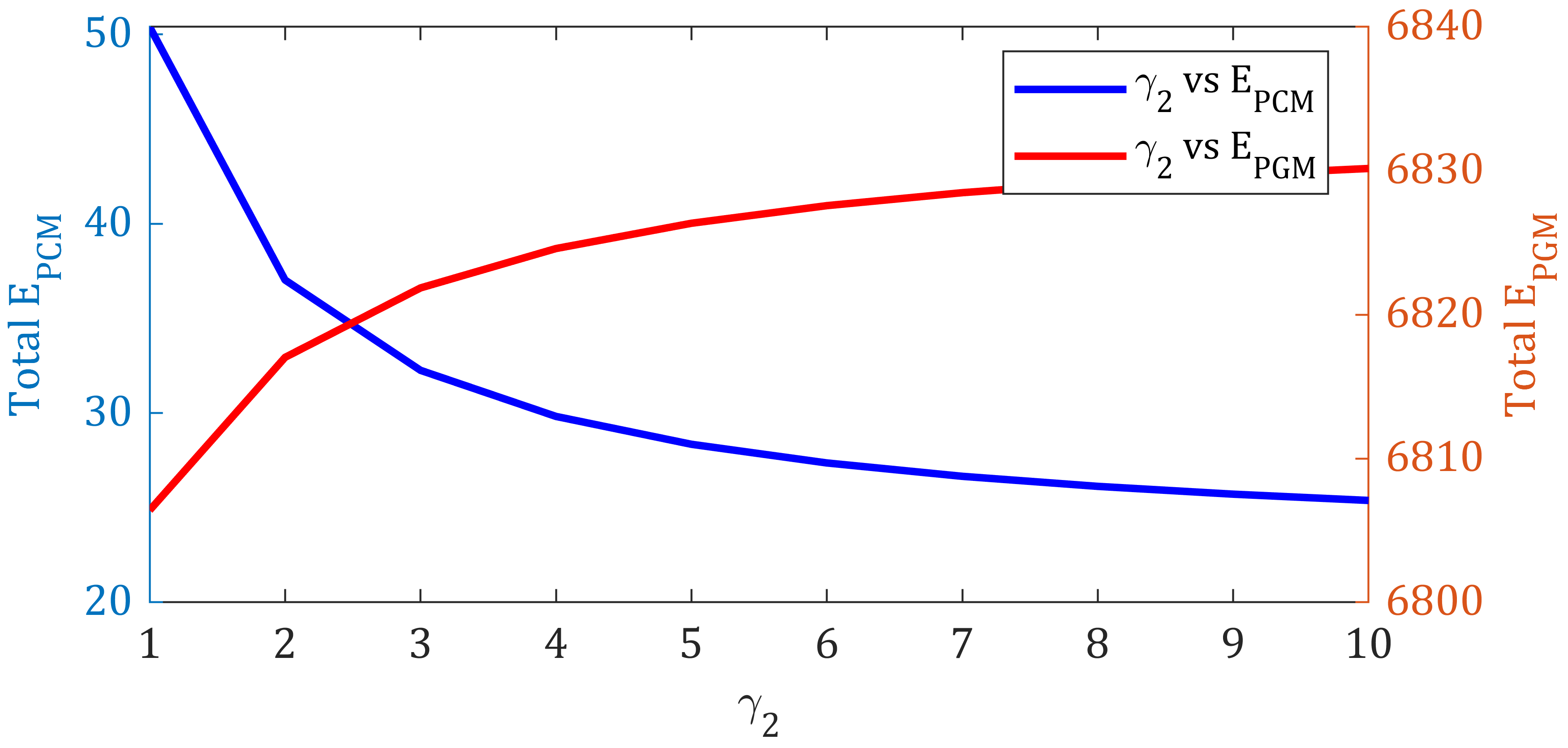} 
	\caption{The impact of the weight $\gamma_2$ keeping $\gamma_1 = 1$ and $\beta = 1$ versus the total energy of the PCMs and PGMs.
    }
	\label{Gamma_vs_Energies}
\end{figure}

Figure \ref{Gamma_vs_Energies} shows the impact of the weight $\gamma_2$ on the total PCMs and PGMs energies. The X-axis represents the weights $\gamma_2$. The left Y-axis represents the total energy of both the PCMs and the right Y-axis represents the total energy of both the PGMs. It can be inferred that, as the weight $\gamma_2$ increases the total energy of the PCMs decreases, to maintain the power balance and to satisfy the power balance constraint, the total energy of the PGMs increases. Thus, from the results, the effectiveness of the designed algorithm and the impact of the weight tuning in battery degradation management while operating the generators at a rated value can be inferred.

\subsubsection{Conclusion} 

A distributed MPC-based energy management strategy considering energy storage degradation is presented for a shipboard power system as a numerical case study. The generator, battery ramp rate limitations, generator-rated conditions, and pulsed power load conditions are considered. Absolute power extracted from, the battery as a model-based battery degradation heuristic capturing the battery usage is proposed and used in the optimization. The trade-off in the power-sharing and the energy exchange for different optimization weights is presented.

\section{Transmission Admissibility Based Power Delivery}

\subsection{Fault-Tolerant Power Delivery}
Integration of Inverter Based Resources (IBRs) which lack the intrinsic characteristics such as the inertial response of the traditional synchronous-generator (SG) based sources presents a new challenge in the form of analyzing the grid stability under their presence. While the dynamic composition of IBRs differs from that of the SGs, the control objective remains similar in terms of tracking the desired active power. This letter presents a decentralized primal-dual-based fault-tolerant control framework for the power allocation in IBRs. Overall, a hierarchical control algorithm is developed with a lower level addressing the current control and the parameter estimation for the IBRs and the higher level acting as the reference power generator to the low level based on the desired active power profile. The decentralized network-based algorithm adaptively splits the desired power between the IBRs taking into consideration the health of the IBRs transmission lines. The proposed framework is tested through a simulation on the network of IBRs and the high-level controller performance is compared against the existing framework in the literature. The proposed algorithm shows significant performance improvement in the magnitude of power deviation and settling time to the nominal value under faulty conditions as compared to the algorithm in the literature.

\subsubsection{Background}
Electrical Power Systems (EPS) are undergoing substantial structural and operational transformation by replacing the existing Synchronous Machines (SMs) with the Inverter Based Resources (IBRs). In contrast to the SMs, the IBRs exhibit low inertia, resulting in a swift response to stochastic events. However, this presents considerable control challenges related to the stability and robustness of the EPS. In the existing literature, several researchers have proposed methods for controlling large-scale IBRs such as discrete-time consensus control using proportional derivative PD \cite{chen2020distributed}, and distributed droop control \cite{schiffer2015voltage}. However, these approaches lack robustness. A distributed model predictive control (MPC) for droop-controlled IBRs was proposed by \cite{anderson2019distributed}. Since this approach is based on MPC, it is computationally heavy. A distributed sliding mode control (SMC) for islanded AC microgrids was presented in ~\cite{alfaro2021distributed}. SMC-based control techniques show great robustness against exogenous disturbances but surfers from chattering effects. Also, observers are often required for disturbance estimation and chattering reduction which may increase with switching gain in inverters thereby increasing the complexity of implementing SMC \cite{9829024}.

Moreover, in a power system network proliferated with IBRs, there is a possibility of losing any of the distributed energy resources (DERs) because of poor response to time-varying load changes or faults. Thus, recovering the aggregated total output power from collective IBRs may be an issue. This presents a serious control problem. In \cite{AMELI}, a robust adaptive control technique is used to track the aggregate output power of the IBRs when one of them is lost. The result gave a perfect output power tracking that matches the total power contributed by each unit of IBRs.  A dispatchable virtual oscillator control is designed for a network of IBRs to track the desired active power, and voltage magnitude in~\cite{subotic2020lyapunov}. The method provides sufficient conditions for voltage stability but does not give the admissible set for desired powers. 

A renowned method for maintaining power sharing in the power grids dominated by IBRs is the droop control technique. The control algorithm is saddled with voltage and frequency deviation issues emanating from uncertainties in the output impedance and poor transient performance. Traditional droop control is very effective in systems with resistive output impedance which results in poor grid stability \cite{en15124439}. In as much as droop control techniques have been modified to tackle these issues, tracking the cumulative output power resulting from the loss of an inverter or the poor output impedance condition is still an open problem. Since the contemporary power droop control method degrades with the line impedance, modified droop control has been developed to improve the active power sharing among the grid-connected inverters \cite{8660486}. 

Nevertheless, the certificate of stability in the presence of the disturbances and the active power tracking is not guaranteed to employ the droop methods. Few authors have proposed a robust control-based approach to address stability and power quality issues arising from the integration of IBRs \cite{Anubi, faiz2020h}. An optimization-based approach was proposed in \cite{Bidram}. However, the active power tracking under IBR failure was not discussed. In \cite{CHANG201685}, the authors proposed a distributed control for IBR coordination in the islanded microgrids. However, the impact of the failure of the IBR was not discussed. Thus, this paper focuses on the tracking of the aggregated active power under IBR failure in a decentralized fault-tolerant framework. The main contributions in the paper are:
\begin{enumerate}
    \item A decentralized fault-tolerant high-level control is proposed that adaptively splits the active power among the IBRs in the presence of faults.
    \item An adaptive estimator for the parameter estimation is designed with the stability analysis. 
    \item A faster response in power shared during faulty conditions is demonstrated compared to the previously proposed approach.
    \item Improvement in the active power tracking under faults is demonstrated compared to previously proposed work.
\end{enumerate}

The section is organized as follows: Section-\ref{Model}, the inverter model is presented followed by the low-level and high-level control design in Section-\ref{Control}. Section-\ref{Sim} presents the simulation results of the designed controller and the results are compared with the results in the literature.

\begin{figure}[h!] 
\centerline{\includegraphics[width=0.5\textwidth]{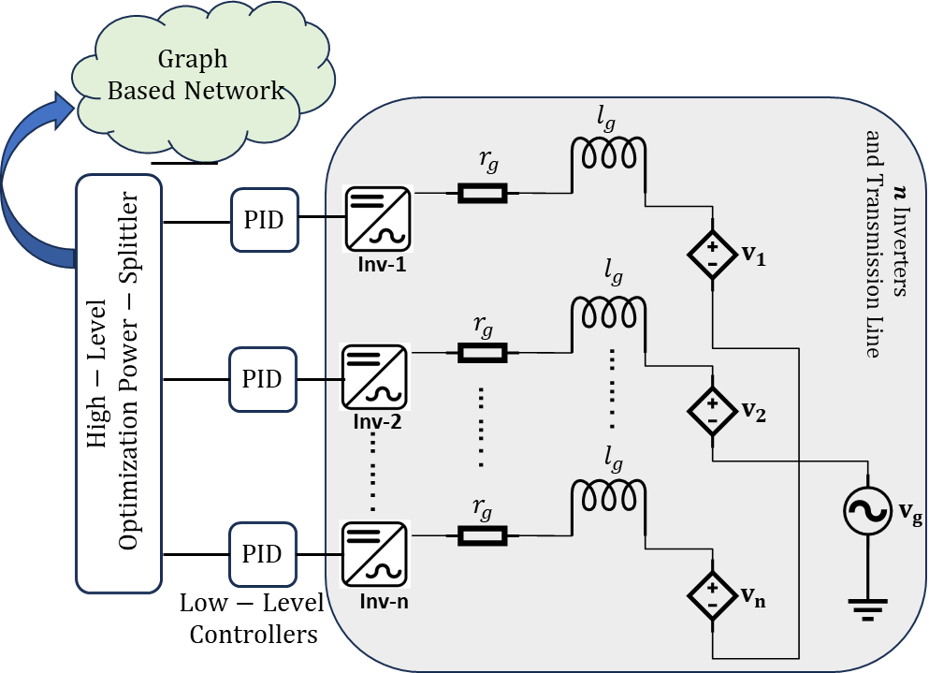}}
\caption{Hierarchical optimization based control design for IBRs}
\label{fig_example25.jpg}
\end{figure}

\subsection{Transmission Model in dq0}\label{Model}
The dynamics of a \emph{single} grid following IBR model connected to the grid via a transmission line in $dq$ coordinates is given as follows \cite{levron}:
\begin{equation}\label{Inverter_Model}
    l_g\frac{d\mathbf{i}}{dt} = -(r_g I-l_g \omega_g J)\mathbf{i}(t) + \mathbf{v}(t) -\mathbf{v}_g,
\end{equation}
where, $\mathbf{i}(t) \in \mathbb{R}^2$, $\mathbf{v}(t) \in \mathbb{R}^2$ and $\mathbf{v}_g \in \mathbb{L}_2$ are the inverter current, voltage (control input) and measured grid voltage respectively. $r_g \in \mathbb{R}_+$ and $l_g \in \mathbb{R}_+$ is the resistance and the inductance of the transmission line connecting to the grid in \textsf{Ohm} and \textsf{Henry}. $\omega_g \in \mathbb{R}_+$ is the grid frequency in \textsf{rad/s}. 
\begin{assumption} 
For simplicity of exposition, the grid frequency $\omega_g$ is considered to be fixed. Also, the switching dynamics of the inverter are neglected \cite{Anubi}. The grid voltage $\mathbf{v}_g$ is regulated to a fixed value since the IBRs are designed in a grid-following mode.
\end{assumption}

\subsection{Fault-Tolerant Control Development}\label{Control}
The overall control development model \cite{vedula2024faulttolerantdecentralizedcontrollargescale} consists of the high-level aggregated power tracking controller acting as a reference current generator for the low-level controller as shown in Figure \ref{fig_example25.jpg}. 
\subsubsection{Low-Level Control}
The objective of the low-level control design is to track a given reference power signal asymptotically while learning the fault indicating uncertain model parameters. Since the voltage is assumed to be regulated, the power tracking problem is converted to a current tracking problem. Namely, given desired active and reactive power pair $(p_{ref},q_{ref})$, the desired reference current is given by
\begin{equation}
    \mathbf{i}_{ref} = \frac{\mathbf{v}_g}{\norm{\mathbf{v}_g}_2^2}p_{ref}+J\frac{\mathbf{v}_g}{\norm{\mathbf{v}_g}_2^2}q_{ref}.
\end{equation}
To achieve the adaptive current tracking objective, we consider the \emph{fault-free} reference model 
   \begin{equation}\label{Model Reference}
   l_m \frac{d\mathbf{i}_m}{dt}=-\bigg((r_m+k)I-l_m\omega_gJ\bigg)\left(\mathbf{i}_m(t)-\mathbf{i}_{ref}(t)\right),
   \end{equation}  
   where $\mathbf{i}_m(t) \in \mathbb{R}^{2}$, $\mathbf{i}_{ref}(t) \in \mathbb{R}^{2}$ is the desired and bounded reference, $r_m >0, l_m >0$ are the known model reference parameters, $k>0$ is a convergence rate tuning parameter. Consequently, the model reference active and reactive powers, $p_m\triangleq\mathbf{v}_g^\top\mathbf{i}_m$ and $q_m\triangleq\mathbf{v}_g^\top J\mathbf{i}_m$, satisfy
   \begin{align}
       \left\|\left[\begin{array}{c}p_m(t)-p_{ref}\\q_m(t)-q_{ref}\end{array}\right]\right\|_2 \le \left\|\left[\begin{array}{c}p_m(0)-p_{ref}\\q_m(0)-q_{ref}\end{array}\right]\right\|_2e^{-\frac{r_m+k}{l_m}t},
   \end{align}
which shows an exponential tracking of the power references with a rate that can be tuned by the control gain $k$.

Next, consider the current tracking error, between the system dynamics in \eqref{Inverter_Model} and the model reference in \eqref{Model Reference}, as follows:
\begin{equation}\label{measurement_tracking_2}
    \tilde{\mathbf{i}}(t) = \mathbf{i}(t)-\mathbf{i}_m(t).
\end{equation}
Taking the first time derivative yields the error dynamics
\begin{equation}
\begin{aligned}
    l_g\dot{\tilde{\mathbf{i}}}(t) = -(r_gI-l_g\omega_gJ)\mathbf{i}(t)+\mathbf{v}(t)-\mathbf{v}_g+((r_m+k)I-l_m\omega_gJ)\mathbf{i}_m(t)-((r_m+k)I-l_m\omega_gJ)\mathbf{i}_{ref}(t), 
\end{aligned}
\end{equation}
adding and subtracting $\hat{r}_g\mathbf{i}(t)$ and setting that the value of the inductance $l_g = l_{g_0} + \Delta l_g$, where $\hat{r}_g$ is an estimate of $r_g$ to be designed and $l_{g_0}$ is the known nominal inductance value with $\Delta l_g$ the associated parametric uncertainty, yields:
\begin{align*}
\l_g\dot{\tilde{\mathbf{i}}}(t) &= -(r_m+k)\tilde{\mathbf{i}}(t)+\tilde{r}_g\mathbf{i}(t)-(\hat{r}_g-r_m-k)\mathbf{i}(t)+(l_{g_0}+\Delta l_g)\omega_g J \mathbf{i}(t)-l_m\omega_g J\mathbf{i}_m(t)-\mathbf{v}_g-\\&\hspace{9.5cm}((r_m+k)I-l_m\omega_gJ)\mathbf{i}_{ref}(t)+\mathbf{v}(t),
\end{align*}
where $\tilde{r}_g=\hat{r}_g-{r}_g$. Consider the control law
\begin{equation}
\begin{aligned}
   \mathbf{v}(t)=(\hat{r}_g-r_m-k)\mathbf{i}(t)-l_{g_0} \omega_g J \mathbf{i}(t)+l_m \omega_g J \mathbf{i}_m(t)+ ((r_m+k)I-l_m\omega_gJ)\mathbf{i}_{ref}(t) + \mathbf{v}_g.
\end{aligned}
\end{equation}
Thus, the error system becomes
\begin{equation}\label{closed_loop_error}
   l_g\dot{\tilde{\mathbf{i}}}(t) = -(r_m+k) \tilde{\mathbf{i}}(t)+\tilde{r}_g\mathbf{i}(t) + \Delta l_g\omega_g J \mathbf{i}(t), 
\end{equation}
with the associated update laws given in the next result
\begin{theorem}
   Given $\varepsilon>0$ and $\bar{r}>0$. Consider the error dynamics in (\ref{closed_loop_error}). If the parameter estimates satisfy the update law
   \begin{equation}\label{update_laws}
   \dot{\hat{r}}_g=\gamma_r \textsf{Proj}_f\left(\hat{r}_g,-\tilde{\mathbf{i}}(t)^\top\mathbf{i}(t)\right),\hspace{2mm}|\hat{r}_g(0)|<\bar{r},
   \end{equation}
   where $\gamma_r > 0$ is the associated adaptation rate for the parameter $r_g$ and $f(r) = \frac{r^2-\bar{r}^2}{2\varepsilon\bar{r}+\varepsilon^2}$, then the origin of the closed-loop system in (\ref{closed_loop_error}) is globally asymptotically stable.  
\end{theorem}
\begin{proof}
    We drop the explicit usage of time for the simplicity of writing the proof. Consider the following Lyapunov candidate function
    $$V(\tilde{\mathbf{i}}, \tilde{r}_g) = \frac{l_g}{2}\tilde{\mathbf{i}}^\top  \tilde{\mathbf{i}} + \frac{1}{2\gamma_r}\tilde{r}_g^2,$$
    taking the first derivative along the time-variables and substituting (\ref{closed_loop_error}) yields
\begin{align*}
    \dot{V} &= \frac{1}{2}\dot{\tilde{\mathbf{i}}}^\top l_g \tilde{\mathbf{i}}+\frac{1}{2}\tilde{\mathbf{i}}^\top l_g \dot{\tilde{\mathbf{i}}}+\frac{1}{\gamma_r}\tilde{r}_g\dot{\hat{r}}_g,\\
    &= -\frac{(r_m+k)}{2}\tilde{\mathbf{i}}^\top\tilde{\mathbf{i}}+\frac{\tilde{r}_g}{2}\mathbf{i}^\top\tilde{\mathbf{i}}+\frac{1}{2}\Delta l_g \omega_g \mathbf{i}^\top J^\top \tilde{\mathbf{i}}-\frac{(r_m+k)}{2}\tilde{\mathbf{i}}^\top \tilde{\mathbf{i}}+\frac{\tilde{r}_g}{2}\tilde{\mathbf{i}}^\top \mathbf{i}+\frac{1}{2}\Delta l_g \omega_g \tilde{\mathbf{i}}^\top J \mathbf{i}+\frac{1}{\gamma_r}\tilde{r}_g\dot{\hat{r}}_g,\\
    &= -(r_m+k)\tilde{\mathbf{i}}^\top\tilde{\mathbf{i}} +\tilde{r}_g\tilde{\mathbf{i}}^\top \mathbf{i}+\Delta l_g \omega_g \mathbf{i}^\top(J+J^\top)\tilde{\mathbf{i}}+\frac{1}{\gamma_r}\tilde{r}_g\dot{\hat{r}}_g,
\end{align*}
using the skew-symmetric property  $(J^\top = -J)$ and substituting the update law (\ref{update_laws}), and using the property $\tilde{r}_g\left(\textsf{Proj}_f(\hat{r}_g,-y)+y\right)\le0 \text{ for all } y\in\mathbb{R}$, yields
\begin{equation}\label{Lyapunov_derivative}
\begin{aligned}
    \dot{V} &= -(r_m+k) \norm{\tilde{\mathbf{i}}}^2 .
\end{aligned}    
\end{equation}
$\dot{V}(t)$ is negative semi-definite (NSD) and $V(t)>0$. Thus, $V(t) \in \mathbb{L}_{\infty}$ which implies $\tilde{\mathbf{i}}(t), \tilde{r}_g(t) \in \mathbb{L}_{\infty}$. Since $\mathbf{i}_m(t)$ is assumed to be bounded, it implies that $\mathbf{i}(t) \in \mathbb{L}_{\infty}$. Consequently, the control input $\mathbf{v}(t) \in \mathbb{L}_{\infty}$. Integrating (\ref{Lyapunov_derivative})
\begin{equation}
    \begin{aligned}
        V(\infty)-V(0) \leq -(r_m+k) \int_{0}^{\infty}\norm{\tilde{\mathbf{i}}(t)}^2 dt,
    \end{aligned}
\end{equation}
It follows that $\tilde{\mathbf{i}}(t) \in \mathbb{L}_2$. From the implications $\tilde{\mathbf{i}}(t)$ is uniformly continuous. Thus, invoking Barbalat's lemma \cite{khalil2002nonlinear} it follows that $\tilde{\mathbf{i}}(t)\longrightarrow \underline{\mathbf{0}}$.
\end{proof}
\begin{proposition}\label{prop-1}
    The current tracking error $\tilde{\mathbf{i}}(t)$ is uniformly bounded according to
    \begin{equation}\label{bound_conv}
        \sup_{t\in\mathbb{R}_+}\norm{\tilde{\mathbf{i}}(t)}_2 \leq \norm{\tilde{\mathbf{i}}(0)}_2 + \frac{\overline{r}}{\sqrt{2\gamma_r}} ,
    \end{equation}
\end{proposition}
\begin{proof}
    Consider (\ref{Lyapunov_derivative}), without loss of generality it can be expressed as
    \begin{align*}
        \dot{V} &= -(r_m+k) \norm{\tilde{\mathbf{i}}}^2, \\
         &\leq \frac{-2k}{l_g}V + \frac{k}{l_g \gamma_r}\tilde{r}_g^2
         \leq \frac{-2k}{l_g}V + \frac{k}{l_g \gamma_r}\overline{r}^2,
    \end{align*}
    using the comparison lemma (\cite{khalil2002nonlinear}), it follows that
    \begin{align*}
        V(t) &\leq e^{\frac{-2k}{l_g}t}V(0)+\frac{\overline{r}^2}{2\gamma_r}\bigg(1-e^{\frac{-2k}{l_g}t}\bigg), \\
        \norm{\tilde{\mathbf{i}}(t)}^2 &\leq e^{\frac{-2k}{l_g}t}\norm{\tilde{\mathbf{i}}(0)}^2+\frac{\overline{r}^2}{2\gamma_r}\bigg(1-e^{\frac{-2k}{l_g}t}\bigg),
    \end{align*}
Thus, $\mathbf{i}(t) \in \mathcal{B}\bigg(\mathbf{i}_m(t),\norm{\tilde{\mathbf{i}}(0)} + \frac{\overline{r}}{\sqrt{2\gamma_r}}\bigg)$. 
\end{proof}

From the Proposition \ref{prop-1}, the resulting power from the \emph{single} IBR satisfies

\begin{align*}
    p(t) &= \mathbf{v}_g(t)^\top \mathbf{i}(t) \in \mathcal{B}\bigg(p_{ref},\norm{\mathbf{v}_g}\left(\norm{\tilde{\mathbf{i}}(0)}+\frac{\overline{r}}{\sqrt{2\gamma_r}}\right)\bigg), \\
    q(t) &= \mathbf{v}_g(t)^\top J \mathbf{i}(t) \in \mathcal{B}\bigg(q_{ref},\norm{\mathbf{v}_g}\left(\norm{\tilde{\mathbf{i}}(0)}+\frac{\overline{r}}{\sqrt{2\gamma_r}}\right)\bigg).
\end{align*}

The total active power supplied by the $n$ number of IBRs satisfies
\begin{equation*}
   \sum_{i=1}^{n} p_i(t)  \in \mathcal{B}\Bigg(\sum_{i=1}^{n} p_{{ref}_i},\norm{\mathbf{v}_g}\sum_{i=1}^{n}\bigg(\norm{\tilde{\mathbf{i}}_{i(0)}}+\frac{\overline{r}_i}{\sqrt{2\gamma_r}_i}\bigg)\Bigg).
\end{equation*}

Moreover, from the update law in (\ref{update_laws}), it is seen that the parameter error is driven by the current tracking error. Thus, if the system is persistently excited \cite{khalil2002nonlinear}, the healthier systems (having small parametric deviations) will result in less error in the active power tracking. Thus, to minimize uncertainty in the total active power, we consider a resource allocation problem using the parameter deviation error to penalize the power requested from the local ($i^{th})$ IBR. 

\subsubsection{High-Level Control}
The objective of the high-level control is to track the total desired power by allocating more reference to healthier IBRs (those with less parameter deviation). This is achieved through the resource allocation problem:
\begin{align}\label{Main_Opt}
    \Minimize_{p_i} & \sum_{i=1}^{n}\frac{\beta_i(\tilde{r}_{g_i})}{2}f_i(p_i) \hspace{1mm} \SubjectTo  \sum_{i=1}^{n} p_i = p_A\\ \label{Reactive_opt}
    \Minimize_{q_i} & \sum_{i=1}^{n}\frac{\beta_i(\tilde{r}_{g_i})}{2}h_i(q_i) \hspace{1mm} \SubjectTo  \sum_{i=1}^{n} q_i = q_A ,
\end{align}
where, $n\in \mathbb{N}$ are the number of IBRs, $\beta_i(\tilde{r}_{g_i})$ is the line parameter dependent penalty associated with the $i^{th}$ IBR. It is determined based on the fault indication and the status of the transmission line. $f_i:\mathbb{R}\longrightarrow\mathbb{R}_+$ is the cost associated with the active power and $h_i: \mathbb{R}\longrightarrow\mathbb{R}_+$ is the cost associated with the reactive power of the individual IBR operation. It is assumed that the functions $f_i, h_i$ are $\mu$-strongly convex and $L$-smooth. $p_i \in \mathbb{R}, q_i \in \mathbb{R}$ are the active and the reactive power of the individual IBR. $p_A, q_A$ are the desired aggregated active and reactive power that needs to be tracked by the IBRs. 
\begin{remark}
        The optimization problems in (\ref{Main_Opt}) and  (\ref{Reactive_opt}) are equality-constrained problems and a closed-loop solution can be obtained offline. However, they depend on the parameter-based penalty parameter $\beta_i(\tilde{r}_{g_i})$ thus the computation needs to be performed online. Moreover, having a centralized controller acts as a point of failure which leads to instabilities in the grid in case the controller fails. Thus, we decentralize the controller to address the single point of failure.
\end{remark}

\textit{Decentralized Development:} The \textit{Lagrangian} for the optimization problems (\ref{Main_Opt}) and (\ref{Reactive_opt}) is given as follows:
\begin{subequations}
\begin{align}
    \mathcal{L}(p,\lambda) = \sum_{i=1}^{n}\frac{\beta_i(\tilde{r}_{g_i})}{2}f_i(p_i)+\lambda (\sum_{i=1}^{n}p_i-p_A),\\
    \mathcal{L}(q,\nu) = \sum_{i=1}^{n}\frac{\beta_i(\tilde{r}_{g_i})}{2}h_i(q_i)+\nu (\sum_{i=1}^{n}q_i-q_A),
\end{align}
\end{subequations}
where $\lambda, \nu \in \mathbb{R}$ are the dual variables for active and reactive powers. Given $\lambda$ and $\nu$, let $p_i^*$ and $q_i^*$ to be the solution of the optimization problem in (\ref{Main_Opt})-(\ref{Reactive_opt}) and is given as
\begin{subequations}
\begin{align}
    p_i^*(\lambda) = \argmin_{p_i} \bigg(\frac{\beta_i(\tilde{r}_{g_i})}{2}f_i(p_i) + \lambda p_i \bigg), \\
    q_i^*(\nu) = \argmin_{q_i} \bigg(\frac{\beta_i(\tilde{r}_{g_i})}{2}h_i(q_i) + \nu q_i \bigg).
\end{align}
\end{subequations} 

\begin{figure}[t!] 
\centerline{\includegraphics[width=0.5\textwidth]{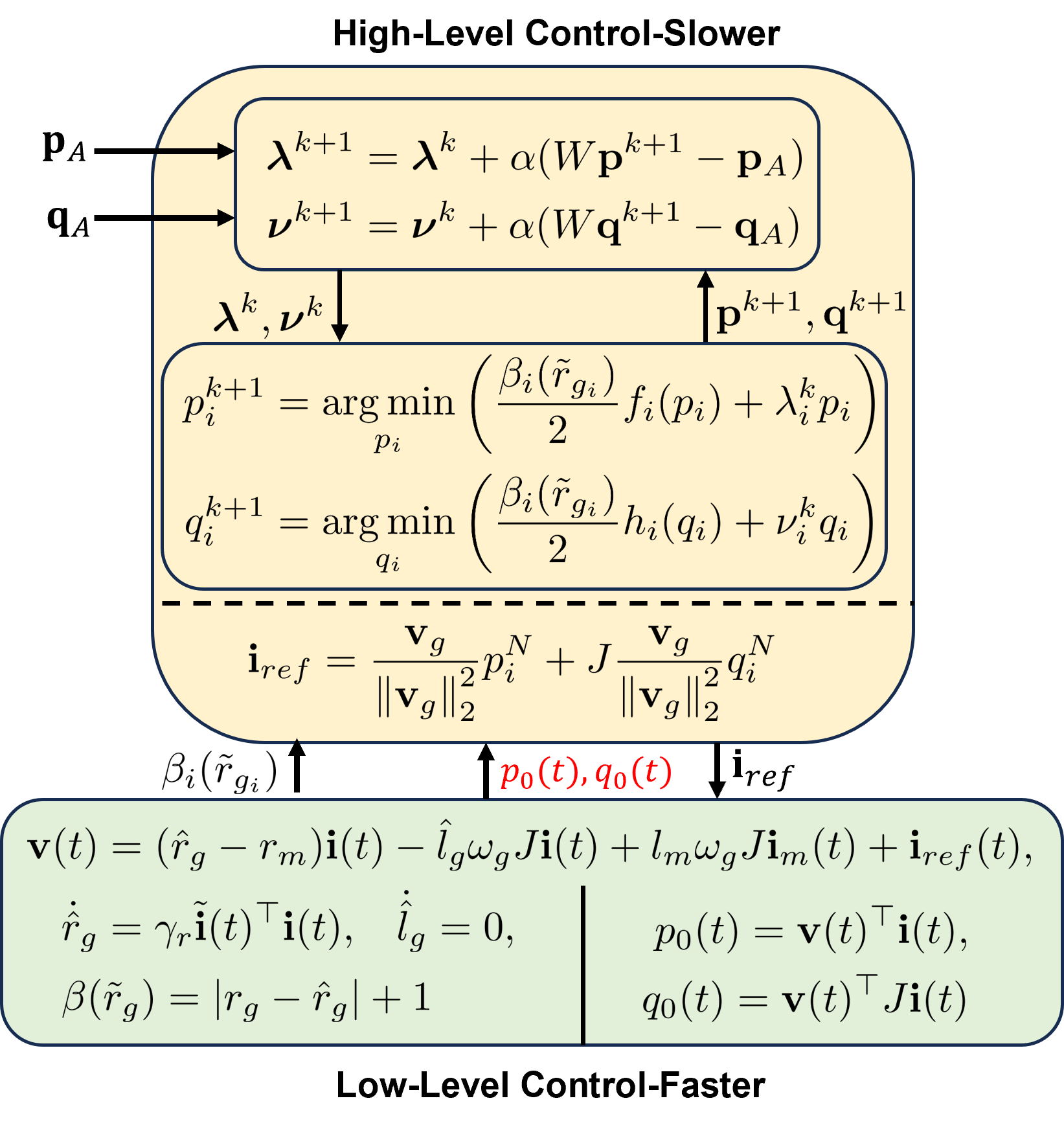}}
\caption{Control architecture}
\label{Control_Architecture}
\end{figure}

Thus, given $p_i^*$ and $q_i^*$, the \emph{Lagrangian dual functions} are defined as the minimum values of the Lagrangian over $p_i$ and $q_i$
\begin{subequations}
\begin{align}
    g(\lambda) \triangleq \inf_{p_i} \mathcal{L}(p_i,\lambda) \equiv \mathcal{L}(p_i^*,\lambda), \\
    d(\nu) \triangleq \inf_{q_i} \mathcal{L}(q_i,\nu) \equiv \mathcal{L}(q_i^*,\nu),
\end{align}
\end{subequations}
Consequently, the dual problems for the active and the reactive power are given as:
\begin{align}\label{Main_Dual}
    \Maximize_\lambda  g(\lambda) && \Maximize_\nu  d(\nu).
\end{align}
Assume the dual problem can be broken down into $n$ sub-problems.  (\ref{Main_Dual}) can be written as:
\begin{subequations}\label{Second_Dual}
\begin{align}
    \Maximize_{\lambda_i=\lambda_j} \quad \sum_{i=1}^{n} g_i(\lambda_i), \hspace{2mm} \forall \{i,j\} \in \mathcal{G} \\
    \Maximize_{\nu_i=\nu_j} \quad \sum_{i=1}^{n} d_i(\nu_i), \hspace{2mm} \forall \{i,j\} \in \mathcal{G}
\end{align}
\end{subequations}
the optimization problem in (\ref{Second_Dual}a-b) are solved over a network $\mathcal{G}$. It is assumed that the network has a \emph{spanning tree}. The information of the reference active and reactive powers $p_A, q_A$ is available only for the first node. Thus, the gradient update step for the $\lambda$ and $\nu$ is given as follows ($k$ here is iteration counter):
\begin{equation}
\lambda_i^{k+1} = \begin{cases} \alpha\sum_{j \in \mathcal{N}_i}w_{ij} \lambda_i^k +  (p_i^* - p_A) \hspace{3mm} \text{for} \hspace{1mm} i=1, \\
    \alpha\sum_{j \in \mathcal{N}_i}w_{ij}  \lambda_i^k + \alpha  p_i^* \hspace{14mm} \text{otherwise},
\end{cases}    
\end{equation}
\begin{equation}
\nu_i^{k+1} = \begin{cases} \alpha\sum_{j \in \mathcal{N}_i}w_{ij} \nu_i^k + \alpha  (q_i^* - q_A) \hspace{3mm} \text{for} \hspace{1mm} i=1, \\
     \alpha\sum_{j \in \mathcal{N}_i}w_{ij} \nu_i^k + \alpha  q_i^* \hspace{14mm} \text{otherwise},
\end{cases}    
\end{equation}
where $\alpha > 0$ is the ascent step size assumed to be a constant for every update step. Consequently, the high-level power splitting decentralized optimization algorithm is given as the following iterative scheme:
\begin{subequations}\label{Network_Opt}
\begin{align}
    p_i^{k+1} &= \argmin_{p_i} \bigg(\frac{\beta_i(\tilde{r}_{g_i})}{2}f_i(p_i) + \lambda_i^k p_i \bigg)  \\
    \boldsymbol{\lambda}^{k+1}  &= W\boldsymbol{\lambda}^k + \alpha (\mathbf{p}^{t+1}-\mathbf{p}_A), \\
    q_i^{k+1} &= \argmin_{q_i} \bigg(\frac{\beta_i(\tilde{r}_{g_i})}{2}h_i(q_i) + \nu_i^k q_i \bigg)  \\
    \boldsymbol{\nu}^{k+1}  &= W\boldsymbol{\nu}^k + \alpha (\mathbf{q}^{k+1}-\mathbf{q}_A),
\end{align}    
\end{subequations}
where $p_i^{k+1} \in \mathbb{R}, q_i^{k+1} \in \mathbb{R}$ are the optimal active and reactive power of the individual IBR, $W \in \mathbb{R}^{n \times n}$ is the graph Laplacian based weigh matrix. $\mathbf{p}^{k+1} \triangleq \left[\begin{array}{cccc} p_1^{k+1}&p_2^{k+1}&\hdots&p_n^{k+1}\end{array}\right]^{\top} \in \mathbb{R}^n$ and \newline $\mathbf{q}^{k+1} \triangleq \left[\begin{array}{cccc} q_1^{k+1}&q_2^{k+1}&\hdots&q_n^{k+1}\end{array}\right]^{\top} \in \mathbb{R}^n$ are the vectors of primal optimal active and reactive power profiles. $\mathbf{p}_A \triangleq \left[\begin{array}{cccc} p_A&0&\hdots&0\end{array}\right]^{\top} \in \mathbb{R}^n$ and $\mathbf{q}_A \triangleq \left[\begin{array}{cccc} q_A&0&\hdots&0\end{array}\right]^{\top} \in \mathbb{R}^n$ are the reference active and reactive power profiles. $f_i$ and $h_i$ are assumed to be $\mu$-strongly convex and $L$-smooth. The primal optimization problems in (\ref{Network_Opt}a) and (\ref{Network_Opt}c) are unconstrained problems and under the assumption that the fixed points $\mathbf{p}^*$ and $\mathbf{q}^*$ exist and the gradient step is a \emph{contraction} the following linear convergence property holds
\begin{subequations}\label{convergence}
\begin{align}
    \norm{\mathbf{p}^{N}-\mathbf{p}^*} \leq L_p^{N-1} \norm{\mathbf{p}_0-\mathbf{p}^*} \\
    \norm{\mathbf{q}^{N}-\mathbf{q}^*} \leq L_q^{N-1} \norm{\mathbf{q}_0-\mathbf{q}^*}
\end{align}    
\end{subequations}
where $N$ is the number of iterations, $L_p, L_q <1 \triangleq max\{\left|1-\alpha \mu \beta\right|,\left|1-\alpha L \beta\right|\}$ are the contraction coefficients. $\mathbf{p}_0(t), \mathbf{q}_0(t)$ are the system active and reactive power measurements acting as the initial values for the optimization problem. The choice of step size $\alpha = 2/(L+\mu \beta)$ minimizes the contraction coefficient. Thus, the condition number defined as $\gamma = L/\mu \beta$ dictates the convergence of the optimization problem and the problem requires at most $\mathcal{O}(\gamma log(\frac{1}{\epsilon}))$ iterations to reach within the $\epsilon$ neighborhood of the primal optimal solutions $\mathbf{p}^*$ and $\mathbf{q}^*$ \cite{boyd2004convex}. The stability of the dual update dynamical systems in (\ref{Network_Opt}b) and (\ref{Network_Opt}d) depends on the weight matrix $W$ and the analysis of it was provided in \cite{Xiao_L}. Given the \emph{optimal} active and reactive powers and the nominal grid voltage, the low-level current reference for $i^{th}$ IBR can be generated based as follows:
\begin{equation}\label{current_ref}
    \mathbf{i}_{ref} = \frac{\mathbf{v}_g}{\norm{\mathbf{v}_g}_2^2}p_i^{N}+J\frac{\mathbf{v}_g}{\norm{\mathbf{v}_g}_2^2}q_i^{N},
\end{equation}
where $\mathbf{i}_{ref} \in \mathbb{R}^2$. Let $p_i^{N}, q_i^{N}$ be the exact solutions to the unconstrained optimization problems in (\ref{Network_Opt}a,\ref{Network_Opt}c), then 
\begin{equation}\label{closed_current_ref}
  \mathbf{i}_{ref} = -\frac{1}{\beta_i(\tilde{r}_{g_i})}\bigg(\frac{\mathbf{v}_g}{\norm{\mathbf{v}_g}_2^2}\lambda_i^N + J \frac{\mathbf{v}_g}{\norm{\mathbf{v}_g}_2^2}\nu_i^N\bigg),   
\end{equation}
from which it is seen that that the reference $\mathbf{i}_{ref}$ to the low-level controller decreases monotonically with $\beta_i(\tilde{r}_{g_i})$. Figure \ref{Control_Architecture} shows the overall control architecture.

\subsection{Numerical Simulation}\label{Sim}
\begin{figure}[t!] 
\centerline{\includegraphics[width=0.6\textwidth]{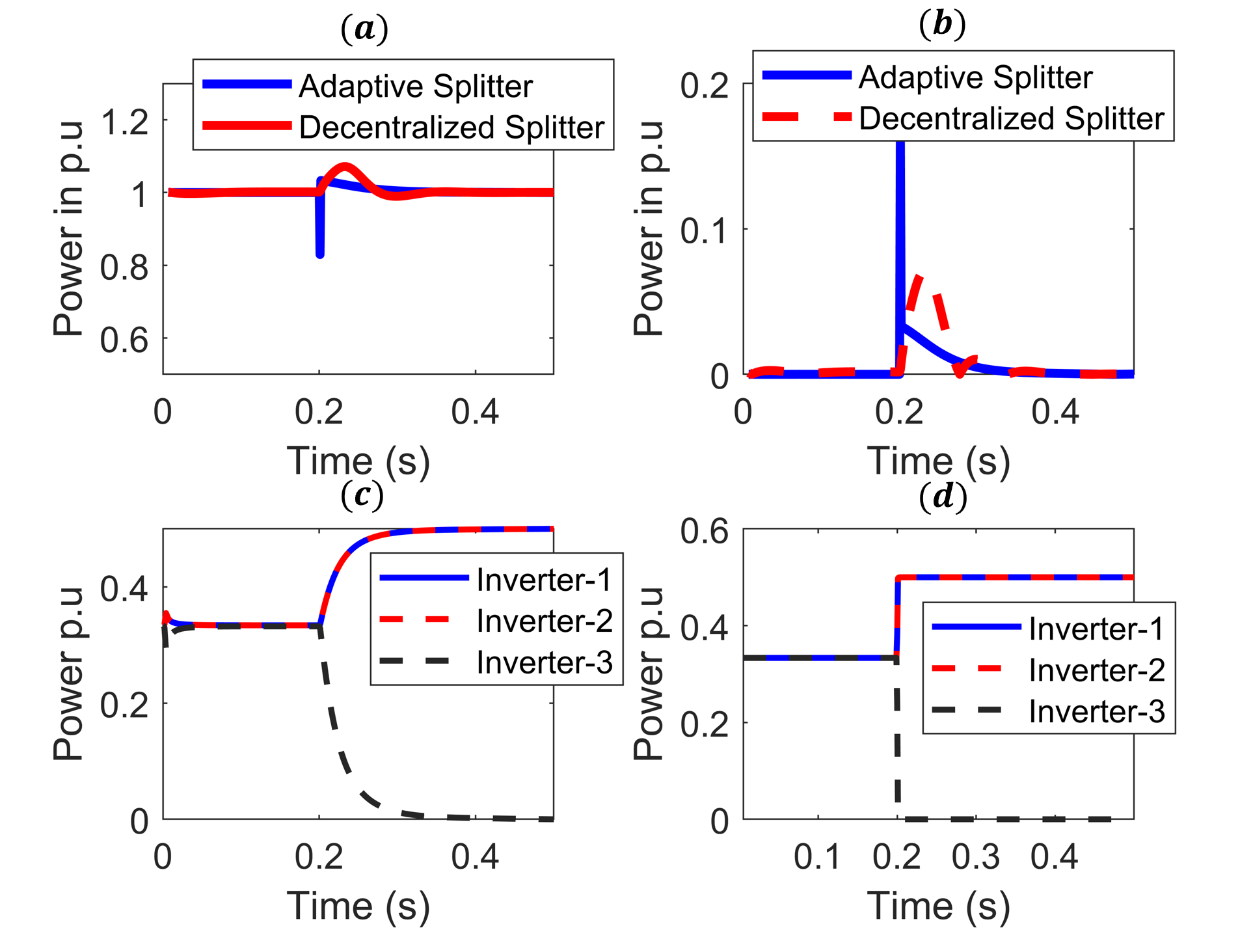}}
\caption{Active power tracking and tracking error (top) and power split (reference generated) for the adaptive splitter-based control design (left) and decentralized control design (right).}
\label{fig_example2.jpg}
\end{figure}

The simulation is performed on a system with 3-IBRs. We consider a scenario where there is a fault in Inverter-3 at $t=0.2s$. Also, a grid voltage swell (as disturbance) of $10\%$ the nominal voltage is induced at $t=0.2s$. The system parameters used in the simulation are the transmission line resistance $r_g=0.027\Omega$, inductance $l_g=0.0367$\textsf{H}, the nominal grid voltage $v_g = 392$\textsf{V} (LN rms), the grid frequency $w_g=60$\textsf{Hz}. The optimization parameters chosen are $\alpha = 0.1$, the optimization weights $\beta_i=1$ under normal operation, and $\beta_i=10^4$ under faulty conditions. The adjacency matrix, the degree matrix, and the Laplacian based weight matrix are determined as follows: $$A=\begin{bmatrix} 0 & 1 &0\\1 & 0 & 1\\0 &1 &0 \end{bmatrix}, D=\begin{bmatrix} 1 & 0 &0\\0 & 2 & 0\\0 &0 &1 \end{bmatrix}, W=\begin{bmatrix} 0.666 & 0.333 & 0\\0.333 & 0.333 & 0.333\\0 &0.333 &0.666 \end{bmatrix}$$
it can be seen that the matrix $W$ is doubly stochastic. The results of the designed \emph{decentralized splitter} controller are compared with the \emph{adaptive splitter} based controller proposed in \cite{AMELI}. The bottom half of Figure \ref{fig_example2.jpg} (c), (d) shows the active power distribution in the adaptive splitter-based and decentralized control-based design. It can be seen that when there is a fault, the split mechanism in the adaptive splitter-based design slowly adjusts to balance the aggregated active power. This is because the low-level control in this approach is a function of the high-level distribution. Whereas, in the decentralized-based approach, it can be seen that the distribution almost instantaneously adjusts to balance the active power references to the low-level control, this is because of the changes in the optimization penalty weight $\beta_i(\tilde{r}_{g_i})$. Moreover, the low-level control designed asymptotically tracks the high-level generated power reference.

\begin{figure}[t!] 
\centerline{\includegraphics[width=0.6\textwidth]{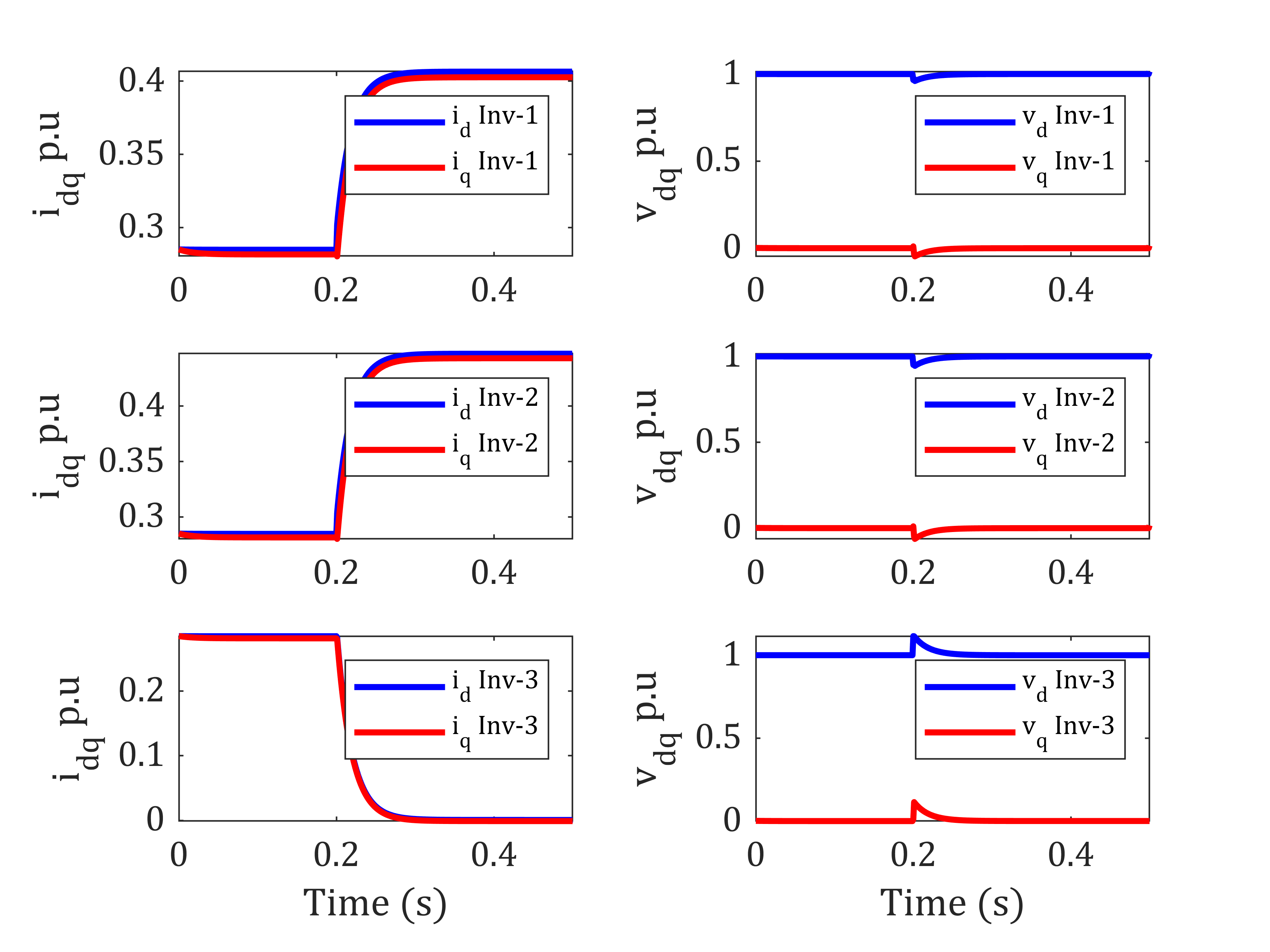}}
\caption{Currents $\mathbf{i}$ and voltages $\mathbf{v}$ of the individual inverters before and after the fault}
\label{fig_example.jpg}
\end{figure}

\begin{figure}
    \centering
    \includegraphics[width=0.6\textwidth]{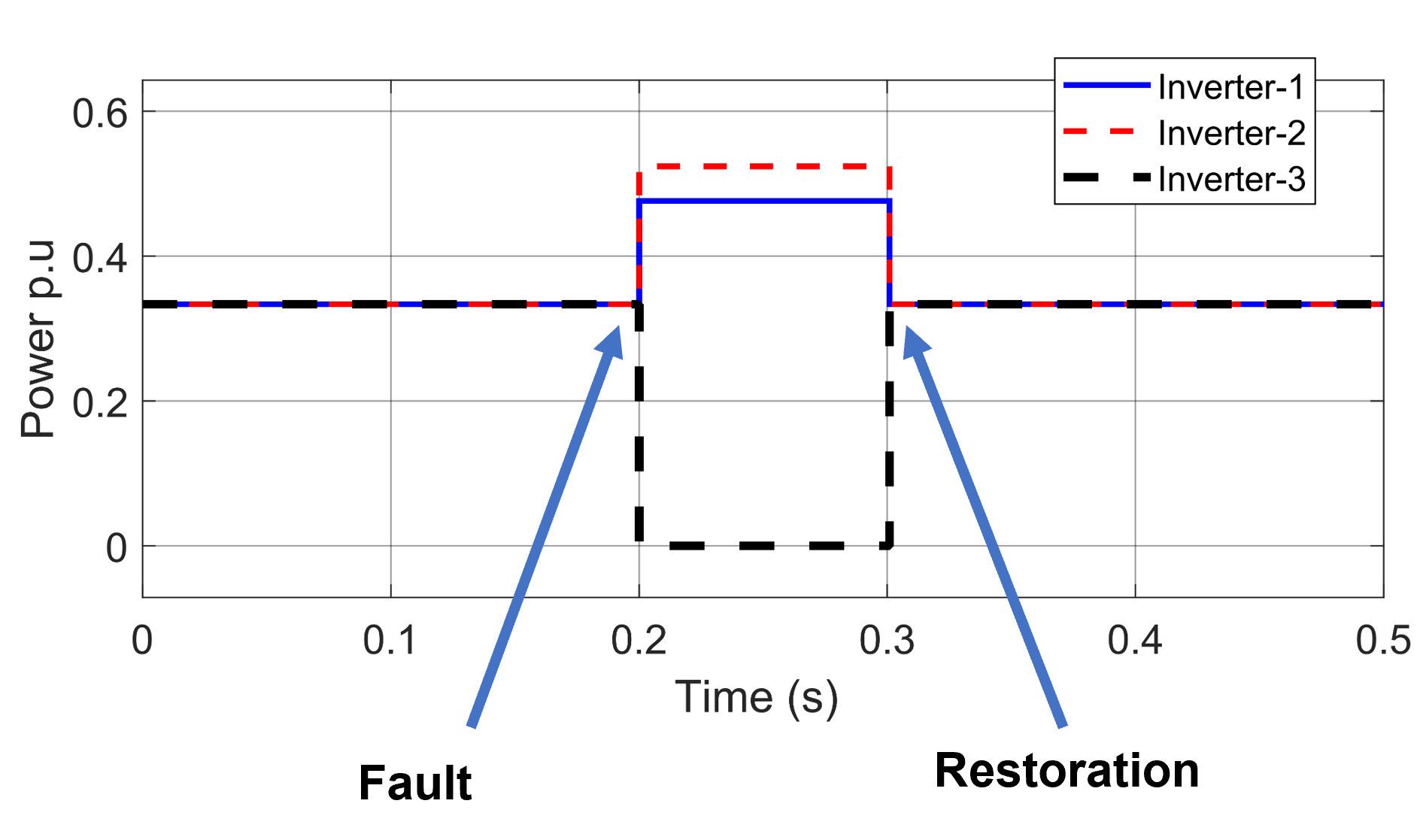}
    \caption{Power split when the fault occurs at 0.2s and restoration at 0.3s.}
    \label{fig:fr_3ibr}
\end{figure}

The aggregated active power tracking shown in the top half of Figure-\ref{fig_example2.jpg} (a), (b) compares the high-level active power distribution between the adaptive control-based splitter design and the decentralized optimization-based splitter design. It can be seen that when the fault and the grid voltage swell occur at $t=0.2s$, the decentralized optimization-based controller provides a better tracking response. The magnitude of active power deviation from the nominal $1$p.u value is significant while using the adaptive control-based power distribution mechanism. The result shows a significant improvement in the context of regulating active power at nominal value from the existing design technique in the literature. 

\begin{figure}
    \centering
    \includegraphics[width=0.6\textwidth]{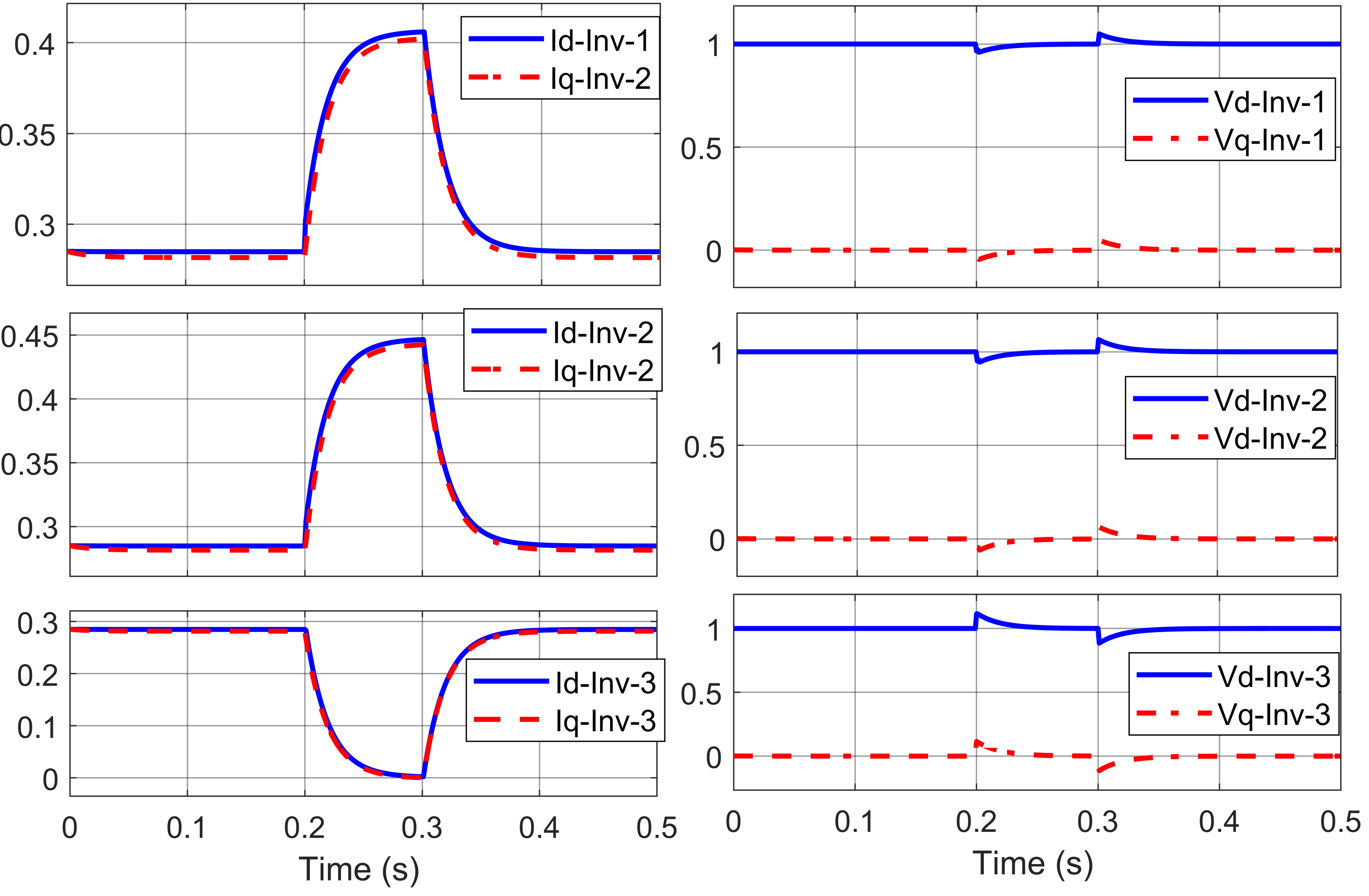}
    \caption{Currents $\mathbf{i}$ and voltages $\mathbf{v}$ of the individual inverters during fault and after restoration.}
    \label{fig:vi_fr_3ibr}
\end{figure}

Finally, the currents ($\mathbf{i}$) and the input voltages ($\mathbf{v}$) for the three inverters are shown in Figure \ref{fig_example.jpg}. It can be seen from the subplot (e) that the d and q-axis currents of Inverter-3 go to zero after the fault. Similarly, the impact of the fault can also be seen on the d and q-axis voltages of the Inverter-3 from subplot (f). Moreover, the current sharing between the IBRs is split between the remaining two IBRs with healthy transmission lines which can be seen from subplots (a) and (c).

Thus, the results demonstrate the impact of the designed high-level decentralized optimization-based control. The comparison with the adaptive control-based power splitter further solidifies the control effectiveness of the designed high-level controller.

\subsection{Conclusion}
This section presents a decentralized high-level control for the aggregated active power tracking in the IBRs under complete loss of inverter due to the faults in the transmission line. The proposed algorithm adaptively splits the active power taking into consideration the condition of the IBR lines. A parameter estimator is developed which determines the healthiness of the transmission line and updates the high-level algorithm on the transmission line's status. The results show an improvement in the active power tracking response compared to the previously proposed method. The developed decentralized algorithm relies on the established concept of duality. The future extensions of this work will focus on the stability of the developed decentralized control. Investigating the transmission limitations and their impact on power-sharing are reported in the next section.

\bibliographystyle{IEEEtran}
\bibliography{myrefs}

\end{document}